\documentclass[12pt,a4paper]{article}
\usepackage{amscd,amsmath,amssymb,amsfonts,times,amsbsy}
\usepackage{mathrsfs}
\usepackage[dvips]{color} 
\usepackage[all]{xy}
\usepackage[T1]{fontenc}
\usepackage{textcomp}
\usepackage{tikz-cd}

\numberwithin{equation}{section}
    \newtheorem{thm}{Theorem}[section]
    \newtheorem{lem}[thm]{Lemma}
    
    \newtheorem{prop}[thm]{Proposition}
    \newtheorem{cor}[thm]{Corollary}

    \newtheorem{conj}[thm]{Conjecture}
    \newtheorem{defn}[thm]{Definition}
    
    \newtheorem{rem}[thm]{Remark}
\DeclareMathAlphabet{\mathpzc}{OT1}{pzc}{m}{it}
\newcommand{\qed}{\mbox{}\nolinebreak$\square$\medbreak\par}
\newenvironment{pf}{\par\smallskip\noindent\emph{Proof.}}{\hfill\qed\par\smallskip}
\newenvironment{pf*}[1]{\par\smallskip\noindent\emph{#1.}}{\hfill\qed\par\smallskip}
\newcommand{\bysame}{\hskip.3em \leavevmode\rule[.5ex]{2.5em}{.3pt}\hskip0.5em}
\begin{document}
\title{A numerical approach toward the $p$-adic Beilinson conjecture for elliptic curves over $\mathbb Q$}
\author{Masanori Asakura,  
 Masataka Chida}
\date\empty
\maketitle
\vspace{-1.3cm}
\begin{center}
{\large with Appendix B by Fran\c{c}ois Brunault}
\end{center}

\def\can{\mathrm{can}}
\def\cano{a}
\def\ch{{\mathrm{ch}}}
\def\Coker{\mathrm{Coker}}
\def\crys{\mathrm{crys}}
\def\dlog{{\mathrm{dlog}}}
\def\dR{{\mathrm{d\hspace{-0.2pt}R}}}            
\def\et{{\mathrm{\acute{e}t}}}  
\def\Frac{{\mathrm{Frac}}}
\def\phami{\phantom{-}}
\def\id{{\mathrm{id}}}              
\def\Image{{\mathrm{Im}}}        
\def\Hom{{\mathrm{Hom}}}  
\def\Ext{{\mathrm{Ext}}}
\def\MHS{{\mathrm{MHS}}}  
\def\GL{{\mathrm{GL}}}  
\def\SL{{\mathrm{SL}}}  
  
\def\ker{{\mathrm{Ker}}}          
\def\rig{{\mathrm{rig}}}
\def\Pic{{\mathrm{Pic}}}
\def\CH{{\mathrm{CH}}}
\def\NS{{\mathrm{NS}}}
\def\an{{\mathrm{an}}}
\def\End{{\mathrm{End}}}
\def\Tr{{\mathrm{Tr}}}
\def\Gal{{\mathrm{Gal}}}
\def\Proj{{\mathrm{Proj}}}
\def\ord{{\mathrm{ord}}}
\def\reg{{\mathrm{reg}}}          %
\def\res{{\mathrm{res}}}          %
\def\Res{\mathrm{Res}}
\def\Spec{{\mathrm{Spec}}}     
\def\syn{{\mathrm{syn}}}
\def\cont{{\mathrm{cont}}}
\def\zar{{\mathrm{zar}}}
\def\bA{{\mathbb A}}
\def\bC{{\mathbb C}}
\def\C{{\mathbb C}}
\def\G{{\mathbb G}}
\def\bE{{\mathbb E}}
\def\bF{{\mathbb F}}
\def\F{{\mathbb F}}
\def\bG{{\mathbb G}}
\def\bH{{\mathbb H}}
\def\bJ{{\mathbb J}}
\def\bL{{\mathbb L}}
\def\cL{{\mathscr L}}
\def\bN{{\mathbb N}}
\def\bP{{\mathbb P}}
\def\P{{\mathbb P}}
\def\bQ{{\mathbb Q}}
\def\Q{{\mathbb Q}}
\def\bR{{\mathbb R}}
\def\R{{\mathbb R}}
\def\bZ{{\mathbb Z}}
\def\Z{{\mathbb Z}}
\def\cH{{\mathscr H}}
\def\cD{{\mathscr D}}
\def\cE{{\mathscr E}}
\def\cF{{\mathscr F}}
\def\O{{\mathscr O}}
\def\cR{{\mathscr R}}
\def\cS{{\mathscr S}}
\def\cX{{\mathscr X}}
\def\cM{{\mathscr M}}
%
\def\ve{\varepsilon}
\def\vG{\varGamma}
\def\vg{\varGamma}
%
%

%
\def\lra{\longrightarrow}
\def\lla{\longleftarrow}
\def\Lra{\Longrightarrow}
\def\hra{\hookrightarrow}
\def\lmt{\longmapsto}
\def\ot{\otimes}
\def\op{\oplus}
\def\Isoc{{\mathrm{Isoc}}}
\def\Fil{{\mathrm{Fil}}}
\def\FIsoc{{F\text{-Isoc}}}
\def\FMIC{{F\text{-MIC}}}
\def\Log{{\mathscr{L}{og}}}
\def\FilFMIC{{\mathrm{Fil}\text{-}F\text{-}\mathrm{MIC}}}

\def\wt#1{\widetilde{#1}}
\def\wh#1{\widehat{#1}}
\def\spt{\sptilde}
\def\ol#1{\overline{#1}}
\def\ul#1{\underline{#1}}
\def\us#1#2{\underset{#1}{#2}}
\def\os#1#2{\overset{#1}{#2}}

\def\s{{(\sigma)}}
\def\EE{2g_2g'_3-3g'_2g_3}
\def\EEE{6g_2g'_3-9g'_2g_3}
\newcommand{\Liminj}[1]{\raisebox{-1.98ex}%
{$\stackrel{\normalsize\mbox{$\varinjlim$}}{\scriptstyle #1}$}}
\newcommand{\Limproj}[1]{\raisebox{-2.00ex}%
{$\stackrel{\normalsize\mbox{$\varprojlim$}}{\scriptstyle #1}$}}

\begin{abstract}
Restricting ourselves to elliptic curves over $\Q$, we reformulate 
the $p$-adic Beilinson conjecture due to Perrin-Riou, which is customized to our computational approach. 
We then develop a new algorithm for numerical verifications of the $p$-adic Beilinson conjecture, which is
based on
the theory of rigid cohomology and $F$-isocrystals. 
\end{abstract}

\section{Introduction}
Since Dirichlet proved his famous class number formula in the 19th century,
the regulator is one of the most important subjects in algebraic number
theory.
In the 20th century, several generalizations of Dirichlet's formula 
were proposed in the context of higher $K$-theory by
Lichtenbaum, Borel, etc.
Especially, Beilinson's conjecture
on special values of motivic $L$-functions
(\cite{beilinson}, \cite{schneider}) is an important milestone, and this mysterious conjecture
is still attracting lots of researchers in arithmetic geometry.

\medskip

The $p$-adic analogue of the Beilinson conjecture
was formulated by Perrin-Riou \cite{PR}, which we call the {\it $p$-adic Beilinson conjecture}.
Roughly speaking, it asserts a certain formula on
the $p$-adic regulators and 
the special values of $p$-adic $L$-functions, which is conceptionally described as
\[
\frac{\text{$p$-adic $L$-value}}{\text{$p$-adic regulator}}
=\frac{\text{$L$-value}}{\text{Beilinson regulator}}.
\]
However, the above involves a serious issue. In the present,
there is no general definition or even candidate
of the $p$-adic counterparts of motivic $L$-functions.
She needed to formulate the statement together with the existence of $p$-adic 
$L$-functions
simultaneously, and that causes the very complicated statement.

On the other hand,
restricting ourselves to the case of elliptic curves over $\Q$, 
the $p$-adic $L$-functions are defined by the following people.
\begin{itemize}
\item Mazur and Swinnerton-Dyer \cite{MS}
\item Vishik \cite{V}, Amice and V\'elu \cite{AV}
\item Pollack and Stevens \cite{PS2}, Bella\"iche \cite{B}
(critical slope $p$-adic $L$-functions).
\end{itemize}
Thus one can re-formulate the $p$-adic Beilinson conjecture in terms of the above
$p$-adic $L$-functions (cf. \cite[4.2]{PR}).
However there still remains a problem.
In general, the regulators are
only determined up to multiplication by $\Q^\times$,
in other words there is no canonical manner to remove the ambiguity of $\Q^\times$.
The $p$-adic Beilinson conjecture is stated by fixing 
basis of several cohomology groups, though
the choice does not matter by replacing the normalization 
of her $p$-adic $L$-functions if necessary.
On the other hand, the normalization of the $p$-adic $L$-functions by
Mazur and Swinnerton-Dyer et al is
fixed and it is a natural one.
Therefore we need to choose a suitable basis of cohomology groups according to this,
and it is not so an obvious question.

\medskip

In this paper
we write down the $p$-adic Beilinson conjecture
for elliptic curves over $\Q$ in terms of the 
$p$-adic $L$-functions by Mazur and Swinnerton-Dyer et al. To be precise, 
for an elliptic curve $E$ over $\Q$ and $\xi\in H^2_\cM(E,\Q(n+2))$,
we introduce $R_{\infty}(h^1(E)(-n), \xi)$ the {\it Beilinson regulator 
without ambiguity of $\Q^\times$}, and
$R_{p,\gamma}(h^1(E)(-n),\xi)$ the {\it $p$-adic regulator without ambiguity of $\Q^\times$}. Then our formulation
of the $p$-adic Beilinson conjecture is described as follows (Conjecture \ref{PREC}):
\begin{equation}\label{int-eq1}
\frac{L_{p,\gamma}(E,\omega^{-n-1},-n)}{R_{p,\gamma}(h^1(E)(-n),\xi)}
=\pm \frac{L'(E,-n)}{R_{\infty}(h^1(E)(-n), \xi)},\quad n\geq 0.
\end{equation}
We also give a formulation for elliptic modular forms of weight $k\geq 2$
with rational Fourier coefficients
(Conjecture \ref{PRMF}).
Our statement features that 
it is more accessible than the original one, especially concerning numerical
verifications.
The latter half of this paper is devoted to give 
numerical verifications of our statement in case $n=0$:
\begin{equation}\label{int-e2}
\frac{L_{p,\gamma}(E,\omega^{-1},0)}{R_{p,\gamma}(h^1(E),\xi)}
=\pm \frac{L'(E,0)}{R_{\infty}(h^1(E), \xi)}.
\end{equation}
To do this, the central role is played by the theory of {\it rigid cohomology} 
and {\it $F$-isocrystals}.
In particular, we use the deformation method by Lauder \cite{Lauder}
and apply the result on effective bounds on overconvergence by Kedlaya-Tuitman
\cite[Theorem 2.1]{KT}.
We shall show,
\begin{thm}\label{intro-thm1}
Let $f:X\to S\subset \P^1_\Q$ be an elliptic fibration which has at least one 
split multiplicative fiber.
Let $X_a=f^{-1}(a)$ be the fiber at a $\Q$-rational point $a\in S(\Q)$.
Let $\xi\in H^2_\cM(X,\Q(2))$. 
Then, under some mild assumption on $p$,
there is an algorithm for computing the $p$-adic regulator
\begin{equation}\label{int-eq3}
R_{p,\gamma}(h^1(X_a),\xi|_{X_a})
\end{equation}
modulo $p^n$ for each $n\geq 1$.
See \S \ref{CSyn-sect} for the down-to-earth algorithm.
\end{thm}
\noindent
For a particular example, there is a theoretical proof of \eqref{int-e2}.

\begin{thm}[Theorem \ref{X4-thm}]\label{intro-thm2}
Let $E$ be the elliptic curve defined by an equation $y^2=x(1-x)(1+3x)$ over $\Q$.
Let
\[
\xi=\left\{\frac{y+1-x}{y-1+x},\frac{4x^2}{(1-x)^2)}\right\}\in K_2(E).
\]
Let $p$ be a prime at which $E$ has a good reduction.
Let $\gamma$ be a root of eigenpolynomial of $E$ modulo $p$ such that $\ord_p(\gamma)<1$.
Then \eqref{int-e2} holds.
\end{thm}
The proof uses a formula of Rogers-Zudilin for computing the right hand side
of \eqref{int-e2}, and
Brunault's formula for computing the left hand side,
together with comparing the symbol $\xi$ with the Beilinson-Kato element (Appendix B).

\medskip

There have been several papers concerning the $p$-adic Beilinson conjecture,
Bannai-Kings \cite{BK}, Kings-Loeffler-Zerbes \cite{KLZ}, Niklas \cite{Niklas}, 
Brunault \cite{Br}, Bertolini-Darmon \cite{BD} and probably more.
A noteworthy point is that our approach is entirely different from theirs. 
For example, many papers rely on Kato's seminal paper \cite{K}
where the Beilinson-Kato elements are introduced.
Ours is less powerful than the method using \cite{K}
while ours allows an approach to new examples such as
the Beilinson-Kato elements vanish in the $p$-adic cohomology
groups $H^1_f(G_\Q,H^1_\et(\ol E,\Q_p(2))$ 
(Remark \ref{remark-BK}).
This also allows us to be convinced of our reformulation of 
the $p$-adic Beilinson conjecture.

\bigskip

The paper is organized as follows.
\S \ref{pLsect} is a review of $p$-adic $L$-functions for
elliptic modular forms, which includes the recent development of
critical slope $L$-functions.
We give a precise statement of the $p$-adic Beilinson conjecture for elliptic curves
in \S \ref{PRsect} (Conjectures \ref{PRMF}, \ref{PREC}).
In \S \ref{Tate-sect}, we provide the algorithm for computing 
the $p$-adic regulators \eqref{int-eq3}, which is summed up in \S \ref{CSyn-sect}.
With use of this algorithm,
we give numerical verifications of the $p$-adic Beilinson conjecture for explicit 
examples of elliptic curves in \S \ref{NVPR-sect} where we 
use SAGE or MAGMA for computing the special values of ($p$-adic) $L$-functions.
The values of critical slope $p$-adic $L$-functions were kindly provided by Professor 
Robert Pollack.
In \S \ref{Appendix-sect} Appendix A, we sum up the standard results on the Gauss-Manin
connection for elliptic fibrations, which are used especially in \S \ref{Tate-sect}.
Most results in Appendix A are probably known to experts.
However since the authors do not find a suitable literature, 
we give a comprehensive exposition
for the convenience of the reader.
Finally we compare the symbol in Theorem \ref{intro-thm2}
with the Beilinson-Kato element in
\S \ref{AppendixB-sect} Appendix B.

\bigskip

\noindent{\bf Acknowledgement.}
The authors very much appreciate Professor Fran\c{c}ois Brunault for kindly providing Appendix B. 
The authors would also express sincere gratitude to Professor Robert Pollack for providing 
numerical values of the critical slope $p$-adic $L$-functions, and 
Professor Shinichi Kobayashi for helpful suggestions and comments.

\section{$p$-adic $L$-functions for elliptic modular forms}\label{pLsect}

Here we recall some basic facts on the cyclotomic $p$-adic $L$-functions for elliptic modular forms.
Let $N$ be a positive integer greater than $3$.
Fix a prime $p>2$ not dividing $N$ and fix embeddings $\overline{\mathbb{Q}} \hookrightarrow \mathbb{C}$ and $\overline{\Q} \hookrightarrow \overline{\Q}_p$.
Let $K$ be a subfield of $\overline{\mathbb{Q}}$.
Let $S_{k+2}(\Gamma_1(N), K)$ be the space of cuspforms of weight $k+2$ for $\Gamma_1(N)$ whose Fourier coefficients belong to $K$.
Let $f$ be a normalized eigenform of weight $k+2$ for $\Gamma_0(N)$ with character $\chi_f$.
We denote the $n$-th Fourier coefficient of $f$ by $a_n(f)$.
Let $\alpha$ and $\beta$ be the two roots of $X^2-a_p(f)X+\chi_f(p)p^{k+1}$ with $\operatorname{ord}_p(\alpha )\leq \operatorname{ord}_p(\beta)$.
Then
$f_{\alpha}(z)=f(z)-\beta f(pz)$ and $f_{\beta}(z)=f(z)-\alpha f(pz)$ are normalized cuspforms for $\Gamma_0=\Gamma_1(N)\cap \Gamma_0(p)$
of weight $k+2$.
Moreover they are also eigenforms for $U_p$-operator and satisfy $U_pf_{\alpha}(z)=\alpha f_{\alpha}(z)$ and $U_pf_{\beta}(z)=\beta f_{\beta}(z)$.
$f_{\alpha}$ and $f_{\beta}$ are called $p$-refinements (or $p$-stabilizations) of $f$.
If $ \operatorname{ord}_p(\beta)=k+1$, the $p$-refinement $f_\beta$ is called a cuspform of critical slope.
If $f$ is non-ordinary, then $\alpha$ and $\beta$ satisfy $\operatorname{ord}_p(\alpha )\leq \operatorname{ord}_p(\beta)< k+1$.
In this case, we have $p$-adic $L$-functions $L_{p,\alpha}(f,\chi,s)$ and $L_{p,\beta}(f,\chi ,s)$ constructed by Vishik \cite{V} and Amice-V\'elu \cite{AV}
for any Dirichlet characters $\chi$.
If $f$ is ordinary, $\alpha$ satisfies $\operatorname{ord}_p(\alpha )=0$. Therefore we also have $p$-adic $L$-function $L_{p,\alpha}(f,\chi,s)$.
For the critical slope case, Pollack-Stevens \cite{PS1, PS2} and Bella\"iche \cite{B} constructed $p$-adic $L$-functions $L_{p,\beta}(f,\chi ,s)$ under certain assumptions.
In this section, we recall the constructions briefly.

\subsection{Preliminaries}

Let $\mathcal{H}=\mathbb{Z}[T_{\ell} \, (\ell \nmid N), \, U_q \, (\ell | N)]$ be the Hecke algebra.
For a commutative ring $R$, let $\mathscr{P}_k (R)$ denote the space of homogeneous polynomials of degree $k$ in $X$ and $Y$ with $R$-coefficients.
Let $\Delta_0=\operatorname{Div}^0(\mathbb{P}^1(\mathbb{Q}))$ be the degree zero divisors on $\mathbb{P}^1(\mathbb{Q})$.
Then $\Delta_0$ is a left $\mathrm{GL}_2(\mathbb{Q})$-module
by the linear fractional transformation.


Let $Y_1(N)$ be the open modular curve over $\mathbb{Q}$ with level $\Gamma_1(N)$-structure and
$\pi : \mathscr{E}\to Y_1(N)$ the universal elliptic curve.
The line bundle $\underline{\omega}=\pi_*\Omega_{\mathscr{E}\slash Y_1(N)}^1$ is naturally extended to the modular curve $X_1(N)$.
Then we have an identification
$$
S_{k+2}(\Gamma_1(N),K)= H^0(X_1(N)_K, \underline{\omega}^k\otimes \Omega_{X_1(N)}^1).
$$
For the local system $\mathscr{L}_k^{\mathrm{B}}=\operatorname{Sym}^k R^1\pi_* \mathbb{Z}$, we also have a canonical identification
$$
H^1_c(Y_1(N)(\mathbb{C}),\mathscr{L}_k^\mathrm{B})
\cong \operatorname{Hom}_{\Gamma_1(N)}(\Delta_0,\mathscr{P}_k(\mathbb{Z})).
$$
For simplicity, we denote $H^1_c(\Gamma_1(N),\mathscr{P}_k(R))= H^1_c(Y_1(N)(\mathbb{C}),\mathscr{L}_k^{\mathrm{B}})\otimes R$ for a commutative ring $R$.
The element $\iota = \begin{pmatrix} -1 & 0 \\ 0 & 1\end{pmatrix}$ defines an involution. If $R$ is a $\mathbb{Z} [1\slash 2]$-module, then we have a decomposition
$$H_c^1(\Gamma_1(N),\mathscr{P}_k(R))=H_c^1(\Gamma_1(N),\mathscr{P}_k(R))^{+} \oplus H_c^1(\Gamma_1(N),\mathscr{P}_k(R))^{-},$$
where $H_c^1(\Gamma_1(N),\mathscr{P}_k(R))^{\pm} =\left\{ x \in H_c^1(\Gamma_1(N),\mathscr{P}_k(R)) \middle| \iota (x) =\pm x \right\}$.
It is known that the action of $\iota$ coincides with the action of the complex conjugation (see \cite[Proof of Lemma 7.19]{K}).
For a cuspform $g \in S_{k+2}(\Gamma_1(N), \mathbb{C})$, we define a cohomology class
$\xi_g \in H_c^1(\Gamma_1(N),\mathscr{P}_k(\mathbb{C}))\cong \operatorname{Hom}_{\Gamma_1(N)}(\Delta_0,\mathscr{P}_k(\mathbb{C}))$ by
$$
\xi_g (\{ x \} - \{ y\} )=2\pi i \int_x^y g(z)(zX+Y)^kdz.
$$
Let $K$ be a subfield of $\overline{\mathbb{Q}}_p$ and $\phi : \mathcal{H}\to K$ be an algebra homomorphim and $V$ a finite dimensional vector space over $K$
with an action of $\mathcal{H}$.
For $\alpha \in K$ and $T\in \mathcal{H}$, let $V[\alpha, T]$ be the eigenspace for $T$ acting on $V$ with the
eigenvalue $\alpha$ and we put
$$
V{[\phi ]}=\bigcap_{T\in \mathcal{H}}V[(\phi (T),T)].
$$
Smilarly, let $V_{(\alpha ,T)}$ be the generalized eigenspace for $T$ acting on $V$ with the
eigenvalue $\alpha$ and we define the $\phi$-isotypic subspace $V_{(\phi)}$ of $V$ by
$$
V_{(\phi)}=\bigcap_{T\in \mathcal{H}}V_{(\phi (T),T)}.
$$
For a normalized eigenform $f(z)$ in $S_{k+2}(\Gamma_1(N),\overline{\mathbb{Q}})$, we put $L=\mathbb{Q}(\{ a_n(f)\}_{n=1}^\infty)$.
It is well-known that $L$ is a finite extension of $\mathbb{Q}$.
Then we can define the homomorphism $\phi_f : \mathcal{H} \to L$
by $\phi_f (T_{\ell})=a_{\ell}(f)$ for $\ell \nmid N$ and $\phi_f (T_q)=a_q (f)$ for $q|N$.
For simplicity, we denote $V[f]=V[\phi_f]$ and $V_{(f)}=V_{(\phi_f)}$.
\subsection{Result of Amice-V\'elu and Vishik}
Let $f$ be a normalized newform in $S_{k+2}(\Gamma_1(N),\overline{\mathbb{Q}})$ and we denote the complex conjugate cuspform
$\sum_{n=1}^\infty \overline{a_n(f)}q^n$ of $f$ by $f^*$.
Let $v$ be the valuation induced by the embedding $L\hookrightarrow \overline{\mathbb{Q}}\hookrightarrow \overline{\mathbb{Q}}_p$
and $\mathscr{O}_{L,v}$ the discrete valuation ring defined by $v$. We define
$$
H_c^1(\Gamma_1(N),\mathscr{P}_k(\mathscr{O}_{L,v}))^{\pm}[f]
=H_c^1(\Gamma_1(N),\mathscr{P}_k(L))^{\pm}[f]\cap H_c^1(\Gamma_1(N),\mathscr{P}_k(\mathscr{O}_{L,v})).
$$
Then it is known that $H_c^1(\Gamma_1(N),\mathscr{P}_k(\mathscr{O}_{L,v}))^{\pm}[f]$ is a free $\mathscr{O}_{L,v}$-module
of rank one.
We fix a generator $u^\pm$ in $H_c^1(\Gamma_1(N),\mathscr{P}_k(\mathscr{O}_{L,v}))^{\pm}[f^*]$.
For $0\neq \eta \in S_{k+2}(\Gamma_1(N),L)[f^*]$ and $u=(u^+,u^-)$, we define a complex number $\Omega_f^\pm=\Omega_f^\pm (\eta ,u)$
by
$$
\xi_\eta^\pm =\Omega_f^\pm \cdot u^\pm .
$$
\begin{rem}\label{normalization}
\begin{enumerate}
\item By the multiplicity one, $\eta$ is a nonzero multiple of $f$ by an element in $L^\times$.
\item The generator $u^\pm$ is determined only up to an element in $\mathscr{O}_{L,v}^\times$.
Since
$$
H_c^1(\Gamma_1(N),\mathscr{P}_k(\mathscr{O}_{L}))^{\pm}[f^*]
=H_c^1(\Gamma_1(N),\mathscr{P}_k(L))^{\pm}[f]\cap H_c^1(\Gamma_1(N),\mathscr{P}_k(\mathscr{O}_{L}))
$$
is not a free $\mathscr{O}_{L}$-module in general,
we do not have a canonical choice of $u^\pm$ to define the periods.
Therefore, our $p$-adic $L$-functions will depend on the choice of $u^\pm$.
\end{enumerate}
\end{rem}

Let $f$ be a normalized eigenform of weight $k+2$ for $\Gamma_0(N)$ with character $\chi_f$.
Let $\alpha$ be a root of $X^2-a_p(f)X+\chi_f (p) p^{k+1}$.
For $a,m \in \mathbb{Q}$ and $P \in \mathscr{P}_k (\mathbb{C})$, we put
$$
\lambda_{f,\alpha}^{\pm}(a,m)(P)=\frac{\pi i}{\Omega_f^{\pm}}\left\{ \int_\infty^{a/m}f_{\alpha}(z)P(-mz+a,1)dz \pm \int_{\infty}^{-a/m}f_{\alpha}(z)P(mz+a,1)dz \right\} .
$$
\begin{thm}(Vishik \cite{V}, Amice-V\'elu \cite{AV})
Let $f$ be a normalized eigenform of weight $k+2$ for $\Gamma_0(N)$ with character $\chi_f$.
Let $\alpha$ be a root of $X^2-a_p(f)X+\chi_f (p) p^{k+1}$ satisfying $\mathrm{ord}_p(\alpha)<k+1$.
Then there exists a unique $p$-adic distribution $\mu^{\pm}_{f,\alpha}$ on $\mathbb{Z}_p^\times$ satisfying
the follwing properties:
\begin{enumerate}
\item
$\displaystyle \int_{a+p^\nu\mathbb{Z}_p}x^j d\mu^\pm_{f,\alpha}=\lambda^\pm_{f,\alpha}(a,p^\nu)(X^jY^{k-j})$ for $\nu >0$ and $0\leq j \leq k$.
\item
For $0\leq i < \operatorname{ord}_p(\alpha)$, $\displaystyle \sup_{a\in \mathbb{Z}_p^\times} \displaystyle \left| \int_{a+p^\nu\mathbb{Z}_p}(x-a)^i d\mu^\pm_{f,\alpha} \right|_p=
o(p^{\nu (\operatorname{ord}_p(\alpha)-i)})$ when $\nu$ tends to infinity
(i.e. $\mu^\pm_{f,\alpha}$ is $\operatorname{ord}_p(\alpha)$-admissible).
\end{enumerate}
\end{thm}
For $x\in \mathbb{Z}_p^\times$, we write $x=\omega(x)\langle x \rangle$ with $\omega (x)\in \mu_{p-1}(\mathbb{Z}_p)$ and $\langle x\rangle \in 1+p\mathbb{Z}_p$.
Let $\chi : \mathbb{Z}_p^\times \to \overline{\mathbb{Q}}^\times$ be a character of conductor $p^n$.
Now we set $\mu_{f,\alpha}=\mu_{f,\alpha}^+ +\mu_{f,\alpha}^-$ and
define the $p$-adic $L$-function $L_{p,\alpha}(f,\chi,s)=L_{p,\alpha,\eta , u}(f,\chi,s)$ by
$$L_{p,\alpha}(f,\chi,s)=\int_{\mathbb{Z}_p^\times}\chi (x) \langle x \rangle^{s-1} d\mu_{f,\alpha}$$
for $s\in \mathbb{Z}_p$. Then we have the following interpolation property.
\begin{thm}\label{interpolation}
Assume that the root $\alpha$ satisfies $\mathrm{ord}_p(\alpha)<k+1$.
For a primitive character $\chi : (\mathbb{Z}\slash p^\nu \mathbb{Z})^\times \to \overline{\mathbb{Q}}^\times$ and $1\leq j \leq k+1$,
we have
$$
L_{p,\alpha}(f,\omega^{j-1} \chi,j)=
\begin{cases}\displaystyle
(1-\chi_f (p) p^{k+1-j}\alpha^{-1}) (1-p^{j-1}\alpha^{-1})\frac{L(f,j)}{\Omega_f^{\epsilon ( \mathbf{1},j)}} & \textup{if } \chi = \mathbf{1},\\
\displaystyle
\left( \frac{p^{j\nu}}{\alpha^\nu}\right)\cdot \frac{(j-1)!}{(-2\pi i )^{j-1}\tau (\overline{\chi})}\cdot \frac{L(f,\overline{\chi},j)}{\Omega_f^{\epsilon (\chi ,j)}}& \textup{if }\chi \neq \mathbf{1},
\end{cases}
$$
where $\epsilon (\chi ,j)= (-1)^{(j-1)}\cdot \chi (-1)$ and $\tau$ is the Gauss sum defined by
$$\tau (\overline{\chi})=\sum_{a\in (\mathbb{Z}\slash p^\nu \mathbb{Z})^\times} \overline{\chi}(a)e^{2\pi i a \slash p^\nu}.$$
\end{thm}

\subsection{Critical slope $p$-adic $L$-functions}
If $\mathrm{ord}_p(\beta)=k+1$, it is said that $\beta$ has a critical slope. Here we briefly recall the construction of critical slope $p$-adic $L$-functions
following Pollack-Stevens \cite{PS1, PS2}.
Let $p$ be a prime and $N$ a positive integer satisfying $(p,N)=1$. We denote $\Gamma_0=\Gamma_1(N)\cap \Gamma_0(p)$
and $S_0(p)=\left\{  g= \begin{pmatrix}  a & b \\ c & d \end{pmatrix} \in M_2 (\mathbb{Z} ) \, \middle| \, p\nmid a, p\mid c, \operatorname{det}g  \neq 0 \right\}$.
For a right $\mathbb{Z}_p[S_0(p)]$-module $M$, we denote the locally constant sheaf associated to $M$ on the open modular curve $Y_{\Gamma_0}(\mathbb{C})=\mathbb{H}\slash \Gamma_0$ by $\widetilde{M}$,
where $\mathbb{H}$ is the upper half plane.
For simplicity, we assume that the order of the torsion elements of $\Gamma_0$ are prime to $p$ or $M$ is a $\mathbb{Z}_p$-module with $p>3$.
Recall that $\Delta_0=\operatorname{Div}^0(\mathbb{P}^1(\mathbb{Q}))$ is the degree zero divisors on $\mathbb{P}^1(\mathbb{Q})$.
For $\varphi \in \operatorname{Hom}(\Delta_0, M)$ and $g \in S_0(p)$, we put $(\varphi | g)(D)=\varphi(g D)|g$ for $D\in \Delta_0$.
Define Hecke operators $T_{\ell}$ $(\ell \nmid N)$ and $U_q$ $(q | N)$
by
$\varphi | T_{\ell}=\varphi | \begin{pmatrix}  \ell & 0 \\ 0 & 1 \end{pmatrix}+ \sum_{a=0}^{\ell-1}\varphi | \begin{pmatrix}  1 & a \\ 0 & \ell \end{pmatrix}$
and $\varphi | U_q=\sum_{a=0}^{q-1}\varphi | \begin{pmatrix}  1 & a \\ 0 & q \end{pmatrix}$.
We denote $H^1_c(\Gamma_0,M)=H^1_c(Y_{\Gamma_0}(\mathbb{C}),\widetilde{M})$ as before.
Then Ash-Stevens \cite[Proposition 4.6]{AS} proved that there exists a canonical Hecke equivariant isomorphism $H^1_c(\Gamma_0,M)\cong \operatorname{Hom}_{\Gamma_0}(\Delta_0,M)$,
where $\operatorname{Hom}_{\Gamma_0}(\Delta_0,M)=\{ \varphi \in \operatorname{Hom}(\Delta_0, M)  \mid  \varphi | g =\varphi \textup{ for all }g \in \Gamma_0 \}$ is the $\Gamma_0$-invariant homomorphisms.
From now on, we will identify these modules by the above isomorphism.
If $M$ is a Banach module and $\Gamma_0$ acts on $M$ as unitary operators, then $H^1_c(\Gamma_0, M)$ is also a Banach space with the norm defined by
$||\Psi ||=\sup_{D\in \Delta_0}|| \Psi (D)||$ for $\Psi \in H^1_c(\Gamma_0,M)$.

For $r\in |\mathbb{C}_p^\times |_p$, we denote
$$
B[\mathbb{Z}_p,r]=\{ z\in \mathbb{C}_p \, | \, \textup{there exists } a \in \mathbb{Z}_p \textup{ such that } |z-a|_p \leq r \}.
$$
Then $B[\mathbb{Z}_p,r]$ is the $\mathbb{C}_p$-valued points  of a $\mathbb{Q}_p$-affinoid variety.
Let $\mathbb{A}[r]$ be the $\mathbb{Q}_p$-Banach algebra of $\mathbb{Q}_p$-affinoid functions on $B[\mathbb{Z}_p,r]$
with the norm $|| \cdot ||_r$ on $\mathbb{A}[r]$ defined by $||f||_r=\sup_{z\in B[\mathbb{Z}_p,r]}|f(z)|_p$.
If $r>r'$, the restriction map $\mathbb{A}[r]\to \mathbb{A}[r']$ is injective, completely continuous and has dense image.
Then we define the topological spaces $\mathscr{A}(\mathbb{Z}_p)=\Liminj{s>0}\mathbb{A}[s]$ and $\mathscr{A}^{\dagger}(\mathbb{Z}_p,r)=\Liminj{s>r}\mathbb{A}[s]$
with the inductive limit topology. We have natural (continuous) inclusions
\begin{equation}\label{inclusions}
\mathscr{A}^{\dagger}(\mathbb{Z}_p,r) \hookrightarrow \mathbb{A}[r] \hookrightarrow \mathscr{A}^{\dagger}(\mathbb{Z}_p,r).
\end{equation}
Let $\mathbb{D}[r]$ be the space of continuous $\mathbb{Q}_p$-linear functionals on $\mathbb{A}[r]$.
Similarly, denote the space of continuous $\mathbb{Q}_p$-linear functionals on $\mathscr{A}(\mathbb{Z}_p)$ (resp. $\mathscr{A}^{\dagger}(\mathbb{Z}_p,r)$) by $\mathscr{D}(\mathbb{Z}_p)$ (resp. $\mathscr{D}^{\dagger}(\mathbb{Z}_p,r)$).
The topologies are given by the strong topology.
Then the space $\mathbb{D}[r]$ becomes a Banach space with the norm defined by
$$
|| \mu ||_r=\sup_{0\neq f \in \mathbb{A}[r]}\frac{| \mu (f) |_p}{||f||_r}.
$$
The family of norms $\{ ||\cdot ||_s\}_{s>0}$ (resp. $\{ ||\cdot ||_s\}_{s>r}$) induces a topology on 
$\mathscr{D}(\mathbb{Z}_p)$ (resp. $\mathscr{D}^{\dagger}(\mathbb{Z}_p,r)$).
By taking the dual of the inclusions (\ref{inclusions}), we have continuous maps
$$
\mathscr{D}(\mathbb{Z}_p)\to \mathbb{D}[r]\to \mathscr{D}^{\dagger}(\mathbb{Z}_p,r).
$$
Fix a positive integer $k$ and denote $\Sigma_0(p)=\left\{   \begin{pmatrix}  a & b \\ c & d \end{pmatrix} \in M_2 (\mathbb{Z}_p ) \, \middle| \, p\nmid a, p\mid c, ad-bc  \neq 0 \right\}$.
For $g =\begin{pmatrix}  a & b \\ c & d \end{pmatrix} \in \Sigma_0(p)$ and $f\in \mathbb{A}[r]$,
we set
$$
(g \cdot_k f)(z)=(a+cz)^k\cdot f\left( \frac{b+dz}{a+cz}\right).
$$
For $g \in \Sigma_0(p)$ and $\mu \in \mathbb{D}[r]$, we also set
$$
(\mu |_k g )(f)=\mu (g \cdot_k f).
$$
These actions induce the actions of $\Sigma_0(p)$ on $\mathscr{A}(\mathbb{Z}_p)$, $\mathscr{A}^{\dagger}(\mathbb{Z}_p,r)$, $\mathscr{D}(\mathbb{Z}_p)$, $\mathbb{D}[r]$ and $\mathscr{D}^{\dagger}(\mathbb{Z}_p,r)$.
To emphasize the role of $k$ in these actions, we write these modules as
$\mathscr{A}_k(\mathbb{Z}_p)$, $\mathscr{A}_k^{\dagger}(\mathbb{Z}_p,r)$, $\mathscr{D}_k(\mathbb{Z}_p)$, $\mathbb{D}_k[r]$  and $\mathscr{D}_k^{\dagger}(\mathbb{Z}_p,r)$.
Moreover we denote $\mathscr{D}_k=\mathscr{D}_k(\mathbb{Z}_p)$, $\mathbb{D}_k=\mathbb{D}_k[1]$ and $\mathscr{D}_k^{\dagger}=\mathscr{D}_k^{\dagger}(\mathbb{Z}_p,1)$
for simplicity.

Let $\mathscr{P}_k=\mathscr{P}_k(\mathbb{Q}_p)$ be the space of homogeneous polynomials with $\mathbb{Q}_p$-coefficients of degree $k$ in $X$ and $Y$.
We endow the space $\mathscr{P}_k$ with the structure of a right $\mathrm{GL}_2(\mathbb{Q}_p)$-module by
$$
(P|g)(x)
=P(dX-cY,-bX+aY),
$$
where $g=\begin{pmatrix}  a & b \\ c & d \end{pmatrix}\in \mathrm{GL}_2(\mathbb{Q}_p)$ and $P\in \mathscr{P}_k$.
Then we have an $\Sigma_0(p)$-equivariant map
$$\rho_k : \mathscr{D}_k^{\dagger}\to \mathscr{P}_k$$
given by
$\displaystyle \mu \mapsto \int (Y-zX)^kd\mu(z)$ and it induces
the specialization map
$$\rho_k^*:H^1_c(\Gamma_0,\mathscr{D}_k^{\dagger})\to H^1_c(\Gamma_0, \mathscr{P}_k).$$

Now we recall the slope decomposition of Banach spaces with a completely continuous operator.
Let $X$ be a Banach space over $\mathbb{Q}_p$ with a completely continuous operator $U$.
Let $Q\in \mathbb{Q}_p[T]$ be an irreducible polynomial with $Q(0)\neq 0$.
By Riesz decomposition theorem \cite[Theorem 4.1]{PS2}, we have the decomposition
into a direct sum of two closed subspaces preserved by the operator $U$:
$X\cong X(Q)\oplus X'(Q)$,
where $Q(U)$ is nilpotent on $X(Q)$ and invertible on $X'(Q)$.
Moreover $X(Q)$ is finite dimensional over $\mathbb{Q}_p$.
For an irreducible polynomial $Q\in \mathbb{Q}_p[T]$, we set $v(Q)=\sup_{\alpha} v_p(\alpha)$, where $\alpha$ runs over the all roots of $Q$ in $\overline{\mathbb{Q}}_p$.
For Banach space $X$ over $\mathbb{Q}_p$ with a completely continuous operator $U$ and $h\in \mathbb{R}$, we denote
$$
X^{(<h)}=\bigoplus_{v_p(Q)<h}X(Q),
$$
where $Q$ runs over all monic irreducible polynomials over $\mathbb{Q}_p$ satisfying $v_p(Q)<h$.
Then it is known that the space $X^{(<h)}$ is finite dimensional.

For $r \in |\mathbb{C}_p^\times |$, $\mathbb{D}_k[r]$ is a Banach space and hence $H^1_c(\Gamma_0, \mathbb{D}_k[r])$ is also a Banach space.
Then the Hecke operator $U=U_p$ acts on the space $H^1_c(\Gamma_0,\mathbb{D}_k[r])$ complete continuously.
Therefore we have the space $H^1_c(\Gamma_0,\mathbb{D}_k[r])^{(<h)}$ for each $h\in \mathbb{R}$.
Since $H^1_c(\Gamma_0,\mathscr{D}_k)$ and $H^1_c(\Gamma_0,\mathscr{D}^\dagger_k)$ are not Banach spaces, we define
$$H^1_c(\Gamma_0,\mathscr{D}_k)^{(<h)}=H^1_c(\Gamma_0,\mathscr{D}_k)\cap H^1_c(\Gamma_0,\mathbb{D}_k)^{(<h)}$$ and 
$$H^1_c(\Gamma_0,\mathscr{D}^\dagger_k)^{(<h)}=H^1_c(\Gamma_0,\mathscr{D}^\dagger_k)\cap H^1_c(\Gamma_0,\mathbb{D}_k[r])^{(<h)}$$
for any $r>1$.
Here we take the intersections as subsets of $\operatorname{Hom}(\Delta_0,\mathbb{D}_k)$ or $\operatorname{Hom}(\Delta_0,\mathbb{D}_k[r])$.
\begin{rem}
By \cite[Lemma 5.3]{PS2}, we have isomorphisms
$$
H^1_c(\Gamma_0,\mathscr{D}_k)^{(<h)} \cong H^1_c(\Gamma_0,\mathbb{D}_k)^{(<h)} \cong H^1_c(\Gamma_0,\mathscr{D}^\dagger_k)^{(<h)} \cong H^1_c(\Gamma_0,\mathbb{D}_k[r])^{(<h)}
$$
for any $h\in \mathbb{R}$ and $r>0$ with $r\in |\mathbb{C}_p^\times|$.
Therefore the definition of $H^1_c(\Gamma_0,\mathscr{D}_k)^{(<h)}$ does not depend on the choice of $r$.
\end{rem}
The following result is an analogue of Coleman's classicality theorem.
\begin{thm}(\cite[Theorem 5.12]{PS1})\label{comp}
The specialization map induces an isomorphism
$$H^1_c(\Gamma_0,\mathscr{D}_k)^{(<k+1)}\overset{\sim}{\longrightarrow} H^1_c(\Gamma_0,\mathscr{P}_k)^{(<k+1)}.$$
\end{thm}
For critical slope eigenforms, we have the following comparison theorem.
\begin{thm}(\cite[Theorem 6.7]{PS2})\label{isom}
Let $f_\beta$ be a critical slope eigenform in $S_{k+2}(\Gamma_0, \overline{\mathbb{Q}}_p)$ (i.e. $f_{\beta}$ is the critical $p$-refinement of a $p$-ordinary
eigenform $f$ in $S_{k+2}(\Gamma_1(N), \overline{\mathbb{Q}}_p)$).
Assume that the order of torsion elements in $\Gamma_0$ are prime to $p$.
Then the map induced by the specialization map
$$
H^1_c(\Gamma_0,\mathscr{D}_k)_{(f_\beta)}\overset{\sim}{\longrightarrow} H^1_c(\Gamma_0,\mathscr{P}_k)_{(f_\beta)}
$$
is an isomorphism if and only if $f_\beta$ is non-$\theta$-critical (i.e. $f_\beta$ does not belong to the image of the $\theta$-operator
$\theta^{k+1}: S_{-k}^{\dagger}(\Gamma_0, \overline{\mathbb{Q}}_p)\to S_{k+2}^{\dagger}(\Gamma_0, \overline{\mathbb{Q}}_p)$
which acts on the $q$-expansion by $\displaystyle \left( q\frac{d}{dq} \right)^{k+1}$, where $S_{r}^{\dagger}(\Gamma_0, \overline{\mathbb{Q}}_p)$
is the space of overconvergent cuspforms of weight $r$).
\end{thm}
For a normalized eigenform $f_\beta(z)$ in $S_{k+2}(\Gamma_0, \overline{\mathbb{Q}}_p)$,
we put $K=\mathbb{Q}_p(\{ a_n(f_\beta) \}_{n=1}^\infty)$. Then it is known that $(H^1_c(\Gamma_0, \mathscr{P}_k)^{\pm}\otimes K)_{(f_\beta)}$
are $1$ dimensional vector spaces over $K$. We choose $\varphi_{f_\beta}^{\pm}=\varphi_{f_\beta,\eta,u}^{\pm}$
in $(H^1_c(\Gamma_0, \mathscr{P}_k)^{\pm}\otimes K)_{(f_\beta)}$ as
$$
\varphi_{f_\beta}^\pm(\{ x\}-\{ y\} )=\frac{\pi i}{\Omega_f^{\pm}}\left( \int_x^y f_{\beta}(z)(zX+Y)^kdz \pm \int_{-x}^{-y}f_\beta (z)(zX-Y)^kdz  \right)
$$
and put $\varphi_{f_\beta}=\varphi_{f_\beta}^+ + \varphi_{f_\beta}^-$.
Recall that $\Omega_f^{\pm}=\Omega_f^{\pm}(\eta,u)$ is the period determined by $\eta$ and $u=(u^+,u^-)$.

Assume that $f_\beta$ is a critical slope eigenform in $S_{k+2}(\Gamma_0, \overline{\mathbb{Q}}_p)$.
Furthermore we suppose that $f_\beta$ is non-$\theta$-critical.
Then by Theorem \ref{isom}, there exists a unique element $\Phi_{f_\beta}$ such that the image of $\Phi_{f_\beta}$ under the specialization map
$
H^1_c(\Gamma_0,\mathscr{D}_k)_{(f_\beta)}\overset{\sim}{\longrightarrow} H^1_c(\Gamma_0,\mathscr{P}_k)_{(f_\beta)}
$
is equal to $\varphi_{f_\beta}$.
Since we have the identification $H^1_c(\Gamma_0,\mathscr{D}_k)=\operatorname{Hom}_{\Gamma_0}(\Delta_0, \mathscr{D}_k)$,
$\mu_{f_\beta} =\Phi_{f_\beta}(\{\infty \}-\{ 0 \})|_{\mathbb{Z}_p^\times}$ defines an element in $\mathscr{D}(\mathbb{Z}_p^\times)$.
Therefore we can define the critical slope $p$-adic $L$-function $L_{p,\beta}(f,\chi,s)=L_{p,\beta, \eta , u}(f,\chi,s)$ by
$$
L_{p,\beta}(f,\chi,s)=\int_{\mathbb{Z}_p^\times}\chi (x) \langle x \rangle^{s-1}d\mu_{f_\beta}
$$
for each finite order character $\chi$ of $\mathbb{Z}_p^\times$.
Then the critical slope $p$-adic $L$-function $L_{p,\beta}(f,\chi,s)$ satisfies the same interpolation property with Theorem \ref{interpolation}.
However the $p$-adic distribution $\mu_{f_\beta} $ is not characterized by the interpolation property in the critical slope case.
\begin{rem}
If a root $\gamma$ of $X^2-a_p(f)+\chi_f(p)p^{k+1}$ satisfies $\operatorname{ord}_p(\gamma )< k+1$, then 
we also have a unique lift $\Phi_{f_\gamma}$ of $\varphi_{f_\gamma}$ by Theorem \ref{comp}.
In this case, we can show that $\Phi_{f_\gamma}(\{ \infty \} - \{ 0 \})|_{\mathbb{Z}_p^\times}=\mu_{f,\gamma}$ \cite[Proposition 6.3]{PS1}.
Therefore we have an alternative construction of the $p$-adic $L$-functions for non-critical slope cuspforms using overconvergent modular symbols.
\end{rem}

\section{$p$-adic Beilinson conjecture for non-critical values of $p$-adic $L$-functions}\label{PRsect}
In this section, we give a reformulation of $p$-adic Beilinson conjecture for the special values of the cyclotomic $p$-adic $L$-functions
in the case of cuspforms with rational Fourier coefficients.
For the case of elliptic curves over $\mathbb{Q}$, the conjecture can be simplified.
Since we will focus on the case of elliptic curves over $\mathbb{Q}$ later, we also explain the formulation of the conjecture for the special case.
First, we begin with the review of Beilinson regulators and the syntomic regulators.

For a smooth scheme $X$ over a commutative ring $A$, we denote by
$\Omega^\bullet_{X/A}$ the algebraic de Rham complex, and
by $H^*_\dR(X/A):={\mathbb H}^*_\zar(X,\Omega^\bullet_{X/A})$ 
the algebraic de Rham cohomology.

\subsection{Deligne-Beilinson cohomology and real regulators}\label{DB-sect}
Let $X$ be a smooth quasi-projective variety over $\C$, and $\ol X$ a smooth
completion such that $D:=\ol X\setminus X$ is a normal crossing divisor
(abbreviated NCD). Let $j:X\hra \ol X$ be the open immersion.
Let
$\Omega^\bullet_{\ol X^\an}(\log D)$ denote the analytic 
de Rham complex with
log poles along $D$, which is quasi-isomorphic to
the de Rham complex $j_*\Omega^\bullet_{X^\an}$.
Let $\operatorname{Fil}^r\Omega_{\ol X^\an}^\bullet(\log D):=\Omega^{\bullet\geq r}_{\ol X^\an}
(\log D)$
be the stupid filtration.

For an integer $r\geq 0$,
the {\it Deligne-Beilinson complex} is
\[
\Z(r)_{\cD,X,\ol X}:=\operatorname{Cone}[Rj_*\Z_X(r)\op 
\operatorname{Fil}^r\Omega^\bullet_{\ol X^\an}(\log D)\os{\epsilon-i}{\lra}
j_*\Omega^\bullet_{X^\an}]
\]
the cone of complexes of sheaves on the analytic site $\ol X^{\an}$ 
where $\Z_X(r):=(2\pi i)^r\Z_X$
and 
$\epsilon:\Z_X(r)\to \Omega^\bullet_{X^\an}$ is the natural map, and
$i:\mathrm{Fil}^r\Omega^\bullet_{\ol X^\an}(\log D)\to \Omega^\bullet_{\ol X^\an}(\log D)
\to j_*\Omega^\bullet_{X^\an}$ is the composition.
The {\it Deligne-Beilinson cohomology} is defined to be the 
cohomology 
\[
H^i_\cD(X,\Z(r)):={\mathbb H}^i(\ol X^{an},\Z(r)_{\cD,X,\ol X}).
\]
This depends only on $X$, and functorial with respect to $X$ (\cite[Lemma 2.8]{EV}).
The Deligne-Beilinson cohomology $H^i_\cD(X,A(r))$ with
coefficients in a ring $A\subset \R$ is defined by replacing $\Z$ with $A$ in the above.
The distinguished triangle
\[
\operatorname{Cone}[
\mathrm{Fil}^r\Omega^\bullet_{\ol X^\an}(\log D)\os{\epsilon}{\to}
j_*\Omega^\bullet_{X^\an}]\lra \Z(r)_{\cD,X,\ol X}
\lra Rj_*\Z_X(r)
\]
gives rise to an exact sequence
\begin{equation}\label{DB-eq1}
\cdots\to H^{i-1}_\dR(X/\C)/\Fil^r\to H^i_\cD(X,\Z(r))\lra H^i_\mathrm{B}(X(\C),\Z(r))\to
H^i_\dR(X/\C)/\Fil^r\to\cdots
\end{equation}
where $\Fil^\bullet =\Fil^\bullet H_\dR(X/\C)$ denotes Deligne's Hodge filtration.
Here we note that the analytic de Rham cohomology ${\mathbb H}^i(X^\an,\Omega_{X^\an}^\bullet)$
agrees with the algebraic de Rham cohomology $H^i_\dR(X/\C)$.

Let $X_\R$ be a smooth variety over $\R$, and $\ol X_\R$ a smooth compactification
such that $D_\R:=\ol X_\R\setminus X_\R$ is a NCD. 
Let $\ol X=\ol X_\R\times_\R\C$ etc. and denote by $\ol X(\C)$ the complex points.
Then the complex conjugation $z\to\bar z$ induces the
anti-holomorphic map $F_\infty:\ol X(\C)\to \ol X(\C)$, which is called
the {\it infinite Frobenius}.
Let $\omega\mapsto \ol\omega$ denote the complex conjugation on
smooth differential forms on $X(\C)$.
The map
\[
\omega\longmapsto \ol{F_\infty^*(\omega)},\quad \omega\in \Omega^i_{\ol X^\an}(\log D)
\]  
induces maps
\[
F_\infty^{-1}\Z_X(r)\lra \Z_X(r),\quad
F_\infty^{-1}\Z(r)_{\cD,X,\ol X}\lra \Z(r)_{\cD,X,\ol X}
\]
of sheaves on the smooth manifold $\ol X(\C)$.
We thus have an action on the cohomology
$H^i_\cD(X,\Z(r))$, $H_\mathrm{B}(X(\C),\Z(r))$ and so on,
 which we write by the same notation
$F_\infty$.
By definition $F_\infty^2=\id$.
We write $H^+$ to be the $F_\infty$-fixed part, and
$H^-$ the anti-fixed part.
We have
\[
H^i_\mathrm{B}(X(\C),\Z(r))^+=
H^i_\mathrm{B}(X(\C),\Z(s))^{(-1)^{r-s}}\ot\Z(r-s).
\]
Define the {\it real Deligne-Beilinson cohomology}
$H^i_\cD(X_\R\slash \R,\R(r))$ to be the $F_\infty$-fixed part (\cite[\S 2]{schneider})
\[
H^i_\cD(X_\R\slash \R,\R(r)):=
H^i_\cD(X,\R(r))^+.
\]

\medskip

By the theory of the universal Chern classes, the {\it Beilinson regulator map}
(the higher Chern class map) 
\begin{equation}
\xymatrix{
\reg_\cD^{i,r}:K_{2r-i}(X)\ar[r]& H^{i}_\cD(X,\Z(r))
}
\end{equation}
is defined (\cite[\S 4, p.29--30]{schneider}).
Let $X_\R$ be a smooth variety over $\R$, then it induces the real regulator map
\begin{equation}
\xymatrix{
\reg_\cD^{i,r}:K_{2r-i}(X_\R)\ar[r]& H^{i}_\cD(X_\R\slash \R,\R(r))
}
\end{equation}
to the real Deligne-Beilinson cohomology group.

\subsection{Syntomic cohomology and Syntomic regulators}\label{syn-sect}
Let $W$ be the Witt ring of a perfect field $k$ of characteristic $p>0$, and put $K:=
\operatorname{Frac} W$.
Let $X$ be a smooth projective $W$-scheme.
We write $X_k:=X\times_Wk$ and $X_K:=X\times_WK$.
Then the {\it syntomic cohomology} of Fontaine and Messing
\[
H^i_\syn(X,\Z_p(r)),\quad r\geq 0
\]
is defined (\cite{FM}, \cite{KaV}). There is an exact sequence
\begin{equation}\label{syn-eq1}
\cdots\to \Fil^rH^{i-1}_\crys(X_k/W)\os{1-\frac{\phi}{p^r}}{\lra} H^{i-1}_\crys(X_k/W)\to
H^i_\syn(X,\Z_p(r))\to \Fil^rH^{i}_\crys(X_k/W)\to\cdots
\end{equation}
which is a counterpart of \eqref{DB-eq1}, where 
$\phi$ is the $p$-th Frobenius on the crystalline cohomology $H^i_\crys(X_k/W)$,
and 
$\Fil^\bullet$ is the Hodge filtration via the comparison
$H^i_\crys(X_k/W)\cong H^i_\dR(X/W)$.
If we drop the assumption that $X$ is proper over $W$, the construction of
the syntomic cohomology in \cite{FM} or \cite{KaV} does not give a ``correct"
cohomology theory.
For such a variety, we need 
the {\it rigid syntomic cohomology} $H^i_{\text{rig-syn}}(X,\Q_p(r))$
by Besser \cite{Be1} (see also \cite[1B]{NN}), which
gives the correct cohomology theory.

\medskip

Let $X$ be a smooth $W$-scheme.
Then the {\it syntomic regulator map}
\begin{equation}\label{syn-eq2}
\xymatrix{
\reg_{\text{rig-syn}}^{i,r}:K_{2r-i}(X)\ar[r]& H^{i}_{\text{rig-syn}}(X,\Q_p(r))
}
\end{equation}
is defined (\cite[Theorem 7.5]{Be1}, \cite[Theorem A]{NN}).
This is compatible with the \'etale regulator map
\begin{equation}\label{syn-eq3}
\xymatrix{
\reg_\et^{i,r}:K_{2r-i}(X_K)\ar[r]& H^{i}_\et(X_K,\Q_p(r))
}
\end{equation}
via the syntomic period map (loc.cit. Theorem A (7))
\[
c:H^{i}_{\text{rig-syn}}(X,\Q_p(r))\lra H^{i}_\et(X_K,\Q_p(r)).
\]
If $X$ is projective over $W$, the rigid syntomic cohomology agrees with the syntomic cohomology of Fontaine-Messing.
In this case we write the subscript ``$\syn$"
instead of ``$\mathrm{rig}$-$\mathrm{syn}$". 

Let $X$ be a smooth projective $W$-scheme.
The composition $\rho$ in the following diagram
\[
\xymatrix{
H^{i-1}_\crys(X_k/W)/(1-p^{-r}\phi)\Fil^r\ar[r]^{\hspace{1cm}\eqref{syn-eq1}}
&H^i_\syn(X,\Q_p(r))\ar[r]^c&H^i_\et(X,\Q_p(r))\\
H^{i-1}_\crys(X_k/W)/\Fil^r\ar[u]^{1-p^{-r}\phi}\ar[rru]_\rho
}
\]
induces a map
\[
H^{i-1}_\crys(X_k/W)/\Fil^r\lra H^1(G_K,H^{i-1}_\et(X_{\ol K},\Q_p(r))).
\]
This agrees with the Bloch-Kato exponential map 
(\cite[Proposition 9.11]{Be1}, \cite[Proposition 1.1]{NN}).

\subsection{$p$-adic Beilinson conjecture for modular forms}
In this section, we assume that the Fourier coefficients of $f$ belong to $\mathbb{Q}$ for simplicity.
For a projective smooth variety over $\mathbb{Q}$, let $H_{\mathscr{M}}^i(X,\mathbb{Q}(j))=K_{2j-i}(X)_{\mathbb{Q}}^{(j)}$ be the motivic cohomology group and
let $H_{\mathscr{M}}^i(X,\mathbb{Q}(j))_{\mathbb{Z}}$ denote the integral part defined by Scholl \cite{Sc2}.
Recall that $\pi : \mathscr{E}\to X_1(N)$ is the universal generalized elliptic curve over the modular curve.
Let $\mathscr{E}^k=\mathscr{E}\times_{X_1(N)}\cdots \times_{X_1(N)} \mathscr{E}$ be the $k$-fold fiber product of $\mathscr{E}$ over
the modular curve $X_1(N)$ and $W_k$ the canonical desingularization of $\mathscr{E}^k$(See \cite{Sc} or an appendix in \cite{BDP}
by Brian Conrad for details).
$W_k$ is a smooth projective variety over $\mathbb{Q}$ of dimension $k+1$.
Fix a prime $p$ not dividing $N$.
Let $\mathcal{W}_k$ be the integral model over $\mathbb{Z}_p$.
For $a\in (\mathbb{Z}\slash N\mathbb{Z})^k$, let $\sigma_a$ denote the automorphism on $\mathscr{E}^k$ obtained from the translation by the sections of order $N$,
which naturally extends to $W_k$. Then we define the projector $\varepsilon_k^{(1)}$ by
$$
\varepsilon_k^{(1)}=\frac{1}{N^k}\sum_{a\in (\mathbb{Z}\slash N\mathbb{Z})^k}\sigma_a.
$$
Let $\mathfrak{S}_k$ be the symmetric group on $k$ letters.
The multiplication by $-1$ on $\mathscr{E}$ and the natural permutation action of $\mathfrak{S}_k$ on $\mathscr{E}^k$
give rise to an action of the semi-direct product $\{ \pm 1 \}^k \rtimes \mathfrak{S}_k$ on $\mathscr{E}^k$.
This action also extends to $W_k$. Let $\mu : \{ \pm 1 \}^k \rtimes \mathfrak{S}_k \to \{ \pm 1 \}$ be the character
which is the multiplication on $\{ \pm 1 \}^k$ and the sign character on $\mathfrak{S}_k$.
Then
$$
\varepsilon_k^{(2)}=\frac{1}{2^k\cdot k!}\sum_{\sigma \in  \{ \pm 1 \}^k \rtimes \mathfrak{S}_k}\mu (\sigma )\sigma 
\in \mathbb{Q} [\operatorname{Aut}(W_k\slash X_1(N))]
$$
defines a projector.
The projectors $\varepsilon_k^{(1)}$ and $\varepsilon_k^{(2)}$ commute, and therefore the composition
$\varepsilon_k =\varepsilon_k^{(1)}\varepsilon_k^{(2)}$ defines a projector.
Then the de Rham realization of the Chow motive $M=(W_k, \varepsilon_k)$ is canonically isomorphic to
the parabolic de Rham cohomology:
$$
H^*_{\mathrm{dR}}(M\slash \mathbb{Q})=\varepsilon_k H^*_{\mathrm{dR}}(W_k\slash \mathbb{Q})=\varepsilon_k H^{k+1}_{\mathrm{dR}}(W_k\slash \mathbb{Q})
\cong H^1_{\mathrm{par}}(X_1(N),\operatorname{Sym}^k\mathscr{L},\nabla),
$$
where $\mathscr{L}=\mathbb{R}^1\pi_*(0\to \mathscr{O}_{\mathscr{E}}\to \Omega^1_{\mathscr{E}\slash X_1(N)}\to 0)$
is the relative de Rham cohomology sheaf on $X_1(N)$.
The space of cuspforms $S_{k+2}(\Gamma_1(N),\Q)= H^0(X_1(N), \underline{\omega}^k\otimes \Omega_{X_1(N)}^1)$
is naturally identified with a subspace of the parabolic de Rham cohomology
$H^1_{\mathrm{par}}(X_1(N),\operatorname{Sym}^k\mathscr{L},\nabla)$.
Moreover, we have
$\operatorname{Fil}^0 \varepsilon_k H^*_{\mathrm{dR}}(W_k\slash \Q)=H^1_{\mathrm{par}}(X_1(N),\operatorname{Sym}^k\mathscr{L},\nabla)$,
$$\operatorname{Fil}^1 \varepsilon_k H^*_{\mathrm{dR}}(W_k\slash \Q)=\cdots =
\operatorname{Fil}^{k+1}\varepsilon_k H^*_{\mathrm{dR}}(W_k\slash \Q)=H^0(X_1(N), \underline{\omega}^k\otimes \Omega_{X_1(N)}^1)$$
and $\operatorname{Fil}^{k+2} \varepsilon_k H^*_{\mathrm{dR}}(W_k\slash \Q)=0$,
where $\operatorname{Fil}^{i}$ denotes the $i$-th step in the Hodge filtration on $\varepsilon_k H^*_{\mathrm{dR}}(W_k\slash \Q)$.
For the Betti realization, we have a canonical identification
$$
H^*_{\mathrm{B}}(M, \Q)=\varepsilon_k H^*_{\mathrm{B}}(W_k(\mathbb{C}), \Q)=\varepsilon_k H^{k+1}_{\mathrm{B}}(W_k(\mathbb{C}), \Q)
\cong H^1_{\mathrm{par}}(X_1(N)(\mathbb{C}),\mathscr{L}_{k,\Q}^{\mathrm{B}}),
$$
where $\mathscr{L}_{k,\Q}^{\mathrm{B}}=\mathscr{L}_{k}^{\mathrm{B}}\otimes_{\mathbb{Z}} \Q$ and
$$H^1_{\mathrm{par}}(X_1(N)(\mathbb{C}),\mathscr{L}_{k,\Q}^{\mathrm{B}})=
\operatorname{Im}\left[ H^1_c(X_1(N)(\mathbb{C}),\mathscr{L}_{k,\Q}^{\mathrm{B}})\to H^1(X_1(N)(\mathbb{C}),\mathscr{L}_{k,\Q}^{\mathrm{B}})\right]$$
is the parabolic Betti cohomology.

Let $f$ be a normalized newform in $S_{k+2}(\Gamma_1(N),\Q)$
(Note that $f=f^*$ in this case).
Then the element $u^\pm$ is determined only up to sign.
Also we fix a non-zero element $\eta$ in $S_{k+2}(\Gamma_1(N),\Q)[f]= H^0(X_1(N), \underline{\omega}^k\otimes \Omega_{X_1(N)})[f]$.
As we mentioned above, $\eta$ can be identified with an element in $\varepsilon_k H^{k+1}_{\mathrm{dR}}(W_k\slash \Q)[f] $
and $u^{\pm}$ can be viewed as an element in 
$$H^*_{\mathrm{B}}(M, \Q)[f]^\pm
=\varepsilon_k H^{*}_{\mathrm{B}}(W_k(\mathbb{C}), \Q)^{\pm}[f]
=\varepsilon_k H^{k+1}_{\mathrm{B}}(W_k(\mathbb{C}), \Q)^{\pm}[f].
$$
Let $\gamma$ be an eigenvalue of the Frobenius $\Phi$ on
$$(\varepsilon_k H^{k+1}_{\mathrm{dR}}(W_k\slash \mathbb{Q}_p)\otimes \overline{\mathbb{Q}}_p) [f] \cong 
(\varepsilon_k H^{k+1}_{\mathrm{crys}}(\mathcal{W}_{k, \mathbb{F}_p}\slash \mathbb{Z}_p)\otimes \overline{\mathbb{Q}}_p)[f]$$
and choose a non-zero eigenvector $v_{\gamma}$ in $(\varepsilon_k H^{k+1}_{\mathrm{dR}}(W_k\slash \mathbb{Q}_p) \otimes {\mathbb{Q}}_p(\gamma))[f]$ with the eigenvalue $\gamma$.
Let $$\mathrm{reg}_{\mathrm{syn}}: \varepsilon_k H_{\mathscr{M}}^{k+2}(W_k ,\mathbb{Q}(n+k+2))_{\mathbb{Z}}\to \varepsilon_k H^{k+2}_{\mathrm{syn}}(W_k,\mathbb{Q}_p(n+k+2))\cong \varepsilon_k H^{k+1}_{\mathrm{dR}}(W_k \slash \mathbb{Q}_p)$$
be the syntomic regulator map.
Let
$$
\mathrm{pr}_{\mathrm{dR},f} : \varepsilon_k H^{k+1}_{\mathrm{dR}}(W_k \slash \mathbb{Q}_p)\to 
\varepsilon_k H^{k+1}_{\mathrm{dR}}(W_k \slash \mathbb{Q}_p)[f]
$$
denote the projection to the $f$-isotypic subspace.
Assume that $\eta$ and $v_{\gamma}$ are linearly independent over $\overline{\mathbb{Q}}_p$.
For an element $\xi \in \varepsilon_k H_{\mathscr{M}}^{k+2}(W_k ,\mathbb{Q}(n+k+2))_{\mathbb{Z}}$, we set
$$
R_{p,\gamma ,\eta}(M(f)(-n), \xi)=\Gamma^*(-n)(1-p^{1+n}\gamma^{-1})\frac{\operatorname{Tr}(\mathrm{pr}_{\mathrm{dR},f} \circ \mathrm{reg}_{\mathrm{syn}}(\xi)\cup v_{\gamma})}{\operatorname{Tr}(\eta \cup v_{\gamma})}
\in {\mathbb{Q}}_p (\gamma ),
$$
where
$$
\Gamma^*(-n)=
\begin{cases}\displaystyle
(-n-1)! & \textup{if }n\leq -1,\\
\displaystyle \frac{(-1)^n}{n!} & \textup{if } n\geq 0.
\end{cases}
$$
Note that $R_{p,\gamma,\eta}(M(f)(-n), \xi)$ is independent of the choice of $v_{\gamma}$ by the definition.
On the other hand, we have the Beilinson regulator map
\begin{align*}
\mathrm{reg}_{\mathscr{D}}:\varepsilon_k H_{\mathscr{M}}^{k+2}(W_{k},\mathbb{Q}(n+k+2))_{\mathbb{Z}}
&\to \varepsilon_k H^{k+2}_{\mathscr{D}}(W_{k,\mathbb{R}}\slash \mathbb{R},\mathbb{R}(n+k+2))\\
&\cong \varepsilon_k H^{k+1}_{\mathrm{B}}(W_k (\mathbb{C}),\mathbb{R}(n+k+1))^+\\
&\cong \varepsilon_k H^{k+1}_{\mathrm{B}}(W_k (\mathbb{C}),\mathbb{R})^{(-1)^{n+k+1}}\otimes \mathbb{R}(n+k+1).
\end{align*}
Choose a non-zero element $\delta^\pm \in \varepsilon_k H^{k+1}_{\mathrm{B}}(W_k (\mathbb{C}),\mathbb{R})^\pm [f]$.
We denote the projection to the $f$-isotypic subspace by
$$
\operatorname{pr}_{\mathrm{B},f}:
\varepsilon_k H^{k+1}_{\mathrm{B}}(W_k (\mathbb{C}),\mathbb{R})^{(-1)^{m}}\otimes \mathbb{R}(m)
\rightarrow \varepsilon_k H^{k+1}_{\mathrm{B}}(W_k (\mathbb{C}),\mathbb{R})^{(-1)^{m}}[f]\otimes \mathbb{R}(m),
$$
where $m=n+k+1$.
Now we define
$$\displaystyle R_{\infty}(M(f)(-n),\xi)=\frac{1}{(2 \pi i)^{n+k+1}}
\frac{\operatorname{Tr}(\mathrm{pr}_{\mathrm{B},f} \circ \mathrm{reg}_{\mathscr{D}}(\xi)\cup \delta^{(-1)^{n+k+1}})}{\operatorname{Tr}(u^{(-1)^{n+k+1}} \cup \delta^{(-1)^{n+k+1}})}
\in \mathbb{R}.$$
By the definition, it is clear that $R_{\infty}(M(f)(-n),\xi)$ is independent of the choice of $\delta^\pm$.
By the Beilinson conjecture, it is expected that the ratio $\displaystyle \frac{L'(f,-n)}{R_{\infty}(M(f)(-n), \xi)}$ is a rational number.
\begin{conj}[$p$-adic Beilinson conjecture for modular forms \cite{PR}]\label{PRMF}
Suppose that one of the following conditions holds:
\begin{enumerate}
\item
$\mathrm{ord}_p(\gamma)<k+1$,
\item
$\mathrm{ord}_p(\gamma)=k+1$ and the $p$-refinement $f_{\gamma}$ is not $\theta$-critical.
\end{enumerate}
Moreover we assume that $\eta$ and $v_{\gamma}$ are linearly independent over $\mathbb{Q}_p$.
Then for $n\geq 0$, there exists a non-zero element
$$\xi \in \varepsilon_k H_{\mathscr{M}}^{k+2}(W_k ,\mathbb{Q}(n+k+2))_{\mathbb{Z}}$$
such that
the ratio $\displaystyle \frac{L_{p,\gamma,\eta , u }(f,\omega^{-n-1},-n)}{R_{p,\gamma,\eta}(M(f)(-n),\xi)}$ is a non-zero rational number
and we have
$$
\frac{L_{p,\gamma,\eta , u }(f,\omega^{-n-1},-n)}{R_{p,\gamma,\eta}(M(f)(-n),\xi)}
=\pm \frac{L'(f,-n)}{R_{\infty}(M(f)(-n), \xi)}.
$$
\end{conj}
\begin{rem}\label{rem-PR-mf}
\begin{enumerate}
\item
In the original formulation of the $p$-adic Beilinson conjecture by Perrin-Riou \cite[4.2.2 Conjecture (CP($M,\tau$))]{PR}, it is expected the existence of $p$-adic $L$-function satisfying
the interpolation property in the critical range and the $p$-adic Beilinson formula in the non-critical range.
For the case of elliptic modular forms, we already have an appropriate $p$-adic $L$-function $L_{p,\gamma}(f,\chi,s)$ satisfying
the interpolation property in the critical range. Therefore we can formulate
the conjecture as more accessible form in this case.
Moreover, Perrin-Riou used the \'etale regulator map to formulate the conjecture. Therefore the formula of the conjecture has a slightly different form.
In general, the comparison between the \'etale regulator and the syntomic regulator is explained by Besser \cite{Be1}.
By \cite[Remark 8.6.3]{Be1}, the difference is given by the factor $(1-\frac{\Phi}{p^n})$.
Using this result, it is easy to see that our formulation is equivalent to Perrin-Riou's formulation.
\item
Our normalization of the $p$-adic $L$-functions is slightly different from Perrin-Riou's normalization.
In fact, Perrin-Riou's $p$-adic $L$-functions should be the half of our $p$-adic $L$-function.
For details, see \cite[4.2.2 Conjecture (CP($M,\tau$)) and 4.3.4]{PR}.
\item
For $c\in \mathbb{Q}^\times$, it is easy to see
$$L_{p,\gamma, c \eta , u }(f,\omega^{-n-1},-n)=c^{-1}\cdot L_{p,\gamma,\eta , u }(f,\omega^{-n-1},-n)$$
and
$$R_{p,\gamma,c\eta}(M(f)(-n),\xi)=c^{-1} \cdot R_{p,\gamma,\eta}(M(f)(-n),\xi).$$
Therefore the ratio $\displaystyle \frac{L_{p,\gamma,\eta , u }(f,\omega^{-n-1},-n)}{R_{p,\gamma,\eta}(M(f)(-n),\xi)}$ does not depend
on the choice of $\eta$.
Moreover, we can show that the conjecture \ref{PRMF} is also independent of the choice of $u =(u^+,u^-)$.
\item
Bannai-Kings \cite{BK} proved a version of the $p$-adic Beilinson conjecture for Hecke motives over imaginary quadratic fields with class number one
when $\operatorname{ord}_p(\gamma)=0$.
\item
In this paper, we only consider the case that $p$ does not divide $N$ as in Perrin-Riou's book \cite{PR}.
If $p$ divides $N$ exactly once, we also have a $p$-adic $L$-function \cite{AV, V}.
More generally, Rodrigues Jacinto \cite{Ro} have constructed a (de Rham cohomology valued) $p$-adic $L$-function
for general case (including the case $p^2$ divides $N$).
Therefore it should be interesting to consider the $p$-adic Beilinson conjecture for the case of bad reduction
using his $p$-adic $L$-function.  
\end{enumerate}
\end{rem}

\subsection{$p$-adic Beilinson conjecture for elliptic curves}\label{PRconjEC}
Let $E$ be an elliptic curve over $\Q$ of conductor $N$ and let $\omega_E$ denote the N\'eron differential.
Let $u^{\pm}$ be the generator of $H^1(E(\mathbb{C}), \mathbb{Z})^{\pm}$ satisfying $\int_{u^{+}}\omega_E>0$ and
$i^{-1}\int_{u^{-}}\omega_E>0$.
Then we put $\Omega^+_E=\int_{u^{+}}\omega_E$ and $\Omega^-_E=\int_{u^{-}}\omega_E$.
It is known that
$$
\frac{L(E, \overline{\chi},1)}{\Omega_{E}^{\chi (-1)}} \in \mathbb{Q}(\chi)
$$
for any Dirichlet character $\chi$.
Let $f=f_E$ be the newform associated to $E$.
For $a,m \in \mathbb{Q}$, we put
$$
\lambda_E^{\pm}(a,m)=\frac{\pi i}{\Omega_E^{\pm}}\left\{  \int^{\frac{a}{m}}_{\infty}f(z)dz \pm \int^{-\frac{a}{m}}_{\infty}f(z)dz  \right\} .
$$
Assume that $E$ has good reduction at $p$.
Let $\gamma$ be a root of $X^2-a_p(E)X+p$.
Using the modular symbol $\lambda^{\pm}_E$, we can construct the $p$-adic distribution $\mu_{E,\gamma}$ for any choice of the root $\gamma$ and we can define the $p$-adic $L$-function
$\displaystyle L_{p,\gamma}(E,\chi,s)=\int_{\mathbb{Z}_p^\times}\chi (x) \langle x \rangle^{s-1} d\mu_{E,\gamma}$
for $s\in \mathbb{Z}_p$. In particular, we have the following interpolation property:
For a primitive Dirichlet character $\chi : (\mathbb{Z}\slash p^\nu \mathbb{Z})^\times \to \overline{\mathbb{Q}}^\times$,
$$
L_{p,\gamma}(E,\chi,1)=
\begin{cases}\displaystyle
(1-\gamma^{-1})^2 \frac{L(E,1)}{\Omega_E^+} & \textup{if } \chi = \mathbf{1},\\
\displaystyle
\left( \frac{p}{\gamma}\right)^{\nu}\cdot \frac{1}{\tau (\overline{\chi})}\cdot \frac{L(E,\overline{\chi},1)}{\Omega_E^{\chi (-1)}}& \textup{if }\chi \neq \mathbf{1}.
\end{cases}
$$

For elliptic curves over $\mathbb{Q}$, the statement of the $p$-adic Beilinson conjecture can be slightly simplified as follows.
Note that $\gamma$ is also an eigenvalue of the Frobenius $\Phi$ on $H^1_{\mathrm{dR}}(E\slash \mathbb{Q}_p)\cong H^1_{\mathrm{crys}}(E_{\mathbb{F}_p}\slash \mathbb{Z}_p)\otimes \mathbb{Q}_p$.
Choose a non-zero eigenvector $v_{\gamma}$ in $H^1_{\mathrm{dR}}(E\slash \mathbb{Q}_p)\otimes \overline{\mathbb{Q}}_p$ with the eigenvalue $\gamma$.
Let $$\mathrm{reg}_{\mathrm{syn}}: H_{\mathscr{M}}^2(E,\mathbb{Q}(n+2))_{\mathbb{Z}}\to H^2_{\mathrm{syn}}(E,\mathbb{Q}_p(n+2))\cong H^1_{\mathrm{dR}}(E\slash \mathbb{Q}_p)$$
be the syntomic regulator map.
Suppose $\Phi (\omega_E)\neq \gamma \omega_E$. Note that this is equivalent to the condition that $\omega_E$ and $v_{\gamma}$ are linearly independent. 
Then we set
$$
R_{p,\gamma}(h^1(E)(-n), \xi)=\Gamma^*(-n)(1-p^{n+1}\gamma^{-1})\frac{\operatorname{Tr}(\mathrm{reg}_{\mathrm{syn}}(\xi)\cup v_{\gamma})}{\operatorname{Tr}(\omega_E \cup v_{\gamma})}
\in \mathbb{Q}_p(\gamma).
$$
Let
\begin{align*}
\mathrm{reg}_{\mathscr{D}}:H_{\mathscr{M}}^2(E,\mathbb{Q}(n+2))_{\mathbb{Z}}&\to H^2_{\mathscr{D}}(E,\mathbb{R}(n+2))\\
&\cong H^1_B(E(\mathbb{C}),\mathbb{R}(n+1))^+\\
&\cong \operatorname{Hom}(H_1(E(\mathbb{C}),\mathbb{Q})^{(-1)^{n+1}},\mathbb{R}(n+1)).
\end{align*}
be the Beilinson regulator map
and we denote
$$\displaystyle R_{\infty}(h^1(E)(-n),\xi)=\frac{1}{(2 \pi i)^{n+1}}\mathrm{reg}_{\mathscr{D}}(\xi)(u^{(-1)^{n+1}})\in \mathbb{R}.$$
By the Beilinson conjecture, it is expected that the ratio
$\displaystyle \frac{L'(E,-n)}{R_{\infty}(h^1(E)(-n), \xi)}$ is a rational number.
\begin{conj}[$p$-adic Beilinson conjecture for elliptic curves over $\mathbb{Q}$]\label{PREC}
Assume that one of the following conditions holds:
\begin{enumerate}
\item
$\mathrm{ord}_p(\gamma)<1$,
\item
$\mathrm{ord}_p(\gamma)=1$ and $f_{\gamma}$ is not $\theta$-critical.
\end{enumerate}
Moreover we suppose that $\Phi (\omega_E)\neq \gamma\omega_E$.
Then for $n\geq 0$, there exists a non-zero element $\xi \in H_{\mathscr{M}}^2(E,\mathbb{Q}(n+2))_{\mathbb{Z}}$
such that
$$
\frac{L_{p,\gamma}(E,\omega^{-n-1},-n)}{R_{p,\gamma}(h^1(E)(-n),\xi)}
=\pm \frac{L'(E,-n)}{R_{\infty}(h^1(E)(-n), \xi)}.
$$
\end{conj}
\begin{rem}\label{rem-PR-ec}
\begin{enumerate}
\item The assumption $\Phi (\omega_E)\neq \gamma\omega_E$ excludes only the case that $E$ has good reduction at $p$, $E$ has complex multiplication
and $\gamma$ has critical slope.
\item By the rank part of the Beilinson conjecture, it is expected that
$$\operatorname{dim}_\mathbb{Q} H_{\mathscr{M}}^2(E,\mathbb{Q}(n+2))_{\mathbb{Z}}=1.$$
\item
Brunault \cite{Br} gave an explicit formula for \'etale regulator of Beilinson-Kato elements.
\item
Bertolini-Darmon \cite{BD} gave an explicit formula for syntomic regulator and Beilinson regulator of Beilinson-Kato elements
in the case of weight two cuspforms under some technical assumptions.
However, their result does not cover Conjecture \ref{PREC}.
\item
In \cite{New}, we proposed conjectures on the relation between the special values of $p$-adic $L$-functions $L_{p,\alpha}(E,\omega^{-1},0)$
and the special values of $p$-adic hypergeometric functions (\cite[Conjecture 4.30 -- 4.35]{New}).
Using the results in \cite{New}, it is easy to see that these conjectures follow from Conjecture \ref{PREC} and the rank part of the Beilinson conjecture.
\end{enumerate}
\end{rem}
For $n=0$, we will denote $R_{p,\gamma}(E,\xi)=R_{p,\gamma}(h^1(E),\xi)$ and $R_{\infty}(E, \xi)=R_{\infty}(h^1(E), \xi)$ for simplicity.


\section{Syntomic regulators and overconvergent functions}\label{Tate-sect}
In this section we work over the Witt ring $W=W(\F)$ 
of a perfect field $\F$ of characteristic $p>0$
endowed with the Frobenius $F$.
Put $K:=\operatorname{Frac}(W)$ the fractional field.
Let $A$ be a faithfully flat $W$-algebra.
We mean by a $p^n$-th Frobenius on $A$ an endomorphism $\sigma$ 
such that $\sigma(x)\equiv x^{p^n}$ mod $pA$ for all $x\in A$ and that 
$\sigma$ is compatible with the 
$p^n$-th Frobenius on $W$. We also write $x^\sigma$ instead of $\sigma(x)$. 
For a $W$-algebra $A$ of finite type,
we denote by $A^{\dag}$ the weak completion.
Namely if $A=W[T_1,\ldots,T_n]$, then 
$A^\dag=W[T_1,\ldots,T_n]^\dag$
is the ring of power series $\sum a_\alpha T^\alpha$ such that
for some $r>1$,
$|a_\alpha|r^{|\alpha|}\to0$ as $|\alpha|\to\infty$, and 
if $A=W[T_1,\ldots,T_n]/I$, then
$A^\dag=W[T_1,\ldots,T_n]^\dag/IW[T_1,\ldots,T_n]^\dag$.

\subsection{Setting and Notation}\label{fib-set-sect}
Let $C$ be a smooth projective scheme over $W$ of relative dimension $1$.
Let
\[
f:Y\lra C
\]
be an elliptic fibration over $W$, which means that
$Y$ is a smooth $W$-scheme and $f$ is a projective flat morphism
with a section $e:C\to Y$
whose general fiber is an elliptic curve.
Let $Z\subset C$ 
be a closed subscheme such that $f$ is smooth over $S:=C\setminus Z$.
Put $X:=f^{-1}(S)$ and $A:=\O(S)$.
Write $A_K:=A\ot_WK$, $Y_K:=Y\times_WK$, $A_\F:=A/pA$, $Y_\F:=Y\times_WW/pW$ etc.

We suppose that the following conditions {\bf A1}, \ldots, {\bf A3} are satisfied.
\begin{description}
\item[A1]
$Z$ is a disjoint union of $W$-rational points $\{P_i\in S(W)\}$, and
each $D_i:=f^{-1}(P_i)$ is a relative NCD (possibly a smooth divisor) over $W$.
\item[A2]
$H^0_\zar(X,\Omega^1_{X/A})$ is a free $A$-module of rank one. 
\end{description}
The condition {\bf A2} is equivalent to that $H^1_\zar(\O_X)$ is a free $A$-module of rank one
by the Serre duality. Then 
$H^1_\dR(X/A)$ is a free $A$-module of rank $2$ by Lemma \ref{A3-lem}.
\begin{description}
\item[A3]
There is a point $P_0$ such that 
$f$ has a split multiplicative reduction at a point $P_0$.
Let $\lambda$ be the uniformizer of the local ring at the point $P_0$ and
let $\cE_{W[[\lambda]]}:=Y\times_C\operatorname{Spec} W[[\lambda]]$ be the Tate curve over the formal neighborhood
of $P_0$. Let 
\begin{equation}\label{global-mult-period}
q=q(\lambda)=\kappa \lambda^n(1+a_1\lambda+\cdots)\in\lambda W[[\lambda]]
\end{equation} 
be the multiplicative period.
Then $\kappa\in W^\times$ and $n$ is prime to $p$.
\end{description}
The multiplicative period $q$ is characterized by the functional $j$-invariant $j(\lambda)$,
namely it is the unique power series satisfying
\[
j(\lambda)=\frac{1}{q}+744+196884q+\cdots.
\]
Let $\tau(\lambda)\in \lambda K[[\lambda]]$ be defined by $q=\kappa \lambda^n\exp(\tau(\lambda))$, 
\begin{equation}\label{tau}
\tau(\lambda)=\log(1+a_1\lambda+\cdots)\in \lambda K[[\lambda]].
\end{equation}
See Proposition \ref{A-sympl-3} \eqref{A-sympl-3-eq1} for differential equations of $q$ and $\tau(\lambda)$.
\subsection{Category of Filtered $F$-isocrystal}
Let $A^\dag$ be the weak completion of $A$.
We write $A^\dag_K=A^\dag\ot_WK$.
Fix a $p$-th Frobenius $\sigma$ on $A^\dag$.
In \cite[2.1]{AM} we introduced
a category $\FilFMIC(S)=\FilFMIC(S,\sigma)$, which 
consists of collections of datum $(H_{\dR}, H_{\rig}, c, \Phi, \nabla, \Fil^{\bullet})$ such that
    \begin{itemize}
            \setlength{\itemsep}{0pt}
        \item $H_{\dR}$ is a finitely generated $A_K$-module,
        \item $H_{\rig}$ is a finitely generated $A^\dag_K$-module,
        \item $c\colon H_\rig\cong H_{\dR}\otimes_{A_K}A^\dag_K$, the comparison,
        \item $\Phi\colon \sigma^{\ast}H_{\rig}\xrightarrow{\,\,\cong\,\,} H_{\rig}$ is an isomorphism of $A_K^\dag$-module,
        \item $\nabla\colon H_{\dR}\to \Omega_{S/K}^1\otimes H_{\dR}$ is an integrable connection
    that satisfies $\Phi\nabla=\nabla\Phi$,
        \item $\Fil^{\bullet}$ is a finite descending filtration on $H_{\dR}$ of locally free 
        $A_K$-module (i.e. each graded piece is locally free),
    that satisfies $\nabla(\Fil^i)\subset \Omega^1\ot \Fil^{i-1}$.
    \end{itemize}

Let $\mathrm{Fil}^\bullet$ denote the Hodge filtration on the de Rham cohomology,
and $\nabla$ the Gauss-Manin connection.
Then one has an object
\begin{equation}\label{F-isoc-eq1}
H^i(X/S)=H^i(X/A):=(H^i_\dR(X_K/S_K),H^i_\rig(X_\F/S_\F),c,\Phi,\nabla,\mathrm{Fil}^\bullet)
\end{equation}
in $\FilFMIC(S)$. 
For an integer $r$, the Tate object $\O_S(r)$ or $A(r)\in \FilFMIC(S)$ 
is defined in a customary way (loc.cit.), and
we define the Tate twist $M(r)$ by 
\[
M(r)=M\ot \O_S(r)=(M_\dR,M_\rig,\nabla,p^{-r}\Phi,\Fil^{\bullet+r})\in \FilFMIC(S).
\]

\subsection{Tate curve}
Let $\cE=\cE_{W[[\lambda]]}$ be the Tate curve as in {\bf A3}.
As is well-known, there is a uniformization
\[
\rho:\wh\bG_m:=\operatorname{Spec}\varprojlim_n \left(W/p^n[u,u^{-1}]\right)\lra \cE
\]
over $W[[\lambda]]$. Let $\wh R:=W((\lambda))^\wedge$ be the $p$-adic completion of
the ring of Laurant power series $W((\lambda))=W[[\lambda]][\lambda^{-1}]$.
Write $\wh R_K:=\wh R\ot_WK$.
The de Rham cohomology $H^1_\dR(\cE/\wh R)$ 
is a free $\wh R$-module of rank $2$.
Following \cite[\S 5.1]{AM}, we define a basis (which we call the {\it de Rham symplectic basis})
\[
\wh\omega,\, \wh\eta\in H^1_\dR(\cE/\wh R)=H^1_\dR(\cE\times_{W[[\lambda]]}\wh R/\wh R)
\]
in the following way.
Firstly $\wh\omega\in \vg(\cE,\Omega^1_{\cE/\wh R})$ is a differential form of
first kind such that $\rho^*\wh\omega=du/u$. Note that this depends on the uniformization
$\rho$.
Let $\nabla:H^1_\dR(\cE/\wh R_K)\to\Omega^1_{\wh R_K/K}\ot 
H^1_\dR(\cE/\wh R_K)$ be the
Gauss-Manin connection, where $\Omega^1_{\wh R_K/K}$ denotes the module of
continuous 1-forms.
Since $\ker(\nabla)\cong K$, one can define the unique element $\wh\eta$
such that $\nabla(\wh\eta)=0$ and $\Tr(\wh\omega\cup\wh\eta)=1$, where 
$\Tr:H^2_\dR(\cE/\wh R_K)\to \wh R_K$ is the trace map.

We shall give an explicit description of $\{\wh\omega,\wh\eta\}$ in case that $p\geq 5$ 
and $\cE$ is defined by $y^2=4x^3-g_2x-g_3$ 
(Proposition \ref{A-sympl-2}, \eqref{S-B-formula-1}).
\begin{prop}\label{local-GM}
\[
\begin{pmatrix}
\nabla(\wh\omega)&\nabla(\wh\eta)
\end{pmatrix}=
\begin{pmatrix}
\wh\omega&\wh\eta
\end{pmatrix}
\begin{pmatrix}
0&0\\
\frac{dq}{q}&0
\end{pmatrix}
\]
\end{prop}
\begin{pf}
\cite[Proposition 5.1]{AM}.
\end{pf}
Let $H^1(X/A)\in \FilFMIC(A_K)$ be the object in \eqref{F-isoc-eq1}.
Suppose that $\sigma$ extends to a $p$-th Frobenius on $\wh R=W((\lambda))^\wedge$ 
(e.g. $\sigma(\lambda)=c\lambda^p$ with $c\in 1+pW$).
Then, according to \cite[5.1]{AM}, it gives rise to an object
\[
H^1(\cE/\wh R)=(H^1_\dR(\cE/\wh R)\ot_WK,\nabla,\Phi)
\]
in the category $\FMIC(\wh R_K)=\FMIC(\wh R_K,\sigma)$.
\begin{prop}\label{local-Phi}
Define
\[
\log^{(\sigma)}(x):=p^{-1}\log\left(\frac{x^p}{x^\sigma}\right)=-p^{-1}
\sum_{n=1}^\infty\frac{1}{n}\left(1-\frac{x^p}{x^\sigma}\right)^n
\]
for $x\in \wh R^\times$. If $x\in W^\times$, then we also write
$\log^{(F)}(x)=\log^{(\sigma)}(x)$.
Put
\begin{equation}\label{tau-eq1}
\tau^{(\sigma)}(\lambda)
:=\log^{(\sigma)}(q)=\log^{(F)}(\kappa)+n\log^\s(\lambda)+\tau(\lambda)-p^{-1}\tau(\lambda)^\sigma\in W[[\lambda]].
\end{equation}
Then we have
\[
\begin{pmatrix}
\Phi(\wh\omega)&\Phi(\wh\eta)
\end{pmatrix}=
\begin{pmatrix}
\wh\omega&\wh\eta
\end{pmatrix}
\begin{pmatrix}
p&0\\
-p\tau^{(\sigma)}(\lambda)&1
\end{pmatrix}.
\]
\end{prop}
\begin{pf}
\cite[Theorem 5.6]{AM}.
\end{pf}

Let $\omega,\eta\in H^1_\dR(X/A)$ be a free $A$-basis such that
$\omega\in \vg(X,\Omega^1_{X/A})\cong A$.
Let $\wh\omega,\wh\eta\in H^1_\dR(\cE/\wh R)$ be the de Rham symplectic
basis as before.
Define $F(\lambda)\in \wh R$ and $f_i(\lambda)\in \wh R$ by 
\begin{equation}\label{F-f1-f2}
\omega=F(\lambda)\wh\omega,\quad \eta=f_1(\lambda)\wh\omega+f_2(\lambda)\wh\eta.
\end{equation}
In case that $p\geq 5$ and $\cE$ is defined by $y^2=4x^3-g_2x-g_3$,
we shall give an explicit description of $F(\lambda)$ and $f_i(\lambda)$ (Proposition \ref{A-sympl-2}
\eqref{S-B-formula-1} together with \eqref{S-B-formula}).
\begin{prop}\label{Frob-formula}
Let $\sigma$ be a $p$-th Frobenius on $A^\dag$ such that it extends on $\wh R$.
Let
\begin{align*}
F_{11}&=\frac{F(\lambda)^\sigma}{F(\lambda)}\left(1+\tau^\s(\lambda)\frac{f_1(\lambda)}{f_2(\lambda)}\right)\\
F_{12}&=\frac{1}{F(\lambda)}\left(pf_1(\lambda)^\sigma-
f_1(\lambda)\frac{f_2(\lambda)^\sigma-p\tau^\s(\lambda)f_1(\lambda)^\sigma}{f_2(\lambda)}\right)\\
F_{21}&=-\tau^\s(\lambda)\frac{F(\lambda)^\sigma}{f_2(\lambda)}\\
F_{22}&=\frac{f_2(\lambda)^\sigma-p\tau^\s(\lambda)f_1(\lambda)^\sigma}{f_2(\lambda)}.
\end{align*}
Then we have
\[
\begin{pmatrix}
\Phi(\omega)&\Phi(\eta)
\end{pmatrix}=
\begin{pmatrix}
\omega&\eta
\end{pmatrix}
\begin{pmatrix}
pF_{11}&F_{12}\\
pF_{21}&F_{22}
\end{pmatrix}.
\]
Moreover $F_{ij}\in A^\dag$.
\end{prop}
\begin{pf}
The representation matrix of $\Phi$ is straightforward from Proposition \ref{local-Phi} together with \eqref{F-f1-f2}. Then this implies that $F_{ij}\in A^\dag_K=A^\dag\ot_WK$.
We further show $F_{ij}\in A^\dag$.
By a commutative diagram
\[
\xymatrix{
0\ar[r]&A^\dag\ar[r]^p\ar[d]&A^\dag\ar[r]\ar[d]&A/pA\ar[r]\ar[d]^\cap&0\\
0\ar[r]&\wh R\ar[r]^p&\wh R\ar[r]&\F((\lambda))\ar[r]&0
}
\]
the $W$-module $\wh R/A^\dag$ is $p$-torsion free, and 
this implies that $A^\dag=A^\dag_K\cap \wh R$.
Therefore it is enough to show 
$F_{ij}\in \wh R$.
To do this, it is enough to show that $F(\lambda),f_2(\lambda)\in \wh R^\times$.
However, since $\{\omega,\eta\}$ is a free $A$-basis of $H_\dR^1(X/A)$,
one has that $\Tr(\omega\cup\eta)=F(\lambda)f_2(\lambda)\in A^\times$ as required.
\end{pf}
\begin{cor}\label{Frob-alpha}
Let $\cano:\operatorname{Spec} W\to S$ be a $W$-rational point of $S$, and ${\frak m}\subset A$
the corresponding prime ideal.
Let
$X_\cano$ be
the fiber at $\cano$, which is an elliptic curve over $W$.
Put $X_{\cano,\ol\F_p}:=X_\cano\times_W\ol\F_p$. 
Suppose $\sigma^{-1}({\frak m}A^\dag)\supset{\frak m}A^\dag$.
Then the representation 
matrix of the $p$-th Frobenius $\Phi$ on $H^1_\rig(X_{\cano,\ol\F_p}/\ol\F_p)$ is 
\[
\begin{pmatrix}
pF_{11}(\cano)&F_{12}(\cano)\\
pF_{21}(\cano)&F_{22}(\cano)
\end{pmatrix},
\]
where for $f\in A^\dag_K$ we write by 
$f(\cano)\in A^\dag_K/{\frak m}A^\dag_K\cong K$
the reduction modulo ${\frak m}A^\dag_K$.
\end{cor}

\subsection{Regulator formula for $K_2$ of elliptic curves}
\label{global-1ext-sect}
Let $U\subset X$ be an open set such that $Z:=X\setminus U$ is smooth over $W$.
Let $\partial:K^M_2(\O(U))\to K_1'(Z)$ be the boundary map (tame symbol).
We denote by $K^M_2(\O(U))_{\partial=0}$ the kernel of the boundary map.
Let
\[\xi\in K_2^M(\O(U))_{\partial=0}\]
be an element of Milnor's $K_2$. 
Let $G_\xi(\lambda)\in \wh R$ be defined by
\[
\dlog(\xi)=G_\xi(\lambda)\frac{d\lambda}{\lambda}\wedge\wh\omega\in \vg(\cE,\Omega^2_{\cE/W}).
\]
Since this has at worst log pole and $\rho^*\wh\omega=du/u$, we have
\begin{equation}\label{local-reg-eq0}
G_\xi(\lambda)=l+a_1\lambda+a_2\lambda^2+\cdots\in W[[\lambda]].
\end{equation}
The constant term $l$ is an integer (possibly zero).
Indeed, let $\delta$ be the composition of the maps
\[
\xymatrix{
K_2(\cE)\ar[r]^{\text{tame}}& K'_1(D)\ar[r]^{j^*}& K_1(\bG_{m,W})
\ar[r]& K_1(\bG_{m,K})\cong K[u,u^{-1}]^\times
\ar[r]^{\hspace{2cm}\ord_u}& \Z
},
\]
where $D\subset \cE$ is the central fiber at $\lambda=0$ which is the N\'eron $n$-gon over $W$,
and $j:\bG_{m,W}\hra D$
is the open immersion. Then $l=\pm \delta(\xi)$.
According to \cite[Proposition 4.2]{AM}, one can associate a 1-extension
\[
0\lra H^1(X/S)(2)\lra M_\xi(X/ S)\lra \O_S\lra 0
\]
in $\FilFMIC(S,\sigma)$ to the element to $\xi$.
Let $e_\xi\in \mathrm{Fil}^0M_\xi(X_K/A_K)_\dR$ be the unique lifting of $1\in A_K$.
Define $\ve^\s_i(\lambda)\in A^\dag_K$ by
\begin{equation}\label{local-reg-eq1}
\Phi(e_\xi)-e_\xi=\ve^\s_1(\lambda)\omega+\ve^\s_2(\lambda)\eta.
\end{equation}
\begin{thm}\label{syn-thm}
Let $\cano:\operatorname{Spec} W\to \operatorname{Spec} A$ be a $W$-rational point of $S$, and 
${\frak m}\subset A$ the corresponding prime ideal.
Let $X_\cano$ be the fiber at $\cano$.
Suppose $\sigma^{-1}({\frak m}A^\dag)\supset{\frak m}A^\dag$.
Let
\[
\reg_\syn:K_2(X_\cano)\lra H^2_\syn(X_\cano,\Q_p(2))\cong H^1_\dR(X_\cano/K)
\]
be the syntomic regulator map. Then
\[
\reg_\syn(\xi|_{X_\cano})=\ve_1^{(\sigma)}(\cano)\omega
+\ve_2^{(\sigma)}(\cano)\eta.
\]
\end{thm}
\begin{pf}
This is immediate from \cite[Theorem 4.4]{AM}.
\end{pf}

Define $E_i^\s(\lambda)\in \wh R_K$ by
\begin{equation}\label{local-reg-eq1-E}
\ve^\s_1(\lambda)\omega+\ve^\s_2(\lambda)\eta=E_1^\s(\lambda)\wh\omega+E_2^\s(\lambda)\wh\eta
\end{equation}
or explicitly
\begin{align}
E^\s_1(\lambda)&=F(\lambda)\ve^\s_1(\lambda)+f_1(\lambda)\ve^\s_2(\lambda)\label{global-1ext-eq1}\\
E^\s_2(\lambda)&=f_2(\lambda)\ve^\s_2(\lambda)\label{global-1ext-eq2},
\end{align}
where $F(\lambda)$ and $f_i(\lambda)$ are as in \eqref{F-f1-f2}.
\begin{thm}\label{global-reg}
Let $\sigma$ be a $p$-th Frobenius on $A^\dag$ such that it extends on $\wh R$.
Define $f'(\lambda)$ by $f'(\lambda)d\lambda=df(\lambda)$. Then we have
\begin{align}
(E^\s_1(\lambda))'&=\frac{G_\xi(\lambda)}{\lambda}-p^{-1}G_\xi(\lambda)^\sigma\frac{(\lambda^\sigma)'}{\lambda^\sigma}
\label{global-reg-eq2}\\
(E_2^\s(\lambda))'&=-E_1^\s(\lambda)\frac{q'}{q}+p^{-1}G_\xi(\lambda)^\sigma
\frac{(\lambda^\sigma)'}{\lambda^\sigma}\log^{(\sigma)}(q).
\label{global-reg-eq3}
\end{align}
If $\lambda^\sigma=c\lambda^p$ for some $c\in 1+pW$, then
$E^\s_1(\lambda), E^\s_2(\lambda)\in W[[\lambda]]\ot_W K$
\footnote{More precisely, we shall see that $E_i^\s(\lambda)\in W[[\lambda]]$ by 
Lemma \ref{int-lem} below.} and
\begin{equation}\label{global-reg-eq4}
E_1^\s(0)=\frac{l}{n}\log^{(F)}(\kappa)-lp^{-1}\log(c).
\end{equation}
\end{thm}
\begin{pf}
Apply $\nabla$ on \eqref{local-reg-eq1} and \eqref{local-reg-eq1-E}. Since $\Phi\nabla=\nabla\Phi$, one has
\[
G_\xi(\lambda)\frac{d\lambda}{\lambda}\wedge\wh\omega-
\Phi\left(G_\xi(\lambda)\frac{d\lambda}{\lambda}\wedge\wh\omega\right)
=(E^\s_1(\lambda))'d\lambda\wedge\wh\omega+\left(E^\s_1(\lambda)\frac{dq}{q}+
(E^\s_2(\lambda))'d\lambda\right)\wedge\wh\eta ,
\]
where we apply Proposition \ref{local-GM} to the right hand side of \eqref{local-reg-eq1-E}.
Apply Proposition \ref{local-Phi} to the second term in the left hand side.
Then we have
\[
\frac{G_\xi(\lambda)}{\lambda}\wh\omega-G_\xi(\lambda)^\sigma\frac{ (\lambda^\sigma)'}{\lambda^\sigma}
(p^{-1}\wh\omega-p^{-1}\log^{(\sigma)}(q)\wh\eta)
=(E^\s_1(\lambda))'\wh\omega+\left(E^\s_1(\lambda)\frac{q'}{q}+(E^\s_2(\lambda))'\right)\wh\eta.
\]
This finishes the proof of \eqref{global-reg-eq2} and \eqref{global-reg-eq3}.
Suppose that $\lambda^\sigma=c\lambda^p$ with $c\in 1+pW$.
It is straightforward to see $E_1^\s(\lambda)\in W[[\lambda]]\ot_W K$ from \eqref{local-reg-eq0}
and \eqref{global-reg-eq2}.
By \eqref{tau-eq1},
the coefficient of $\lambda^{-1}$ in the right hand side of \eqref{global-reg-eq3} is
\[
-nE^\s_1(0)+p^{-1}\cdot lp(\log^{(F)}(\kappa)-np^{-1}\log(c)).
\]
Since $E^\s_2(\lambda)\in \wh R_K$, this must be zero. 
This shows \eqref{global-reg-eq4}.
The right hand side of \eqref{global-reg-eq3} belongs to $K[[\lambda]]$, 
so that we have $E^\s_2(\lambda)\in \wh R_K\cap K[[\lambda]]
=W[[\lambda]]\ot_W K$.
\end{pf}
We have a complete description of the constant $E_1^{(\sigma)}(0)$, while
the constant $E_2^{(\sigma)}(0)$ seems more delicate.
\begin{conj}\label{const-conj}
Suppose $\lambda^\sigma=c\lambda^p$ for some $c\in 1+pW$.
Let $D\subset \cE$ be the central fiber at $\lambda=0$ which is the N\'eron $n$-gon over $W$, 
and $D_\bullet\to D$ the simplicial scheme.
Then there exists a map (specialization map)
\[
\mathrm{Sp}:H^2_\cM(\cE\times_{W[[\lambda]]}W((\lambda)),\Q(2))\lra H^2_\cM(D_\bullet,\Q(2))
\]
which is induced by the ``nearby cycle functor with respect to a parameter
$\kappa^{\frac{1}{n}}\lambda$'', and we have
\[
E^\s_2(0)=-\frac{l}{2n}(\log^{(F)}(\kappa)-np^{-1}\log(c))^2+R_\syn(\mathrm{Sp}(\xi)),
\]
where $R_\syn$ is the composition of the following arrows
\[
\xymatrix{
H^2_\cM(D_\bullet,\Q(2))\ar[r]^{\reg_\syn}&H^2_\syn(D_\bullet,\Q_p(2))&
H^1_\dR(D_\bullet/W)\ar[l]_{\quad\cong}
\ar[r]^{\hspace{0.8cm}\wh\eta\mapsto 1}_{\hspace{0.8cm}\cong}
&W.
}
\]
In particular, $R_\syn(\mathrm{Sp}(\xi))$ is a $\Q$-linear combination of 
$\ln_2^{(p)}(\epsilon)$'s with $\epsilon\in W\cap \ol\Q$, where $\ln^{(p)}_r(z)
=\sum_{p \not{\hspace{0.7mm}|}\, n}z^n/n^r$ is the $p$-adic dilogarithm function.
\end{conj}
\subsection{Unit root fomula}
\begin{lem}\label{Put-lem}
Let $\sigma$, $\Phi$ and $F_{ij}$ be as in Proposition \ref{Frob-formula}.
Put
\[
\wh B:=\varprojlim_n(A^\dag/p^nA^\dag[F_{22}^{-1}])=\varprojlim_n(A/p^nA[F_{22}^{-1}])
\]
and $\widehat{B}_K:=\widehat{B}\ot_WK$.
Note that $\wh B\ne\{0\}$ as $F_{22}\equiv f_2(\lambda)^{p-1}\not\equiv 0$ mod $p\wh R$
\footnote{
Indeed $F(\lambda)f_2(\lambda)\in A^\times$ as is in shown in the proof of Proposition \ref{Frob-formula}. Note that $f_2(\lambda)^{p-1}\in \wh R/p\wh R=k((\lambda))$ belongs to $A/pA$ as so does
($F_{22}$ mod $pA^\dag$). 
However this does not imply $f_2(\lambda)^{p-1}\in (A/pA)^\times$
because 
$F(\lambda)^{p-1}\not\in A/pA$.}.
Then
\[
\frac{F(\lambda)}{F(\lambda)^\sigma},\,\frac{E^\s_1(\lambda)}{F(\lambda)},\,
\frac{f_1(\lambda)}{F(\lambda)}\in \widehat{B}_K.
\]
\end{lem}
\begin{pf}
We prove our lemma by the same argument as that of \cite[Proposition (7.12)]{Put}.
For $f(\lambda)\in \wh B$, one has 
\[
\Phi(f(\lambda)\omega+\eta)=(F_{12}+pf(\lambda)^\sigma F_{11})\omega
+(F_{22}+pf(\lambda)^\sigma F_{21})\eta
\]
by Proposition \ref{Frob-formula}.
Let $V:=\{f(\lambda)\omega+\eta\mid f(\lambda)\in \widehat{B}\}$.
Consider a map
\[
S:V\lra V,\quad 
S(f(\lambda)\omega+\eta)=\frac{F_{12}+pf(\lambda)^\sigma F_{11}}{F_{22}+pf(\lambda)^\sigma F_{21}}\omega+\eta.
\]
This is well-defined as $F_{22}+pf(\lambda)^\sigma F_{21}\in \wh B^\times$ and we have
\[
S(f_1(\lambda)\omega+\eta)-S(f_2(\lambda)\omega+\eta)
=\frac{p(f_1(\lambda)^\sigma-f_2(\lambda)^\sigma)D}
{(F_{22}+pf_1(\lambda)^\sigma F_{21})(F_{22}+pf_2(\lambda)^\sigma F_{21})}\omega,
\]
where $D:=F_{11}F_{22}-F_{21}F_{12}\in (A^\dag)^\times$.
This shows that
$S$ is a contraction map with respect to the norm $||\cdot||$
on $\wh R=W((\lambda))^\wedge$ which is defined by
$||\sum a_i\lambda^i||:=\sup\{|a_i|\mid i\in\Z\}$.
Hence there is a unique fixed point $u_0=f_0(\lambda)\omega+\eta\in V$, namely
\[
S(u_0)=u_0\quad \Longleftrightarrow\quad \Phi(u_0)=(F_{22}+pf_0(\lambda)^\sigma F_{21})u_0.
\]
Replace $\wh B$ with $\wh R$ and consider the similar map
$\wh S$ on $\wh V:=\{f(\lambda)\omega+\eta\mid f(\lambda)\in \widehat{R}\}$. 
Then one also has a unique fixed point
$\wh u_0=\wh f_0(\lambda)\omega+\eta\in \wh V$ of $\wh S$, and then it turns out
that $u_0=\wh u_0$ by the uniqueness. 
From \eqref{F-f1-f2}, we have 
\[
f_2(\lambda)\wh \eta=-\frac{f_1(\lambda)}{F(\lambda)}\omega+\eta.
\]
By Proposition \ref{local-Phi}, one has that
$\Phi(f_2(\lambda)\wh \eta)=f_2(\lambda)^\sigma \wh\eta$. This means that the above
coincides with
the unique fixed point $\wh u_0$ of $\wh S$, and hence
\[
-\frac{f_1(\lambda)}{F(\lambda)}=\wh f_0(\lambda)=f_0(\lambda).
\]
Moreover, 
\[
\Phi(f_2(\lambda)\wh \eta)=f_2(\lambda)^\sigma \wh\eta
\quad\Longleftrightarrow\quad
\Phi(f_0(\lambda)\omega+\eta)
=\frac{f_2(\lambda)^\sigma}{f_2(\lambda)}(f_0(\lambda)\omega+\eta)
\]
implies 
\[
\quad \frac{f_2(\lambda)^\sigma}{f_2(\lambda)}=F_{22}+pf_0(\lambda)^\sigma F_{21}.
\]
We thus have $f_1(\lambda)/F(\lambda),f_2(\lambda)^\sigma/f_2(\lambda)\in\wh B$.
Moreover since $\Tr(\omega\cup\eta)=F(\lambda)f_2(\lambda)\in A^\times$ (see 
the proof of Proposition \ref{Frob-formula}), one also has 
$F(\lambda)/F(\lambda)^\sigma\in\wh B$.
Finally, we have 
\[
\Tr(E^\s_1(\lambda)\wh\omega+E^\s_2(\lambda)\wh\eta)\cup 
u_0)=\Tr(E^\s_1(\lambda)\wh\omega\cup u_0)=E^\s_1(\lambda)/F(\lambda) \cdot
\Tr(\omega\cup \eta)\in \wh B_K
\]
from \eqref{local-reg-eq1-E}.
This implies $E^\s_1(\lambda)/F(\lambda)\in \widehat{B}_K$.
\end{pf}
\begin{cor}[Unit root formula]\label{unit-root-formual}
Suppose $W=\Z_p$ for simplicity. 
Let $\cano:\operatorname{Spec} \Z_p\to\operatorname{Spec} A$ be a $\Z_p$-rational point of $S$, and 
${\frak m}\subset A$ the prime ideal.
Let $X_\cano$ be
the fiber at $\cano$.
Suppose $\sigma^{-1}({\frak m}A^\dag)\supset{\frak m}A^\dag$.
Suppose further that 
$X_{\cano,\F_p}:=X_\cano\times_{\Z_p} \F_p$ is an ordinary elliptic curve. 
Then the unit root of $X_{\cano,\F_p}$ is
\[
\frac{F(\lambda)}{F(\lambda)^\sigma}\bigg|_{\lambda=\cano}.
\]
\end{cor}
\begin{pf}
We first note that $X_{\cano,\F_p}$ is ordinary if and only if 
the absolute Frobenius on $H^1(\O_{X_{\cano,\F_p}})$ does not vanish, and this is equivalent to that $F_{22}(\lambda)|_{\lambda=\cano}\not\equiv 0$ mod $p$ in Proposition \ref{Frob-formula}.
Hence $\cano$ is a point of $\operatorname{Spec} \wh B_K$ and the special value $F(\lambda)/F(\lambda)^\sigma|_{\lambda=\cano}$ makes sense.
As is shown in the proof of Lemma \ref{Put-lem}, the unit root vector is $u_0|_{\lambda=\cano}$ and 
the unit root is
\[
\frac{f_2(\lambda)^\sigma}{f_2(\lambda)}\bigg|_{\lambda=\cano}
=\frac{F(\lambda)}{F(\lambda)^\sigma}\bigg|_{\lambda=\cano}\times
\frac{d(\lambda)^\sigma}{d(\lambda)}\bigg|_{\lambda=\cano}
=\frac{F(\lambda)}{F(\lambda)^\sigma}\bigg|_{\lambda=\cano},
\]
where $d(\lambda):=F(\lambda)f_2(\lambda)\in A^\times$.
\end{pf}

\subsection{Bounds on $p$-adic expansions of $\ve_k^{(\sigma)}(\lambda)$}\label{exp-sect}
In this section we give effective bounds of rational functions
$\ve_k^{(\sigma)}(\lambda)$ modulo $p^n$ according to \cite{KT}.
The main results are Theorems \ref{exp-cor} and \ref{exp-Frob},
which play key roles in later sections.
We note that Kedlaya and Tuitman give a theorem 
in more general situation (\cite[Theorem 2.1]{KT}).
However, as we need a little more precise statement, 
we give a self-contained proof for the sake of completeness.

\medskip

We assume $p\geq3$.
Let $f:Y\to C$ be as in \S \ref{fib-set-sect} which satisfies {\bf A1}, {\bf A2} and {\bf A3}.
For a $W$-rational point $P_i\in Z$, let $\lambda_i$ be a uniformizer of the local ring at $P_i$.
Put $Y_i:=Y\times_{C}\operatorname{Spec} W[[\lambda_i]]$ and $D_i:=f^{-1}(P_i)$ the central fiber which is
a relative NCD over $W$ by {\bf A1}. Put $X_i=Y_i\setminus D_i$.
We consider the following condition.
\begin{description}
\item[A4]
The multiplicity of each component of $D_i$ is prime to $p$. Moreover
put \[H_i:=\operatorname{Im}[H^1_\zar(Y_i,\Omega^\bullet_{Y_i/W[[\lambda_i]]}(\log D_i))\to
H^1_\dR(X_i/W((\lambda_i)))]\]
and
\[\Fil^1 H_i:=\operatorname{Im}[H^0_\zar(Y_i,\Omega^1_{Y_i/W[[\lambda_i]]}(\log D_i))\to
H^1_\dR(X_i/W((\lambda_i)))].\]
Then, both of $\Fil^1H_i$ and 
$H_i/\Fil^1H_i$ are free $W[[\lambda_i]]$-modules of rank $1$.
\end{description} 
We endow the log structures on $Y_i$ and $\operatorname{Spec} W[[\lambda_i]]$ arising from the divisors
$D_i$ and $\lambda_i=0$ respectively. 
Since the multiplicity of each component of $D_i$ is prime to $p$,
 the morphism $(Y_i,D_i)\to (\operatorname{Spec} W[[\lambda_i]],(\lambda_i))$ is a smooth morphism
of log schemes in the sense of
\cite[\S 6]{Kato-log-crys}. 
By \cite[Theorem (6.4)]{Kato-log-crys}, there is the natural isomorphism
\[
H^k_\crys((Y_{i,\ol\F_p},D_{i,\ol\F_p})/(W[[\lambda_i]],(\lambda_i)))
\cong H^k_\zar(Y_i,\Omega^\bullet_{Y_i/W[[\lambda_i]]}(\log D_i)),
\]
where the left hand side is the log crystalline cohomology.
In particular, 
 if $\sigma(\lambda_i)=c\lambda_i^p$ with $c\in 1+pW$, then the $\sigma$-linear Frobenius
$\Phi$ is defined.
\begin{lem}
\label{int-lem}
Assume that {\bf A4} is satisfied.
Let $\{\omega_i,\eta_i\}$ be a $W[[\lambda_i]]$-basis of $H_i$.
Suppose that $\sigma$ is defined by $\lambda_i^\sigma=c\lambda_i^p$ for some $c\in 1+pW$.
Define $\ve_j^\s(\lambda_i)\in A^\dag_K$ ($j=1,2$) by
\begin{equation}\label{int-lem-eq1}
\Phi(e_\xi)-e_\xi=\ve^\s_1(\lambda_i)\omega_i+\ve^\s_2(\lambda_i)\eta_i\in
H^1_\dR(X_{i}/W((\lambda_i)))\ot_{W((\lambda_i))}W((\lambda_i))^\wedge ,
\end{equation}
where $W((\lambda_i))^\wedge$ denotes the $p$-adic completion.
Then this element lies in the image of
$H^1_\zar(Y_i,\Omega^\bullet_{Y_i/W[[\lambda_i]]}(\log D_i))$.
Hence
$\ve_1^{(\sigma)}(\lambda_i)$ and $\ve_2^{(\sigma)}(\lambda_i)$ belong to the subspace 
$W[[\lambda_i]]\cap A_K^\dag=W[[\lambda_i]]\cap A^\dag$.
\end{lem}
\begin{pf}
\cite[Theorem 4.4]{AM}.
\end{pf}
\begin{lem}[Transformation of Frobenius]\label{F-lem}
Let $\sigma$ be a $p$-th Frobenius on $A^\dag$.
Let $\FMIC(A_K^\dag,\sigma)$ be the category of objects $(H,\nabla,\Phi)$, where 
$(H,\nabla)$ is an integrable connection of free $A^\dag_K$-module, and $\Phi$ is 
a $p$-th Frobenius on $H$ compatible with $\sigma$.
Let $\partial=\nabla_{d/d\lambda_i}$ be the differential operator which is defined to be
the composition
\[
\xymatrix{
H\ar[r]^-\nabla&
A\cdot d\lambda_i\ot H\ar[rr]^-{d\lambda_i\ot\omega\mapsto\omega}&&H.
}\]
Let $\sigma'$ be another Frobenius on $A^\dag$.
Thanks to the natural equivalence $\FMIC(A^\dag_K,\sigma)\cong \FMIC(A^\dag_K,\sigma')$
one can associate $(H,\nabla,\Phi')\in \FMIC(A^\dag_K,\sigma')$ to an object
$(H,\nabla,\Phi)\in \FMIC(A^\dag_K,\sigma)$.
Then, for $x\in H$ we have
\[
\Phi'(x)-\Phi(x)=\sum_{k=1}^\infty\frac{(\lambda_i^{\sigma'}-\lambda_i^\sigma)^k}{k!}\Phi\partial^kx.
\]
\end{lem}
\begin{pf}
See \cite[6.1]{EK} or \cite[17.3.1]{Ke}.
\end{pf}
\begin{lem}
\label{est-lem}
Assume that {\bf A4} is satisfied.
Let $\{\omega_i,\eta_i\}$ be a $W[[\lambda_i]]$-basis
of $H_i$
such that $\omega_i\in \Fil^1H_i$.
Let $\sigma$ be a $p$-th Frobenius on $A^\dag$.
Assume that $\lambda_i^\sigma\in pW[[\lambda_i]]+\lambda_i
W[[\lambda_i]]$ so that it naturally induces a $p$-th Frobenius
on $W[[\lambda_i]]$.
Let $\ve^\s_1(\lambda_i)$ and $\ve^\s_2(\lambda_i)$ be defined by \eqref{int-lem-eq1}.
Put
\begin{equation}\label{exp-lem-eq1}
e_n:=\max\{i\in\Z_{\geq 0}\mid \ord_p(p^i/i!)<n\}\footnote{Notice that
$e_n<\infty$ as we assumed $p>2$.}.
\end{equation}
Put $l_i:=\min\{\ord_{\lambda_i}(\lambda^\sigma_i-\lambda_i^p),p\}\in \Z_{\geq0}$, where $\ord_{\lambda_i}:K((\lambda_i))\to \Z$ is the
order function with respect to $\lambda_i$.
Then
\[
\ve^\s_1(\lambda_i),\,\ve^\s_2(\lambda_i)\in \lambda_i^{(l_i-p)e_{n+2}}W[[\lambda_i]]+p^nW[[\lambda_i,\lambda_i^{-1}]]
\]
for any $n\geq 1$.
\end{lem}
\begin{pf}
Let $\sigma_0$ be the $p$-th Frobenius on $W[[\lambda_i]]$ such that $\sigma_0(\lambda_i)=\lambda_i^p$, and $\Phi_0$ the Frobenius with respect to $\sigma_0$.
By Lemma \ref{F-lem}, we have
\begin{equation}\label{exp-lem-eq2}
\Phi(e_\xi)-\Phi_0(e_\xi)=\sum_{k=1}^\infty\frac{(\lambda_i^{\sigma}-\lambda_i^p)^k}{k!}\Phi_0\partial^ke_\xi.
\end{equation}
Write the right hand side by $G_1(\lambda_i)\omega_i+G_2(\lambda_i)\eta_i$, so that we have
\[
\ve^{(\sigma_0)}_1(\lambda_i)-\ve^{(\sigma)}_1(\lambda_i)=G_1(\lambda_i),
\quad
\ve^{(\sigma_0)}_2(\lambda_i)-\ve^{(\sigma)}_2(\lambda_i)=G_2(\lambda_i).
\]
It follows from Lemma \ref{int-lem} that $\ve_1^{(\sigma_0)}(\lambda_i)$
and $\ve_2^{(\sigma_0)}(\lambda_i)$ belong to $W[[\lambda_i]]$.
Hence it is enough to show
\begin{equation}\label{exp-lem-eq3}
G_1(\lambda_i),\, G_2(\lambda_i)\in \lambda_i^{(l_i-p)e_n}W[[\lambda_i]]+p^nW[[\lambda_i,\lambda_i^{-1}]].
\end{equation}
The differential operator $\partial$ is induced from the connection
\[
\nabla:H^1(Y_i,\Omega^\bullet_{Y_i/W[[\lambda_i]]}(\log D_i))\lra \frac{d\lambda_i}{\lambda_i}
\ot H^1(Y_i,\Omega^\bullet_{Y_i/W[[\lambda_i]]}(\log D_i)).
\]
Moreover since $\nabla e_\xi=\dlog(\xi)\in \frac{d\lambda_i}{\lambda_i}\ot
\vg(Y_i,\Omega^1_{Y_i/W[[\lambda_i]]}(\log D_i))$, 
we have
\[
\partial e_\xi \in\lambda_i^{-1}W[[\lambda_i]]\omega_i,\quad
\partial^k e_\xi \in\lambda_i^{-k}(W[[\lambda_i]]\omega_i+W[[\lambda_i]]\eta_i)
\]
for all $k\geq 2$.
Since $\Phi_0$ acts on the cohomology group
$H^1(Y_i,\Omega^\bullet_{Y_i/W[[\lambda_i]]}(\log D_i))$, 
there are $u_k(\lambda_i),v_k(\lambda_i)\in W[[\lambda_i]]$ such that
\begin{equation}\label{exp-lem-eq4}
\Phi_0\partial^k e_\xi=\lambda_i^{-kp}p^{-2}(u_k(\lambda_i)\omega_i+v_k(\lambda_i)\eta_i)
\end{equation}
for each $k\in \Z_{\geq 1}$,
where ``$p^{-2}$'' comes from the Tate twist $H^1(X/A)\ot A(2)$.
Notice that since $\Phi_0\omega_i\equiv 0$ modulo $p$, one has
\[
\Phi_0\partial e_\xi
\in \lambda_i^{-p}p^{-2}W[[\lambda_i]]\Phi_0\omega_i\subset
\lambda_i^{-p}p^{-1}(W[[\lambda_i]]\omega_i+W[[\lambda_i]]\eta_i)
\]
and hence \[
u_1(\lambda_i),\,v_1(\lambda_i)\in pW[[\lambda_i]].
\]
Applying \eqref{exp-lem-eq4} to \eqref{exp-lem-eq2}, we have
\begin{align*}
G_1(\lambda_i)\omega_i+G_2(\lambda_i)\eta_i
&=
\sum_{k=1}^\infty\frac{(\lambda_i^{\sigma}-\lambda_i^p)^k}{k!}
\lambda_i^{-kp}p^{-2}(u_k(\lambda_i)\omega_i+v_k(\lambda_i)\eta_i)\\
&=
\sum_{k=1}^\infty\frac{p^{k-2}}{k!}
\left(\frac{p^{-1}(\lambda^\sigma_i-\lambda_i^p)}{\lambda_i^p}\right)^k(u_k(\lambda_i)\omega_i+v_k(\lambda_i)\eta_i).
\end{align*}
Now \eqref{exp-lem-eq3} is immediate from this
on noticing $p^{-1}(\lambda^\sigma_i-\lambda_i^p)/\lambda_i^p\in \lambda_i^{l_i-p}W[[\lambda_i]]$.
\end{pf}
\begin{thm}\label{exp-cor}
Suppose that $C=\P^1_W$ and 
$f:Y\to \P^1$ satisfies {\bf A1},\ldots,{\bf A3} and {\bf A4} for all $i$.
Let $\lambda$ be an inhomogeneous coordinates of $\P^1$.
We assume that $f$ has a split multiplicative reduction at $\lambda=0$ as in {\bf A3}.
Let $\{\cano_0,\cano_1,\ldots,\cano_m\}\subset W$ 
be the subset such that $D_i:=f^{-1}(\cano_i)$ is a singular fiber, and set $\cano_0=0$.
Put $D_\infty:=f^{-1}(\infty)$ (possibly a smooth fiber).
Put $\lambda_i:=\lambda-\cano_i$ for $i\in \{0,1,\ldots,m\}$ and $\lambda_\infty:=\lambda^{-1}$.
Suppose further that the following condition is satisfied. 
\begin{description}
\item[\bf A5]
There are $\omega\in \vg(\Omega^1_{X/A})$ and $\eta\in H^1_\dR(X/A)$
such that for each $i\in\{0,1,\ldots,m,\infty\}$, $\{\lambda^{a_i}_i\omega,\lambda_i^{b_i}\eta\}$ 
is a $W[[\lambda_i]]$-basis of 
\[
H_i=\operatorname{Im}\left[H^1(Y_i,\Omega^\bullet_{Y_i/W[[\lambda_i]]}(\log D_i))\to H^1_\dR(X_i/W((\lambda_i)))\right]
\]
for some $a_i,b_i\in \Z$.
\end{description}
Let $\sigma$ be a $p$-th Frobenius on $A^\dag$ such that $\lambda^\sigma=c\lambda^p$
with $c\in 1+pW$.
Define $\ve^\s_k(\lambda)\in A^\dag_K$ by
\[
e_\xi-\Phi(e_\xi)=\ve^\s_1(\lambda)\omega+\ve^\s_2(\lambda)\eta.
\]
Put $l_i:=\min\{\ord_{\lambda_i}(\lambda_i^\sigma-\lambda_i^p),p\}$.
Let $e_n$ be as in \eqref{exp-lem-eq1}, and 
\begin{equation}\label{exp-cor-eq0}
d_{1,n}:=\sum_{i=0}^m(-a_i+(p-l_i)e_{n+2})-a_\infty,\quad
d_{2,n}:=\sum_{i=0}^m(-b_i+(p-l_i)e_{n+2})-b_\infty.
\end{equation}
Then $\ve_k^{(\sigma)}(\lambda)\in A^\dag$ and there are polynomials $f_{k,n}(\lambda)$ of degree $\leq d_{k,n}$ such that
\begin{equation}\label{exp-cor-eq1}
\ve_1^{(\sigma)}(\lambda)\equiv\frac{f_{1,n}(\lambda)}
{\prod_{i=0}^m\lambda_i^{-a_i+(p-l_i)e_{n+2}}}\mod p^nA^\dag,
\end{equation}
\begin{equation}\label{exp-cor-eq2}
\ve_2^{(\sigma)}(\lambda)\equiv\frac{f_{2,n}(\lambda)}
{\prod_{i=0}^m\lambda_i^{-b_i+(p-l_i)e_{n+2}}}\mod p^nA^\dag.
\end{equation}
\end{thm}
\begin{pf}
Put $\omega_i:=\lambda_i^{a_i}\omega$ and $\eta_i:=\lambda_i^{b_i}\eta$. 
Then
\[
\ve^\s_1(\lambda)\omega+\ve^\s_2(\lambda)\eta=\lambda_i^{-a_i}\ve^\s_1(\lambda)\omega_i+\lambda_i^{-b_i}\ve^\s_2(\lambda)\eta_i
\]
for $i\in \{0,1,\ldots,m,\infty\}$.
Now the assertion is immediate from Lemma \ref{est-lem} (estimate lemma).
\end{pf}
In a similar way to the above, we can obtain an estimate of 
the $p$-adic expansion of the Frobenius matrix on $H^1_\rig(X_\F/S_\F)$.
\begin{thm}\label{exp-Frob}
Let the notation and assumption be as in Theorem \ref{exp-cor}.
Let
\[
\begin{pmatrix}
\Phi(\omega)&\Phi(\eta)
\end{pmatrix}
=\begin{pmatrix}
\omega&\eta
\end{pmatrix}
\begin{pmatrix}
pF_{11}&F_{12}\\
pF_{21}&F_{22}
\end{pmatrix}.
\]
Let $e_n$ be as in \eqref{exp-lem-eq1}, and put
\begin{align*}
d_{11,n}:&=\sum_{i=0}^m(pa_i-a_i+(p-l_i)e_n)+pa_\infty-a_\infty,\\
d_{21,n}:&=\sum_{i=0}^m(pa_i-b_i+(p-l_i)e_n)+pa_\infty-b_\infty,\\
d_{12,n}:&=\sum_{i=0}^m(pb_i-a_i+(p-l_i)e_n)+pb_\infty-a_\infty,\\
d_{22,n}:&=\sum_{i=0}^m(pb_i-b_i+(p-l_i)e_n)+pb_\infty-b_\infty.
\end{align*}
Then there are polynomials $f_{ij,n}(\lambda)\in W[\lambda]$ of degree $\leq d_{ij,n}$ such that
\begin{align}
pF_{11}&\equiv \left(\prod_{i=0}^m(\lambda_i^\sigma)^{a_i}\lambda_i^{-a_i+(p-l_i)e_n}\right)^{-1}f_{11,n}(\lambda)
\mod p^n A^\dag\label{exp-Frob-eq4},\\
pF_{21}&\equiv \left(\prod_{i=0}^m(\lambda_i^\sigma)^{a_i}\lambda_i^{-b_i+(p-l_i)e_n}\right)^{-1}f_{21,n}(\lambda)
\mod p^n A^\dag\label{exp-Frob-eq5},\\
F_{12}&\equiv \left(\prod_{i=0}^m(\lambda_i^\sigma)^{b_i}\lambda_i^{-a_i+(p-l_i)e_n}\right)^{-1}f_{12,n}(\lambda)
\mod p^n A^\dag\label{exp-Frob-eq6},\\
F_{22}&\equiv \left(\prod_{i=0}^m(\lambda_i^\sigma)^{b_i}\lambda_i^{-b_i+(p-l_i)e_n}\right)^{-1}f_{22,n}(\lambda)
\mod p^n A^\dag\label{exp-Frob-eq7}.
\end{align}
\end{thm}
\begin{pf}
Put $\omega_i:=\lambda_i^{a_i}\omega$ and $\eta_i:=\lambda_i^{b_i}\eta$ for $i\in \{0,1,\ldots,m,\infty\}$. Let $F^{(i)}_{ab}$ be defined by
\[
pF^{(i)}_{11}\omega_i+pF^{(i)}_{21}\eta_i
:=\Phi(\omega_i)=(\lambda_i^\sigma)^{a_i}\lambda_i^{-a_i}pF_{11}\cdot\omega_i
+(\lambda_i^\sigma)^{a_i}\lambda_i^{-b_i}pF_{21}\cdot\eta_i,
\]
\[
F^{(i)}_{12}\omega_i+F^{(i)}_{22}\eta_i
:=\Phi(\eta_i)=(\lambda_i^\sigma)^{b_i}\lambda_i^{-a_i}F_{12}\cdot\omega_i
+(\lambda_i^\sigma)^{b_i}\lambda_i^{-b_i}F_{22}\cdot\eta_i.
\]
Let $\sigma_0(\lambda_i)=\lambda_i^p$, and $\Phi_0$ the associated Frobenius on $H^1_\rig(X_0/S_0)$.
Since $\Phi_0$ acts on $H^1(Y_i,\Omega^\bullet_{Y_i/W[[\lambda_i]]}(\log D_i))$, we have
the ``integrality''
\begin{equation}\label{exp-Frob-eq1}
\Phi_0(\omega_i),\,
\Phi_0(\eta_i)\in W[[\lambda_i]]\omega_i+W[[\lambda_i]]\eta_i.
\end{equation}
It follows from Lemma \ref{F-lem} (transformation of Frobenius) that
there are $u_k(\lambda_i),v_k(\lambda_i)\in W[[\lambda_i]]$ such that
\begin{align}
\Phi(\omega_i)-\Phi_0(\omega_i)
&=
\sum_{k=1}^\infty\frac{(\lambda_i^{\sigma}-\lambda_i^p)^k}{k!}
\lambda_i^{-kp}(u_k(\lambda_i)\omega_i+v_k(\lambda_i)\eta_i)\label{exp-Frob-eq2}\\
&=
\sum_{k=1}^\infty\frac{p^k}{k!}
\left(\frac{l(\lambda_i)}{\lambda_i^p}\right)^k(u_k(\lambda_i)\omega_i+v_k(\lambda_i)\eta_i)
\label{exp-Frob-eq3}
\end{align}
as $\partial^k\omega_i\in\lambda_i^{-k}W[[\lambda_i]]\omega_i+\lambda_i^{-k}W[[\lambda_i]]\eta_i$, where $l(\lambda_i):=p^{-1}(\lambda_i^\sigma-\lambda_i^p)\in W[\lambda_i]$. 
By \eqref{exp-Frob-eq1} and \eqref{exp-Frob-eq3}, we have ``estimates''
\[
pF^{(i)}_{11},\,pF^{(i)}_{21}\in \lambda_i^{(l_i-p)e_n}W[[\lambda_i]]+p^nW[[\lambda_i,\lambda^{-1}_i]]
\]
or equivalently
\begin{align*}
pF_{11}&\in \frac{1}{(\lambda_i^\sigma)^{a_i}\lambda_i^{-a_i+(p-l_i)e_n}}W[[\lambda_i]]+p^nW[[\lambda_i,\lambda^{-1}_i]]\\
pF_{21}&\in \frac{1}{(\lambda_i^\sigma)^{a_i}\lambda_i^{-b_i+(p-l_i)e_n}}W[[\lambda_i]]+p^nW[[\lambda_i,\lambda^{-1}_i]]
\end{align*}
for all $i\in \{0,1,\ldots,m,\infty\}$. Now
\eqref{exp-Frob-eq4} and \eqref{exp-Frob-eq5} follows from this.
In the same way, we have
\begin{align*}
F_{12}&\in \frac{1}{(\lambda_i^\sigma)^{b_i}\lambda_i^{-a_i+(p-l_i)e_n}}W[[\lambda_i]]+p^nW[[\lambda_i,\lambda^{-1}_i]]\\
F_{22}&\in \frac{1}{(\lambda_i^\sigma)^{b_i}\lambda_i^{-b_i+(p-l_i)e_n}}W[[\lambda_i]]+p^nW[[\lambda_i,\lambda^{-1}_i]]
\end{align*}
for all $i\in \{0,1,\ldots,m,\infty\}$, and hence \eqref{exp-Frob-eq6} 
and \eqref{exp-Frob-eq7} follow. This completes the proof.
\end{pf}
\subsection{Computing Syntomic Regulators}\label{CSyn-sect}
In this section, we demonstrate how to compute the syntomic regulators of $K_2$ of elliptic curves. The key is Theorem \ref{exp-cor}. 
In what follows, $p\geq 3$ and $W=W(\F_{p^r})$ is the Witt ring of the finite 
field $\F_{p^r}$ with $p^r$ elements. Put $K:=\operatorname{Frac} (W)$.

\medskip

\noindent{\bf (Situation)}.
Let $f:Y\to \P^1$ be an elliptic fibration whose generic fiber is defined by
an affine equation \[y^2=4x^3-g_2 x-g_3,\quad g_2,g_3\in K(\lambda),\]
where $\lambda$ is an inhomogeneous coordinate of $\P^1$. 
Put
\[
\omega=\frac{dx}{y},\quad \eta=\frac{xdx}{y}.
\]
We assume that $f$ satisfies
the conditions {\bf A1},\ldots,{\bf A3} in \S \ref{fib-set-sect}, and the multiplicity 
of each component of an arbitrary singular fiber is prime to $p$.
Then the conditions
{\bf A4} and {\bf A5} in \S \ref{exp-sect} are satisfied by Proposition \ref{A3}.
Let $\{\cano_0,\cano_1,\ldots,\cano_m\}\subset W$ 
be the subset such that $D_i:=f^{-1}(\cano_i)$ is a singular fiber
($f^{-1}(\infty)$ may be a singular fiber).
We suppose that $\cano_0=0$ and $f$ has a split multiplicative reduction at $\lambda=0$
as in {\bf A3}.
Put $S:=\operatorname{Spec} W[\lambda,\prod_i (\lambda-\cano_i)^{-1}]$ and
$X:=f^{-1}(S)$.

\medskip

Let $\xi\in K_2(X)$. 
Let $\sigma$ be a $p$-th Frobenius on 
$W[\lambda,\prod(\lambda-\cano_i)^{-1}]^\dag$
given by 
\begin{equation}\label{CSyn-eq00}
\sigma(\lambda)=c\lambda^p,\quad c\in 1+pW.
\end{equation} 
Let
\[
e_\xi-\Phi(e_\xi)=\ve_1^{(\sigma)}(\lambda)\omega+\ve_2^{(\sigma)}(\lambda)\eta.
\]
Let $s\geq 1$ be an integer, and
$\cano\in W$ an element such that $\cano\not\equiv \cano_i$ mod $p$ for all $i$.
Our goal is to compute
\begin{equation}\label{CSyn-eq0}
\ve_1^{(\sigma)}(\cano),\,
\ve_2^{(\sigma)}(\cano)\mod p^sW
\end{equation}
explicitly. In case $\sigma(\lambda)=\cano^{F-p}\lambda^p$, this gives the syntomic regulator
\[
\reg_\syn(\xi|_{X_\cano})=\ve_1^{(\sigma)}(\cano)\omega
+\ve_2^{(\sigma)}(\cano)\eta
\in H^2_\syn(X_\cano,\Q_p(2))\cong H^1_\dR(X_\cano/K)
\]
for the fiber $X_\cano:=f^{-1}(\cano)$ 
(see Theorem \ref{syn-thm}).

\medskip

\noindent{\bf (Notation)}.
Replacing the variables $x,y$ with $\lambda^{2k}x,\lambda^{3k}y$ if necessary,
we take an affine equation $y^2=4x^3-g_2x-g_3$ with
$\ord_\lambda(g_2)=\ord_\lambda(g_3)=0$ and $n:=\ord_\lambda(\Delta)>0$, where $\Delta:=g_2^3-27g_3^2$.
For $f(\lambda)\in K((\lambda))$, we denote by $f'(\lambda)=df(\lambda)/d\lambda$ the derivative with respect to $\lambda$.

\begin{enumerate}
\renewcommand{\labelenumi}{(\roman{enumi})}
\item
Put $\lambda_i:=\lambda-\cano_i$ for $i\in \{0,1,\ldots,m\}$ and $\lambda_\infty:=\lambda^{-1}$.
For $i\in \{0,1,\ldots,m\}$, put $l_i:=\min\{\ord_{\lambda_i}(\lambda_i^{\sigma}-\lambda_i^p),p\}$.
It follows from Proposition \ref{A3}
that for each $i\in \{0,1,\ldots,m,\infty\}$, there are integers $a_i,b_i$ such that 
\[
\lambda_i^{a_i}\omega,\quad \lambda_i^{b_i}\eta
\]
is a free $W[[\lambda_i]]$-basis of $H_i$, where $H_i$ is as in Theorem \ref{exp-cor}, {\bf A5}. 
\item(cf. \eqref{global-mult-period}).
Put $J:=g_2^3/\Delta$.
The $j$-invariant is $1728J$. 
Put
\[
\kappa:=\frac{\lambda^{-n}}{1728J}\bigg|_{\lambda=0}, \quad \text{where }n:=-\ord_\lambda(J),
\]
which belongs to $W^\times$ by {\bf A3}.
\item(cf. Proposition \ref{A-sympl-1}).
We fix a square root
\[
c_0:=\sqrt{-\frac{g_2(0)}{18g_3(0)}}.
\]
By {\bf A3}, $f$ has a split multiplicative reduction at $\lambda=0$ and the fiber
$f^{-1}(0)$ is a relative NCD over $W$, so that we have $c_0\in W^\times$. 
\item(cf. \eqref{F-f1-f2} and Proposition \ref{A-sympl-1}, \eqref{S-B-formula}).
Let $(12g_2)^{-1/4}\in W[[\lambda]]$ be the power series 
with the constant term $c_0$. Put
\[
F(\lambda):=(12g_2)^{-\frac{1}{4}}{}_2F_1\left({\frac{1}{12},\frac{5}{12}\atop 1};J^{-1}\right).\]
\item(cf. \eqref{tau} and Proposition \ref{A-sympl-3}).
Let $\tau(\lambda)\in \lambda K[[\lambda]]$ be the power series determined by
\[
\frac{d\tau(\lambda)}{d\lambda}=\frac{3(2g_2g'_3-3g'_2g_3)}{2\Delta F(\lambda)^2}-\frac{n}{\lambda}.
\]
By Proposition \ref{A-sympl-3}, $q:=\kappa\lambda^n\exp(\tau(\lambda))$ is the multiplicative period
of the Tate curve $\cE=X\times_A\operatorname{Spec} W((\lambda))$. 
\item(cf. \eqref{tau-eq1}).
Put
\[
\tau^{(\sigma)}(\lambda):=\log^{(\sigma)}(q)=\log^{(F)}(\kappa)-np^{-1}\log(c)+
\tau(\lambda)-p^{-1}\tau(\lambda)^{\sigma}.
\]
\item(cf. \S \ref{global-1ext-sect}).
For $\xi\in K_2(X)$, let $g_\xi(\lambda)\in A$ be defined by
$\dlog\xi=g_\xi(\lambda)\frac{d\lambda}{\lambda}\wedge\frac{dx}{y}$. 
Put
\[
G_\xi(\lambda):=g_\xi(\lambda)F(\lambda)=l+a_1\lambda+a_2\lambda^2+\cdots\in W[[\lambda]].
\]
Note that $l$ is an integer.
\item(cf. Proposition \ref{A-sympl-2}).
Put
\[
H(\lambda):=\frac{2\Delta}{3(2g_2g'_3-3g'_2g_3)}\left(F'(\lambda)+\frac{\Delta'}{12\Delta}F(\lambda)\right).
\]
\item(cf. Theorem \ref{global-reg} \eqref{global-reg-eq2}).
Let $E^\s_1(\lambda)\in W[[\lambda]]$ be defined by
\[
(E^\s_1(\lambda))'=\frac{G_\xi(\lambda)}{\lambda}-\frac{G_\xi(\lambda)^\sigma}{\lambda},\quad
E^\s_1(0):=\frac{l}{n}\log^{(F)}(\kappa)-\frac{l}{p}\log(c).
\]
\end{enumerate}

\medskip

\noindent{\bf Step 1}.
Compute the power series
\[
F(\lambda)=f_0+f_1\lambda+f_2\lambda^2+\cdots+f_N\lambda^N+O(\lambda^{N+1}),
\] 
\[
H(\lambda)=h_0+h_1\lambda+h_2\lambda^2+\cdots+h_N\lambda^N+O(\lambda^{N+1}),
\] 
\[
E^\s_1(\lambda)=A_0+A_1\lambda+A_2\lambda^2+\cdots+A_N\lambda^N+O(\lambda^{N+1}).
\] 
Moreover, let 
\[
E^\s_2(\lambda)=B_0+B_1\lambda+B_2\lambda^2+\cdots+B_N\lambda^N+O(\lambda^{N+1})
\] 
be the power series satisfying
\begin{align*}
(E^\s_2(\lambda))'
&=-E^\s_1(\lambda)\frac{q'}{q}+\frac{G_\xi(\lambda)^\sigma}{\lambda}\tau^{(\sigma)}(\lambda)\\
&=-E^\s_1(\lambda)\left(\frac{n}{\lambda}+\tau'(\lambda)\right)+\frac{G_\xi(\lambda)^\sigma}{\lambda}\tau^{(\sigma)}(\lambda)
\end{align*}
(cf. Theorem \ref{global-reg} \eqref{global-reg-eq3}).
Each $B_i$ is determined except the constant term $B_0$ which shall be computed
in the next step.
Here we take $N$ to be a large integer according to the Step 2 and 3.
See 
Remark \ref{rem-est} \eqref{rem-est-eq2} for an explicit lower bound of $N$.

\medskip

\noindent{\bf Step 2 (Computing $E_2^\s(0)$)}.
Let
\[
\ve_1^{(\sigma)}(\lambda):=\frac{E^\s_1(\lambda)}{F(\lambda)}-H(\lambda)E^\s_2(\lambda),\quad
\ve_2^{(\sigma)}(\lambda):=F(\lambda)E^\s_2(\lambda)
\]
(cf. \eqref{global-1ext-eq1}, \eqref{global-1ext-eq2} and \eqref{S-B-formula-1}).
We compute the constant $E_2^\s(0)=B_0$ modulo $p^s$ (without assuming Conjecture \ref{const-conj}).
We first note that $B_0$ is uniquely determined by the fact that
$\ve_2^\s(\lambda)=F(\lambda)E^\s_2(\lambda)$ is an overconvergent function (i.e. 
$F(\lambda)E^\s_2(\lambda)\in A_K^\dag$). 
Indeed, if there is another constant $B^*_0(\ne B_0)$ such that 
$F(\lambda)(B^*_0+B_1\lambda+\cdots)$ is an overconvergent function, then 
$F(\lambda)E_2^\s(\lambda)-F(\lambda)(B^*_0+B_1\lambda+\cdots)=(B_0-B^*_0)F(\lambda)$ is also
overconvergent. This is impossible. Indeed, 
let $\zeta\in W^\times$ be a root of unity such that
$X_\zeta\times\F_{p^r}$ is an ordinary elliptic curve.
If $F(\lambda)$ were overconvergent, 
then the special value
of $F(\lambda)$ at $\lambda=\zeta$ makes sense.  Then one has
\[
\gamma_\zeta
:=\frac{F(\lambda)}{F(\lambda^{p^r})}\bigg|_{\lambda=\zeta}
=\frac{F(\zeta)}{F(\zeta^{p^r})}
=\frac{F(\zeta)}{F(\zeta)}=1.
\] 
This contradicts with the fact that $\gamma_\zeta$ is the unit root of 
the $p^r$-th Frobenius of $X_\zeta\times\F_{p^r}$ 
(Corollary \ref{unit-root-formual}).
In a practical way, one can use Theorem \ref{exp-cor} to compute $B_0$.
Namely letting
$e_s:=\max\{i\in\Z_{\geq 0}\mid \ord_p(p^i/i!)<s\}$ be as in \eqref{exp-lem-eq1} and 
$d_{2,s}:=\sum_{i=0}^m(-b_i+(p-l_i)e_{s+2})-b_\infty$ as in \eqref{exp-cor-eq0},
one finds the unique 
\[
B_0\text{ mod }p^sW
\] 
satisfying that for a large $s'\geq s$,
\[
\ve_2^{(\sigma)}(\lambda)\equiv\frac{\text{polynomial of degree }\leq d_{2,s'} }
{\prod_{i=0}^m(\lambda-\cano_i)^{-b_i+(p-l_i)e_{s'+2}}}\mod p^{s'}A^\dag
\]
or equivalently, the power series
\[
F(\lambda)(B_0+B_1\lambda+\cdots)\prod_{i=0}^m(\lambda-\cano_i)^{-b_i+(p-l_i)e_{s'+2}}\in(W/p^{s'}W)[[\lambda]]
\]
modulo $p^{s'}$ terminates at least at the degree $d_{2,s'}$.
See Remark \ref{rem-est} \eqref{rem-est-eq3} for an explicit lower bound of $s'$.

\medskip

\noindent{\bf Step 3 (Computing $\ve_k^{(\sigma)}(\cano)$)}.
Compute the power series $\ve_k^{(\sigma)}(\lambda)$.
We again apply Theorem \ref{exp-cor}.
Compute the power series
\[
f_{1,s}(\lambda):=\ve_1^{(\sigma)}(\lambda)\prod_{i=0}^m(\lambda-\cano_i)^{-a_i+(p-l_i)e_{s+2}},\quad
f_{2,s}(\lambda):=\ve_2^{(\sigma)}(\lambda)\prod_{i=0}^m(\lambda-\cano_i)^{-b_i+(p-l_i)e_{s+2}}
\]
modulo $p^sW[[\lambda]]$ which terminates at least at the degree $d_{1,s}$ and $d_{2,s}$
respectively.
Then one has
\[
\ve_1^{(\sigma)}(\cano)\equiv \frac{f_{1,s}(\cano)}
{\prod_{i=0}^m(\cano-\cano_i)^{-a_i+(p-l_i)e_{s+2}}}\mod p^sW,
\]
\[
\ve_2^{(\sigma)}(\cano)\equiv \frac{f_{2,s}(\cano)}
{\prod_{i=0}^m(\cano-\cano_i)^{-b_i+(p-l_i)e_{s+2}}}\mod p^sW.
\]
\begin{rem}\label{rem-est}
Suppose that there are two constants $B_0,B^*_0$ such that
$B_0\not\equiv B^*_0$ mod $p^s$ and
\[
F(\lambda)(B_0+B_1\lambda+\cdots)\prod_{i=0}^m(\lambda-\cano_i)^{-b_i+(p-l_i)e_{s'+2}}\in(W/p^{s'}W)[[\lambda]]
\]
and
\[
F(\lambda)(B^*_0+B_1\lambda+\cdots)\prod_{i=0}^m(\lambda-\cano_i)^{-b_i+(p-l_i)e_{s'+2}}\in(W/p^{s'}W)[[\lambda]]
\]
terminate at the degree $d_{2,s'}$.
Since $\ord_p(B_0-B^*_0)<s$, 
\begin{equation}\label{rem-est-eq1}
F(\lambda)\prod_{i=0}^m(\lambda-\cano_i)^{-b_i+(p-l_i)e_{s'+2}}\in(W/p^{s'-s+1}W)[[\lambda]]
\end{equation}
terminates at the degree $d_{2,s'}$.
This implies that the unit root $\gamma_\zeta$ of the $p^r$-th Frobenius on  $X_\zeta\times_W \F_{p^r}$ satisfies
\[
\gamma_\zeta
=\frac{F(\lambda)}{F(\lambda^{p^r})}\bigg|_{\lambda=\zeta}
\equiv 1\mod p^{s'-s+1}W.
\]
Conversely letting 
\[
N_\zeta:=\min\{i\in \Z_{\geq 1}\mid \gamma_\zeta\not\equiv1\mod p^iW\},
\]
for an arbitrary $s'\geq N_\zeta+s-1$,
the power series \eqref{rem-est-eq1} cannot terminate at the degree $d_{2,s'}$.
Therefore, 
the constant $B_0$ mod $p^s$ is characterized by the condition that
\[
F(\lambda)(B_0+B_1\lambda+\cdots)\prod_{i=0}^m(\lambda-\cano_i)^{-b_i+(p-l_i)e_{s'+2}}\in(W/p^{s'}W)[[\lambda]]
\]
terminates at least at the degree $d_{2,s'}$ for all $s'\geq N_\zeta+s-1$.

\medskip

Summing up the above, we get a bound of $N$ in {\bf Step 1} and a bound of $s'$ in {\bf Step 2}
as follows.
Choose $\zeta\in W$ such that $X_\zeta\times\F_{p^r}$ is ordinary.
Then take $s'$ such that
\begin{equation}\label{rem-est-eq3}
s'\geq N_\zeta+s-1.
\end{equation}
Let
$F(\lambda)\prod_{i=0}^m(\lambda-\cano_i)^{-b_i+(p-l_i)e_{s'+2}}=\sum_{i=0}^\infty C_i\lambda^i$.
Find an integer $M>d_{2,s'}$ such that
\[
C_M\not\equiv0\mod p^{s'-s+1}W.
\]
Then take $N$ such that
\begin{equation}\label{rem-est-eq2}
N>\max\{d_{1,s},d_{2,s},M\}.
\end{equation}

\end{rem}
\subsection{Computing the Frobenius matrix on $H^1_\rig(X_\F/S_\F)$}\label{Frob-sect}
We keep the situation and notation in \S \ref{CSyn-sect}.
Let
\[
\begin{pmatrix}
\Phi(\omega)&\Phi(\eta)
\end{pmatrix}=
\begin{pmatrix}
\omega&\eta
\end{pmatrix}
\begin{pmatrix}
pF_{11}(\lambda)&F_{12}(\lambda)\\
pF_{21}(\lambda)&F_{22}(\lambda)
\end{pmatrix}
\]
be the matrix of the $p$-th Frobenius on $H^1_\rig(X_\F/S_\F)$.
It follows from Proposition \ref{Frob-formula} together with
Proposition \ref{A-sympl-2} \eqref{S-B-formula-1} that one has
\begin{align*}
F_{11}(\lambda)&=\frac{F(\lambda)^\sigma}{F(\lambda)}\left(1+\tau^\s(\lambda)H(\lambda)F(\lambda)\right),\\
F_{12}(\lambda)&=\frac{pH(\lambda)^\sigma}{F(\lambda)}-\frac{H(\lambda)}{F(\lambda)^\sigma}
+p\tau^\s(\lambda)H(\lambda)H(\lambda)^\sigma ,\\
F_{21}(\lambda)&=-\tau^\s(\lambda)F(\lambda)F(\lambda)^\sigma ,\\
F_{22}(\lambda)&=\frac{F(\lambda)}{F(\lambda)^\sigma}-p\tau^\s(\lambda)H(\lambda)^\sigma F(\lambda).
\end{align*}
One can apply Theorem \ref{exp-Frob} to compute 
\[
\begin{pmatrix}
pF_{11}(\cano)&F_{12}(\cano)\\
pF_{21}(\cano)&F_{22}(\cano)
\end{pmatrix}
\mod p^s.
\]
Namely compute the power series
\[
f_{11,s}(\lambda):=pF_{11}(\lambda)\prod_{i=0}^m(\lambda_i^\sigma)^{a_i}\lambda_i^{-a_i+(p-l_i)e_s}\in (W/p^sW)[[\lambda]]
\]
with coefficients in $W/p^sW$ which terminates at least at the degree $d_{11,s}$,
where we note $\lambda_i^\sigma=(\lambda-\cano_i)^\sigma=c\lambda^p-\cano^F_i$
(see \eqref{CSyn-eq00}).
Then one computes $F_{11}(\cano)$ by
\[
pF_{11}(\cano)\equiv f_{11,s}(\cano) \cdot
\left(\prod_{i=0}^m(c\cano^p-\cano_i^F)^{a_i}
(\cano-\cano_i)^{-a_i+(p-l_i)e_s}\right)^{-1}\mod p^sW.
\]
One can also compute $pF_{12}$, $F_{21}$ and $F_{22}$ in the same way.

\section{Numerical Verifications of the $p$-adic Beilinson conjecture
for $K_2$}\label{NVPR-sect}
We apply the algorithm in \S \ref{CSyn-sect} to compute the
syntomic regulators, and give numerical verifications of the $p$-adic Beilinson conjecture
for $K_2$ (Conjecture \ref{PREC}, $n=0$) 
of the following two examples
\begin{itemize}
\item
$y^2=x^3-2x^2+(1-\cano)x$ \text{ with } $\cano\in \Q\setminus\{0,1\}$,
\item
$y^2=x(1-x)(1-(1-\cano)x)$\text{ with }  $\cano\in \Q\setminus\{0,1\}$.
\end{itemize}

\subsection{$y^2=x^3-2x^2+(1-\lambda)x$}\label{exp2-sect}
We first discuss the elliptic fibration in \cite[Theorem 4.29]{New}.
 By a simple replacement of variables, it is a fibration 
$f:Y\to\P^1$
defined by a Weierstrass equation
\[
y^2=4x^3-g_2x-g_3, \quad g_2:=4\lambda+\frac{4}{3},\,g_3:=\frac{8}{3}\lambda-\frac{8}{27}
\]
over the base ring $\Z_{(p)}=\Q\cap\Z_p$ with $p\geq 5$.
Put
\[
\omega=\frac{dx}{y},\quad\eta=\frac{xdx}{y}.
\] 
The functional $j$-invariant is $64(3\lambda+1)^3/(\lambda(1-\lambda)^2)$.
$f$ has a multiplicative reduction of Kodaira type I$_1$ at $\lambda=0$ 
which is split over $\Z_{(p)}$,
a multiplicative reduction of type I$_2$ at $\lambda=1$, and an additive reduction of
type III$^*$ at $\lambda=\infty$.
One can check that $f$ satisfies
all the conditions {\bf A1},\ldots,{\bf A3} in \S \ref{fib-set-sect} and
{\bf A4}, {\bf A5} in \S \ref{exp-sect} (the detail is left to the reader).
Put $S=\operatorname{Spec} \Z_{(p)}[\lambda,(\lambda-\lambda^2)^{-1}]$ and $X=f^{-1}(S)$. 
Write $X_{\Z_p}=X\times_{\Z_{(p)}}\Z_p$, $X_{\F_p}=X\times_{\Z_{(p)}}\F_p$ etc.
Let $\sigma:\Z_p[[\lambda]]\to\Z_p[[\lambda]]$ be the $p$-th Frobenius given by
$\sigma(\lambda)=c\lambda^p$ with $c\in 1+p\Z_p$.
Let
\begin{equation}\label{exp2-symbol}
\xi:=\left\{\frac{3y+6x-2}{3y-6x+2},\frac{-9\lambda (3x+2)}{(3x-1)^3}\right\}\in K_2(X)
\end{equation}
be a symbol (cf. \cite[Theorem 4.29]{New}).
A direct calculation yields
\begin{equation}\label{exp2-sect-eq1}
\dlog(\xi)
=2\frac{d\lambda}{\lambda}\frac{dx}{y}.
\end{equation}
Put 
\[
F_{\frac{1}{4},\frac34}(\lambda)
={}_2F_1\left({\frac{1}{4},\frac{3}{4}\atop 1};\lambda\right).
\]
The notation in \S \ref{CSyn-sect} in our case is explicitly given as follows.
\begin{enumerate}
\item[(i)]
$(a_0,b_0)=(a_1,b_1)=(0,0)$ and
$(a_\infty,b_\infty)=(0,1)$.
\item[(ii)]
Since $1728J=64(3\lambda+1)^3/(\lambda(1-\lambda)^2))$, one has $n=-\ord_\lambda(J)=1$ and
$\kappa=1/64$.
\item[(iii)]
One has $c_0^2=-g_2(0)/(18g_3(0))=1/4$. We take $c_0=1/2$.
\item[(iv)]
One can show
\[
F(\lambda):=(12g_2)^{-\frac{1}{4}}{}_2F_1\left({\frac{1}{12},\frac{5}{12}\atop 1};J^{-1}\right)=\frac{1}{2}F_{\frac{1}{4},\frac34}(\lambda)
\]
in the following way.
Recall from the proof of Proposition \ref{A-sympl-1} that $F(\lambda)$ is characterized by
\[
F(\lambda)=\frac{1}{2\pi i}\int_{\delta_\lambda}\omega.
\]
By Theorem \ref{A1} we have
\[
\begin{pmatrix}
\partial_\lambda(\omega)&
\partial_\lambda(\eta)
\end{pmatrix}
=\begin{pmatrix}
\omega&\eta
\end{pmatrix}
\begin{pmatrix}
\frac{3\lambda-1}{12(\lambda-\lambda^2)}&
-\frac{3\lambda+1}{36(\lambda-\lambda^2)}\\
\frac{1}{4(\lambda-\lambda^2)}&-\frac{3\lambda-1}{12(\lambda-\lambda^2)}
\end{pmatrix},
\]
where $\partial_\lambda:=\nabla_{d/d\lambda}$. Hence one has
\[
((\lambda-\lambda^2)\partial_\lambda^2+(1-2\lambda)\partial_\lambda-\frac34)(\omega)=0.
\]
This implies that $F(\lambda)$ satisfies 
\begin{equation}\label{HGdiff}
(\lambda-\lambda^2)\frac{d^2}{d\lambda^2}F(\lambda)+(1-2\lambda)\frac{d}{d\lambda}F(\lambda)-\frac34 F(\lambda)=0.
\end{equation}
This is the hypergeometric differential equation (e.g. \cite[15.10]{NIST}).
Hence $F(\lambda)=cF_{\frac14,\frac34}(\lambda)$ for some constant $c$.
Since the constant term of $F(\lambda)$ is $c_0=1/2$ by definition, one has $c=1/2$.
\item[(v)]
The multiplicative period is $q=64^{-1}\lambda\exp(\tau(\lambda))$ with
\[
\lambda\frac{d}{d\lambda}\tau(\lambda)=-1+(1-\lambda)^{-1}F_{\frac14,\frac34}(\lambda)^{-2},\quad \tau(0)=0.
\]
Explicitly, 
\[
\tau(\lambda)=-\frac{5}{8} \lambda-\frac{269}{1024} \lambda^{2}-\frac{1939}{12288} \lambda^{3}
-\frac{922253}{8388608} \lambda^4+\cdots.
\]
\item[(vi)]
$
\tau^{(\sigma)}(\lambda)=-p^{-1}\log(64^{p-1}c)+
\tau(\lambda)-p^{-1}\tau(\lambda)^{\sigma}.
$
\item[(vii)]
By \eqref{exp2-sect-eq1}, one has $g_\xi(\lambda)=2$ and
$G_\xi(\lambda)=2F(\lambda)=F_{\frac{1}{4},\frac34}(\lambda)$ (hence $l=1$).
\item[(viii)]
\begin{align*}
H(\lambda)&=2(\lambda-\lambda^2)F_{\frac{1}{4},\frac34}^\prime(\lambda)+\frac{1}{6}
(1-3\lambda)F_{\frac{1}{4},\frac34}(\lambda)\\
&=\frac{1}{6}-\frac{3}{32} \lambda-\frac{85}{2048} \lambda^{2}-\frac{875}{32768} \lambda^{3}
-\frac{165165}{8388608} \lambda^{4}+\cdots.
\end{align*}
\item[(ix)]
Let $E^\s_1(\lambda)\in \Z_p[[\lambda]]$ be defined by
\[
(E^\s_1(\lambda))'=\frac{G_\xi(\lambda)}{\lambda}-\frac{G_\xi(\lambda)^\sigma}{\lambda},\quad
E^\s_1(0):=-p^{-1}\log(64^{p-1}c).
\]
\end{enumerate}
Let $E^\s_2(\lambda)=C+a_1\lambda+\cdots$ be the power series satisfying
\begin{align*}
\frac{d}{d\lambda}E^\s_2(\lambda)
&=-E^\s_1(\lambda)\left(\frac{1}{\lambda}+\tau'(\lambda)\right)+\frac{G_\xi(\lambda)^\sigma}{\lambda}\tau^{(\sigma)}(\lambda)\\
&=-\frac{E^\s_1(\lambda)}{\lambda}\left((1-\lambda)^{-1}F_{\frac{1}{4},\frac34}(\lambda)^{-2}\right)+\frac{G_\xi(\lambda)^\sigma}{\lambda}\tau^{(\sigma)}(\lambda).
\end{align*}
This determines $E^\s_2(\lambda)$ except the constant term $C$.
We expect
\[
C=-\frac{1}{2}\left(p^{-1}\log(64^{p-1}c)\right)^2
\]
according to Conjecture \ref{const-conj}.
The authors do not know a proof of this, while it does not matter toward
numerical computation thanks to the method in {\bf Step 2} in \S \ref{CSyn-sect}.

\medskip

Let $\cano\in \Z_{(p)}$ such that $\cano\not\equiv0,1$ mod $p$, and
$X_\cano$ the fiber at $\lambda=\cano$.
Let $\sigma_\cano(\lambda)=\cano^{1-p}\lambda^p$, and
\[
\reg_\syn(\xi|_{\lambda=\cano})=\ve^{(\sigma_\cano)}_1(\cano)\frac{dx}{y}
+\ve^{(\sigma_\cano)}_2(\cano)\frac{xdx}{y}.
\]
If $\cano=2,4,8$, then the symbol $\xi$ is integral.
Running the steps in \S \ref{CSyn-sect}, we obtain the following table:
\begin{itemize}
\item[\fbox{$\cano=2$}]
 \begin{tabular}{c|c|c}
 $p$&$\ve^{(\sigma_\cano)}_1(\cano)$&$\ve^{(\sigma_\cano)}_2(\cano)$\\
 \hline
$5$&$  19674787784(p^{15})$&$ 14084682764(p^{15})$\\
$7$&$2590802242826(p^{15})$&$ 1117519637850(p^{15})$\\
$11$&$ 7965734565(p^{10})$&
$ 13678445088(p^{10})$\\
$13$&$ 11926785562 (p^{10})$&$ 72813360585(p^{10})$\\
$17$&
$ 278698244(p^8)$&
$3681687094(p^8)$
\end{tabular} ``$(p^s)$'' means ``mod $p^s$''.
\item[\fbox{$\cano=4$}]
 \begin{tabular}{c|c|c}
 $p$&$\ve^{(\sigma_\cano)}_1(\cano)$&$\ve^{(\sigma_\cano)}_2(\cano)$\\
 \hline
$5$&$ 10899377190(p^{15})$&$ 18907289609(p^{15})$\\
$7$&$1079749670299(p^{15})$&$ 3800305332560(p^{15})$\\
$11$&$1981716464(p^{10})$&
$ 23499057595(p^{10})$\\
$13$&$ 99992661864(p^{10})$&$ 54272943535(p^{10})$\\
$17$&
$5624041300(p^8)$&
$ 6557787442(p^8)$
\end{tabular} 
\item[\fbox{$\cano=8$}]
 \begin{tabular}{c|c|c}
 $p$&$\ve^{(\sigma_\cano)}_1(\cano)$&$\ve^{(\sigma_\cano)}_2(\cano)$\\
 \hline
$5$&$ 25874535623(p^{15})$&$ 28249557885(p^{15})$\\
$7$& bad reduction\\
$11$&$ 11113392970(p^{10})$&
$ 10168663913(p^{10})$\\
$13$&$50210889322(p^{10})$&$ 1878713357(p^{10})$\\
$17$&
$4803397460(p^8)$&
$ 6578109195(p^8)$
\end{tabular}
 \end{itemize}

\medskip

Next we compute the $p$-adic regulators
\[
R_{p,\gamma}(X_\cano,\xi)=R_{p,\gamma}(h^1(X_\cano),\xi)=(1-p\gamma^{-1})\frac{\Tr(\reg_\syn(\xi|_{\lambda=\cano})\cup v_{\gamma})}
{\Tr(\omega_{X_\cano}\cup v_{\gamma})}
\]
introduced in \S \ref{PRconjEC},
where $\gamma$ is a unit or non-unit root of $X_\cano$, 
$v_\gamma\in H^1_\dR(X_\cano/\Q_p)\ot\ol\Q_p$ is the eigenvector of
$\Phi\ot1_{\ol\Q_p}$ with
eigenvalue $\gamma$ and 
$\omega_{X_\cano}$ is the N\'eron differential.
To do this we first compute
the representation matrix 
\[
\begin{pmatrix}
\Phi(\omega)&\Phi(\eta)
\end{pmatrix}=
\begin{pmatrix}
\omega&\eta
\end{pmatrix}
\begin{pmatrix}
pF_{11}(\lambda)&F_{12}(\lambda)\\
pF_{21}(\lambda)&F_{22}(\lambda)
\end{pmatrix}\bigg|_{\lambda=\cano}
\]
of the $p$-th Frobenius 
on $H_\crys(X_{\cano,\F_p}/\Z_p)$.
Let us see the case $\cano=4$ and $p=5$ (the other cases are done in the same way).
Let $\alpha$ be the unit root, and $\beta=p\alpha^{-1}$ the non-unit root.
By counting the $\F_5$-rational points of $X_{\cano,\F_5}$,
one has that the trace of the Frobenius is $-2$, and hence $\alpha=-1+2\sqrt{-1}\equiv 22219310363$ mod $5^{15}$.
Thanks to the method in \S \ref{Frob-sect}, one has
\begin{equation}\label{exp2-sect-eq2}
\begin{pmatrix}
pF_{11}(\cano)&F_{12}(\cano)\\
pF_{21}(\cano)&F_{22}(\cano)
\end{pmatrix}
\equiv
\begin{pmatrix}30465757535&
       3584355845\\
       11560944585&
       51820588\\
       \end{pmatrix}
\mod 5^{15}.
\end{equation}
Hence
the eigenvector $v_\alpha$ satisfies
\[
v_\alpha\equiv710786365\omega+\eta\mod 5^{14},
\]
and the eigenvector $v_\beta$
for the non-unit root $\beta$ satisfies
\[
v_{\beta}\equiv\omega+ 28450379180\eta\mod 5^{15}.
\]
In cases $\cano=2,4,8$, it follows by computing the minimal
Weierstrass equation (e.g. \textsc{Magma} \cite{magma}) that
one has $\omega_{X_\cano}=\omega$.
We thus have
\begin{align*}
\frac{\Tr(\reg_\syn(\xi|_{\lambda=\cano})\cup v_{\alpha})}
{\Tr(\omega_{X_\cano}\cup v_{\alpha})}
&\equiv
(\ve^{(\sigma_\cano)}_1(\cano)-710786365\ve^{(\sigma_\cano)}_2(\cano))
\frac{\Tr(\omega\cup\eta)}
{\Tr(\omega_{X_\cano}\cup\eta)}\\
&\equiv   3956230280\mod 5^{14},
\end{align*}
\begin{align*}
\frac{\Tr(\reg_\syn(\xi|_{\lambda=\cano})\cup v_{\beta})}
{\Tr(\omega_{X_\cano}\cup v_{\beta})}
&\equiv
\frac{28450379180\ve^{(\sigma_\cano)}_1(\cano)-
\ve^{(\sigma_\cano)}_2(\cano)}{28450379180}
\frac{\Tr(\omega\cup\eta)}
{\Tr(\omega_{X_\cano}\cup\eta)}\\
&\equiv  5^{-1}\cdot5576309542
\mod 5^{13},
\end{align*}
where $a\equiv b$ mod $p^s$ for $a,b\in \Q_p$ means that $a-b\in p^s\Z_p$.
Hence
\begin{align*}
R_{5,\alpha}(X_\cano,\xi)&\equiv  5110910605\mod 5^{14},\\
R_{5,5\alpha^{-1}}(X_\cano,\xi)
&\equiv  5^{-1}\cdot2495229428\mod 5^{13}.
\end{align*}
If $\cano=4$ and $p=7$, then $X_\cano$ has a supersingular reduction.
In this case the Frobenius is
\begin{equation}\label{exp2-sect-eq3}
\begin{pmatrix}
pF_{11}(\cano)&F_{12}(\cano)\\
pF_{21}(\cano)&F_{22}(\cano)
\end{pmatrix}
\equiv
\begin{pmatrix}3603314253994&
      2631033376372\\
      899507565369&
      1144247255949\\
      \end{pmatrix}
\mod 7^{15},
\end{equation}
and the eigenvector $v_{\sqrt{-7}}$ for $\sqrt{-7}$ satisfies
\[
v_{\sqrt{-7}}\equiv\omega+ 
\overbrace{(2066673311059+  1689081159560\sqrt{-7})}^C
\eta\mod 7^{15}\Z_p[\sqrt{-7}].
\]
Hence
\begin{align*}
R_{7,\sqrt{-7}}(X_\cano,\xi)
&=(1+\sqrt{-7})\frac{\Tr(\reg_\syn(\xi|_{t=\cano})\cup v_{\sqrt{-7}})}
{\Tr(\omega_{X_\cano}\cup v_{\sqrt{-7}})}\\
&\equiv (1+\sqrt{-7})\frac{C\ve^{(\sigma_\cano)}_1(\cano)-
\ve^{(\sigma_\cano)}_2(\cano)}{C}\\
&\equiv  224156767500+7^{-1}\cdot577165355498\sqrt{-7}
\end{align*}
modulo $7^{13}\Z_p[\sqrt{-7}]$.
Other results are given in the following table,
where $\alpha$ denotes the unit root or $\sqrt{-p}$:
\begin{center}
{\bf Table of $p$-adic regulators}
\end{center}
\begin{itemize}
\item[\fbox{$\cano=2$}]
 \begin{tabular}{c|c|c}
 $p$&$R_{p,\alpha}(X_\cano,\xi)$&$R_{p,\beta}(X_\cano,\xi)$\\
 \hline
$5$&$5162207171(p^{15})$&$5^{-1}\cdot  4698651004(p^{13})$\\
$7$&$1898050197813(p^{15})$ &$  47494881888(p^{13})$\\
$11$&$4688873516(p^{10})$&
$4077822 (p^{8})$\\
$13$&$  109756622389(p^{10})$&$13^{-1}\cdot  1014980225(p^{8})$\\
$17$&
$ 1874620233(p^8)$&$17^{-1}\cdot  239203868(p^6)$
\end{tabular} \quad
$(p^s)=$ mod $p^s$
\item[\fbox{$\cano=4$}]
 \begin{tabular}{c|c|c}
 $p$&$R_{p,\alpha}(X_\cano,\xi)$&$R_{p,\beta}(X_\cano,\xi)$\\
 \hline
$5$&$ 5110910605(p^{14})$&$5^{-1}\cdot2495229428(p^{13})$\\
$7$&$224156767500+7^{-1}\cdot577165355498\sqrt{-7} (p^{13})$ &(supersingular)\\
$11$&$23123035879(p^{10})$&
$ 11^{-1}\cdot1673260489 (p^{8})$\\
$13$&$  131267814396(p^{10})$&$13^{-1}\cdot8897873775(p^{8})$\\
$17$&
$ 4659831644(p^8)$&$17^{-1}\cdot  401683799(p^6)$
\end{tabular} 
\item[\fbox{$\cano=8$}]
 \begin{tabular}{c|c|c}
 $p$&$R_{p,\alpha}(X_\cano,\xi)$&$R_{p,\beta}(X_\cano,\xi)$\\
 \hline
$5$&$299743657+ 14397739353\sqrt{-5}(p^{14})$&(supersingular)\\
$7$& bad reduction\\
$11$&$ 7356580887(p^{10})$&
$11^{-1}\cdot1944330879(p^{8})$\\
$13$&$ 89811548583(p^{10})$&$ 13^{-1}\cdot7256384468(p^{8})$\\
$17$&
$4687529300(p^8)$&
$17^{-1}\cdot   270140891(p^6)$
\end{tabular}
 \end{itemize}

 \medskip
 
We compare the above regulators with the $p$-adic $L$-values.
One can compute the special values $L_{p,\alpha}(X_\cano,\omega^{-1},0)$
by SAGE. The values $L_{p,\beta}(X_\cano,\omega^{-1},0)$ 
of the critical slope $p$-adic $L$-function are kindly provided by Professor Robert Pollack:

\medskip

 \begin{itemize}
\item[\fbox{$\cano=2$}]
   \begin{tabular}{c|c}
&$L_{p,\beta}(X_\cano,\omega^{-1},0)$\\
 \hline
$p=5$& $2208371/5(p^{10})$\\
$p=7$&$1125190415157500314434008856467880(p^{30})$\\
$p=11$&$88350501278121136930215849171019600321077(p^{30})$
\end{tabular}
\item[\fbox{$\cano=4$}]
   \begin{tabular}{c|c}
$a=4$&$L_{p,\beta}(X_\cano,\omega^{-1},0)$\\
 \hline
$p=5$& $ 19846135286/5(p^{15})$\\
$p=11$&$ 312624858906775755190263231723045879608690/11(p^{30})$\\
$p=13$&$475670741377552429015509835731352129856903404/13(p^{30})$\\
\end{tabular}
 \end{itemize}

The above indicate the following:

\medskip

\begin{tabular}{c|c|c}
 $\cano$
 &$R_{p,\alpha}/L_{p,\alpha}(X_\cano,\omega^{-1},0)$
 &$R_{p,\beta}/L_{p,\beta}(X_\cano,\omega^{-1},0)$\\
 \hline
$2$&$-1$&$-1$\\
$4$&$-2$&$-2$\\
$8$&$1/2$&\bysame
\end{tabular}

\bigskip

Next we compute the Beilinson regulators. 
Let $X_\lambda:=f^{-1}(\lambda)$ be the fiber at
$\lambda\in \R_{>0}\setminus\{1\}$, an elliptic curve over $\R$.
Let $F_\infty$ denote the infinite Frobenius
on $X_\lambda(\C)$.
Fix generators $u^\pm\in H_1(X_\lambda(\C),\Z)^{F_\infty=\pm1}\cong\Z$.
Since the elliptic curve $X_\lambda$ is the rectangle case
(i.e. $4x^3-g_2x-g_3$ has 3 distinct
real roots), one has $H_1(X_\lambda(\C),\Z)=\Z u^++\Z u^-$.
Let
\[
\reg_\R:K_2(X_\lambda)\lra H^2_\cD(X_\lambda,\R(2))\cong \Hom(H_1(X_\lambda(\C),\Z),\R(1))
\]
be the Beilinson regulator map. Then
\[
R_\infty(X_\lambda,\xi)=\frac{1}{2\pi i}\langle\reg_\R(\xi),u^-\rangle.
\]
Note that in general
\begin{equation}\label{exp2-sect-eq4}
\langle\reg_\R(\xi),u^+\rangle\in\R,\quad
\langle\reg_\R(\xi),u^-\rangle\in\R(1).
\end{equation}
We claim
\begin{equation}\label{exp2-sect-eq5}
\frac{1}{2\pi i}\langle\reg_\R(\xi),u^-\rangle=\pm
\mathrm{Re}\left[-\log 64+\log \lambda+\frac{3\lambda}{16}\,{}_4F_3\left(
{\frac54,\frac74,1,1\atop 2,2,2};\lambda\right)\right].
\end{equation}
If $0<\lambda<1$, the equation $4x^3-g_2x-g_3$ has 3 distinct
real roots
\[
-\frac23<\frac13-\sqrt\lambda<\frac13+\sqrt\lambda.
\]
Hence $u^\pm$ are the cycles in the $x$-plane displayed as follows:
\begin{center}
{\unitlength 0.1in%
\begin{picture}( 23.7500,  7.8500)(  3.0500,-13.8000)%
\put(18.1000,-7.4000){\makebox(0,0)[lb]{$u^-$}}%
\put(4.7800,-10.5800){\makebox(0,0)[lb]{$-\frac23$}}%
\put(13.9200,-10.5200){\makebox(0,0)[lb]{$Q_1$}}%
\put(21.5400,-10.5400){\makebox(0,0)[lb]{$Q_2$}}%
%
\special{pn 13}%
\special{pa 1050 1210}%
\special{pa 1100 1212}%
\special{fp}%
\special{sh 1}%
\special{pa 1100 1212}%
\special{pa 1034 1189}%
\special{pa 1047 1210}%
\special{pa 1033 1229}%
\special{pa 1100 1212}%
\special{fp}%
\put(8.4000,-7.4000){\makebox(0,0)[lb]{$u^+$}}%
%
\special{pn 8}%
\special{ar 1940 1000 740 210  0.0000000  6.2831853}%
%
\special{pn 8}%
\special{ar 1090 1000 740 210  0.0000000  6.2831853}%
%
\special{pn 13}%
\special{pa 1910 1210}%
\special{pa 1960 1212}%
\special{fp}%
\special{sh 1}%
\special{pa 1960 1212}%
\special{pa 1894 1189}%
\special{pa 1907 1210}%
\special{pa 1893 1229}%
\special{pa 1960 1212}%
\special{fp}%
%
\special{pn 4}%
\special{sh 1}%
\special{ar 588 1080 16 16 0  6.28318530717959E+0000}%
%
\special{pn 4}%
\special{sh 1}%
\special{ar 1474 1078 16 16 0  6.28318530717959E+0000}%
\special{sh 1}%
\special{ar 2226 1080 16 16 0  6.28318530717959E+0000}%
\put(3.0500,-15.2500){\makebox(0,0)[lb]{$Q_1=\frac13-\sqrt\lambda$, $Q_2=\frac13+\sqrt\lambda$, $0<\lambda<1$}}%
\end{picture}}%

\end{center}
If $\lambda>1$, the equation $4x^3-g_2x-g_3$ has 3 distinct
real roots
\[
\frac13-\sqrt\lambda<-\frac23<\frac13+\sqrt\lambda,
\]
the generator $u^\pm$ are as follows:
\begin{center}
{\unitlength 0.1in%
\begin{picture}( 24.6000,  8.0000)(  2.2000,-13.9500)%
\put(18.1000,-7.4000){\makebox(0,0)[lb]{$u^-$}}%
\put(13.4600,-10.5400){\makebox(0,0)[lb]{$-\frac23$}}%
\put(4.9000,-10.4000){\makebox(0,0)[lb]{$Q_1$}}%
\put(20.8000,-10.3400){\makebox(0,0)[lb]{$Q_2$}}%
%
\special{pn 13}%
\special{pa 1050 1210}%
\special{pa 1100 1212}%
\special{fp}%
\special{sh 1}%
\special{pa 1100 1212}%
\special{pa 1034 1189}%
\special{pa 1047 1210}%
\special{pa 1033 1229}%
\special{pa 1100 1212}%
\special{fp}%
\put(8.4000,-7.4000){\makebox(0,0)[lb]{$u^+$}}%
%
\special{pn 8}%
\special{ar 1940 1000 740 210  0.0000000  6.2831853}%
%
\special{pn 8}%
\special{ar 1090 1000 740 210  0.0000000  6.2831853}%
%
\special{pn 13}%
\special{pa 1910 1210}%
\special{pa 1960 1212}%
\special{fp}%
\special{sh 1}%
\special{pa 1960 1212}%
\special{pa 1894 1189}%
\special{pa 1907 1210}%
\special{pa 1893 1229}%
\special{pa 1960 1212}%
\special{fp}%
%
\special{pn 4}%
\special{sh 1}%
\special{ar 578 1058 16 16 0  6.28318530717959E+0000}%
%
\special{pn 4}%
\special{sh 1}%
\special{ar 1464 1058 16 16 0  6.28318530717959E+0000}%
%
\special{pn 4}%
\special{sh 1}%
\special{ar 2154 1058 16 16 0  6.28318530717959E+0000}%
\put(2.2000,-15.4000){\makebox(0,0)[lb]{$Q_1=\frac13-\sqrt\lambda$, $Q_2=\frac13+\sqrt\lambda$, $\lambda>1$}}%
\end{picture}}%

\end{center}
If $0<\lambda<1$, the cycle $u^-$ agrees with the vanishing cycle $\delta_\lambda$ at $\lambda=0$.
We move $\lambda$ to a point in the interval $(1,\infty)$. 
Then $\delta_\lambda=u^-+mu^+$ with some $m\in\Z$.
Hence it follows from \eqref{exp2-sect-eq4} that
\begin{equation}\label{exp2-sect-eq6}
\frac{1}{2\pi i}\langle\reg_\R(\xi),u^-\rangle
=\mathrm{Re}\left(\frac{1}{2\pi i}\langle\reg_\R(\xi),\delta_\lambda\rangle\right)
\end{equation}
in both cases of $0<\lambda<1$ and $\lambda>1$.
By \eqref{exp2-sect-eq1} one has
\[
f(\lambda):=\lambda\frac{d}{d\lambda}\langle\reg_\R(\xi),\delta_\lambda\rangle=2\left\langle \frac{dx}{y},
\delta_\lambda\right\rangle=2\int_{\delta_\lambda}\frac{dx}{y}=4\int_{\frac13-\sqrt\lambda}^{\frac13+\sqrt\lambda}\frac{dx}{y}
\]
which is a single valued analytic function
on $0\leq \mathrm{arg}(z-1)<2\pi$.
Since $\delta_\lambda$ is the vanishing cycle,
$f(\lambda)$ is holomorphic at $\lambda=0$.
Moreover it satisfies the hypergeometric differential
equation \eqref{HGdiff}.
Therefore $f(\lambda)=cF_{\frac14,\frac34}(\lambda)$ with some constant $c$.
Since
\[f(0)=\lim_{\lambda\to+0} 4\int_{\frac13-\sqrt\lambda}^{\frac13+\sqrt\lambda}\frac{dx}{y}=\pm2\pi i,\]
one concludes
$f(\lambda)=\pm2\pi i\,F_{\frac14,\frac34}(\lambda)$.
Hence we have
\[
\frac{1}{2\pi i}\langle\reg_\R(\xi),u^-\rangle\os{\eqref{exp2-sect-eq6}}{=}
\mathrm{Re}\left(\frac{1}{2\pi i}\langle\reg_\R(\xi),\delta_\lambda\rangle\right)=
\pm\mathrm{Re}
\left[C+\log \lambda+\frac{3\lambda}{16}\,{}_4F_3\left(
{\frac54,\frac74,1,1\atop 2,2,2};\lambda\right)\right]
\]
with some constant $C$. One can determine $C$ by a similar argument
to the proof of \cite[Theorem 3.2]{A}. 
We thus have \eqref{exp2-sect-eq5}.

\medskip

We may take
\[
R_\infty(X_\lambda,\xi)=
\mathrm{Re}\left[-\log 64+\log \lambda+\frac{3\lambda}{16}\,{}_4F_3\left(
{\frac54,\frac74,1,1\atop 2,2,2};\lambda\right)\right]
\]
by replacing $u^-$ with $-u^-$ if necessary.
Here is the table of numerical computations (of precision 20)
by MATHEMATICA for $_4F_3$
and MAGMA for $L$-values:

 \medskip
 
\begin{center}
{\bf Table of Beilinson regulators and $L$-values}
\end{center}

   \begin{tabular}{c|c|c|c}
 $\cano$&$R_{\infty}(X_\cano,\xi)$&$L'(X_\cano,0)$&$R_{\infty}/L'(X_\cano,0)$\\
 \hline
$2$&$- 3.01582754200021$&$3.01582754200021$&$-1$\\
$4$& $ -2.47055987085139$&$ 1.23527993542569$&$-2$\\
$8$&$-2.04611026378501$&$ -4.09222052757003$&$0.5$
\end{tabular}

\bigskip

We thus have numerical verifications on the $p$-adic Beilinson conjecture for $K_2$ of
the above elliptic curves.
\begin{rem}\label{remark-BK}
Let $z_{\mathrm{BK}}\in \varprojlim_n K_2(Y(Np^n))_\Q$ be the Beilinson-Kato element.
We denote by $z_{\mathrm{BK}}(2)\in K_2(E)_\Q$ its image for $E$ an elliptic curve
of conductor $N$.
Brunault\footnote{Brunault \cite[Th\'eor\`eme 38 (115)]{Br}
obtains a formula in terms of 
the composition of the \'etale regulator map and
the Bloch-Kato exponential map
\[
\log\reg_f:K_2(E)\lra H^1_\dR(E/\Q_p).
\]
By \cite[Proposition 9.11]{Be1}, this agrees with the syntomic regulator map,
more precisely
$(1-p^{-2}\phi)\log\reg_f=\reg_\syn$, where $\phi$ is the $p$-th Frobenius on 
$H^1_\dR(E/\Q_p)\cong H^1_\crys(E_{\F_p}/\Z_p)\ot\Q$ (see also \S \ref{syn-sect}).
Hence one has $
\Tr(\reg_\syn(z)\cup v_\alpha)=
(1-p^{-1}\alpha^{-1})\Tr(\log\reg_f(z)\cup v_\alpha)
$ and then \cite[Th\'eor\`em 38 (115)]{Br} reads \eqref{remark-BK-eq1}.} 
showed that if $E$ does not have CM, then 
\begin{equation}\label{remark-BK-eq1}
R_{p,\alpha}(E,z_{\mathrm{BK}}(2))=\prod_{l|N}(1-a_l(E))\frac{L(E,1)\Omega^-_E}{i\langle
f,f\rangle}
L_{p,\alpha}(E,\omega^{-1},0).
\end{equation}
In particular the $p$-adic regulator vanishes when $a_l(E)=1$ for some $l|N$
or $L(E,1)=0$ (he further expects $z_{\mathrm{BK}}(2)=0$, loc.cit.
Question 12).
Concerning $y^2=x^3-2x^2+(1-a)x$, if $a=4$ then $a_3(E)=1$ and if $a=8$
then $L(E,1)=0$.
In these cases one cannot use $z_{\mathrm{BK}}$ to
discuss the $p$-adic Beilinson conjecture for $K_2$.
Instead we use the symbol $\xi$ in \eqref{exp2-symbol}, and obtain the numerical verifications.
\end{rem}

 \subsection{Legendre family $y^2=x(1-x)(1-(1-\lambda)x)$}\label{exp3-sect}
Next example is the Legendre family given by a Weierstrass equation
$y_0^2=x_0(1-x_0)(1-(1-\lambda)x_0)$.
Let $p\geq 5$ be a prime.
Replace the variables $(x,y)=(x_0^{-1}+(\lambda-2)/3,2y_0x_0^{-2})$.
We have a family
\[
f:Y\lra \P^1
\]
of elliptic curves over $\Z_{(p)}=\Z_p\cap \Q$
whose generic fiber is given by a Weierstrass equation 
\[
y^2=4x^3-g_2x-g_3, \quad g_2:=\frac{4}{3}(\lambda^2-\lambda+1),\,g_3:=-\frac{4}{27}(\lambda+1)(\lambda-2)(2\lambda-1).
\] 
Put $S:=\operatorname{Spec}\Z_{(p)}[\lambda,(\lambda-\lambda^2)^{-1}]$ and $X:=f^{-1}(S)$.
The functional $j$-invariant is
$256(\lambda^2-\lambda+1)^3/(\lambda^2(1-\lambda)^2)$.
The singular fibers appear at $\lambda=0,1,\infty$, and the Kodaira symbols
are I$_2$, I$_2$ and I$_2^*$ respectively.
Moreover the fiber at $\lambda=0$ is split multiplicative over $\Z_{(p)}$, while 
so is not the fiber at $\lambda=1$.
One can check that $f$ satisfies
all the conditions {\bf A1},\ldots,{\bf A3} in \S \ref{fib-set-sect} and
{\bf A4}, {\bf A5} in \S \ref{exp-sect} (the detail is left to the reader).
Let
\[
\omega:=\frac{dx}{y}=-\frac{dx_0}{2y_0},\quad\eta:=\frac{xdx}{y}
=-\left(x_0^{-1}+\frac{\lambda-2}{3}\right)\frac{dx_0}{2y_0}
\] 
be a $\O(S)$-basis of $H^1_\dR(X/S)$. 
Let
\begin{equation}\label{RZ-eq2}
\xi=\left\{\frac{y_0+1-x_0}{y_0-1+x_0},\frac{\lambda x_0^2}{(1-x_0)^2)}\right\}\in K_2(X)
\end{equation}
be a symbol (cf. \cite[(4.52)]{New}). One immediately has
\[
\dlog(\xi)=-\frac{d\lambda}{\lambda}\frac{dx_0}{y_0}=2\frac{d\lambda}{\lambda}\frac{dx}{y}.
\]
Put $\Delta:=g_2^3-27g_3^2=16\lambda^2(1-\lambda)^2$, and 
\[
F_{\frac12,\frac12}(\lambda)
={}_2F_1\left({\frac12,\frac12\atop 1};\lambda\right).
\]
The notation in \S \ref{CSyn-sect} in our case is explicitly given as follows.
\begin{enumerate}
\item[(i)]
$(a_0,b_0)=(a_1,b_1)=(0,0)$ and
$(a_\infty,b_\infty)=(0,1)$.
\item[(ii)]
$n=-\ord_\lambda(J)=2$,
$\kappa=1/256$.
\item[(iii)]
One has $c_0^2=-g_2(0)/(18g_3(0))=1/4$. We take $c_0=1/2$.
\item[(iv)]
One can show
\[
F(\lambda)=\frac{1}{2}F_{\frac12,\frac12}(\lambda)
\]
in the same way as in \S \ref{exp2-sect} (iv).
\item[(v)]
The multiplicative period is $q=256^{-1}t^2\exp(\tau(\lambda))$ with
\[
\lambda\frac{d}{d\lambda}\tau(\lambda)=-2+2(1-\lambda)^{-1}F_{\frac12,\frac12}(\lambda)^{-2},\quad \tau(0)=0.
\]
Explicitly, 
\[
\tau(\lambda)=-\lambda-\frac{13}{32}\lambda^2-\frac{23}{96}\lambda^3-\frac{2701}{16384} \lambda^{4}+\cdots.
\]
\item[(vi)]
$
\tau^{(\sigma)}(\lambda)=-p^{-1}\log(256^{p-1}c^2)+
\tau(\lambda)-p^{-1}\tau(\lambda)^{\sigma}.
$
\item[(vii)]
By \eqref{exp2-sect-eq1}, one has $g_\xi(\lambda)=2$ and
$G_\xi(\lambda)=2F(\lambda)=F_{\frac12,\frac12}(\lambda)$ (hence $l=1$).
\item[(viii)]
\begin{align*}
H(\lambda)&=(\lambda-\lambda^2)F_{\frac12,\frac12}^\prime(\lambda)+\frac{1}{6}
(1-2\lambda)F_{\frac12,\frac12}(\lambda)\\
&=\frac{1}{6}-\frac{1}{24}\lambda-\frac{11}{384} \lambda^2-\frac{29}{1536} \lambda^3
-\frac{1375}{98304} \lambda^4+\cdots.
\end{align*}
\item[(ix)]
Let $E^\s_1(\lambda)\in \Z_p[[\lambda]]$ be defined by
\[
(E^\s_1(\lambda))'=\frac{G_\xi(\lambda)}{\lambda}-\frac{G_\xi(\lambda)^\sigma}{\lambda},\quad
E^\s_1(0):=-p^{-1}\log(16^{p-1}c).
\]
\end{enumerate}
Let $E^\s_2(\lambda)=C+a_1\lambda+\cdots$ be the power series satisfying
\begin{align*}
\frac{d}{d\lambda}E^\s_2(\lambda)
&=-E^\s_1(\lambda)\left(\frac{2}{\lambda}+\tau'(\lambda)\right)+\frac{G_\xi(\lambda)^\sigma}{\lambda}\tau^{(\sigma)}(\lambda)\\
&=-2\frac{E^\s_1(\lambda)}{\lambda}\left((1-\lambda)^{-1}F_{\frac12,\frac12}(\lambda)^{-2}\right)+\frac{G_\xi(\lambda)^\sigma}{\lambda}\tau^{(\sigma)}(\lambda).
\end{align*}
This determines $E^\s_2(\lambda)$ except the constant term $C$.
We expect
\[
C=-\frac{1}{4}\left(p^{-1}\log(256^{p-1}c^2)\right)^2=-\left(p^{-1}\log(16^{p-1}c)\right)^2
\]
according to Conjecture \ref{const-conj}, while it does not matter toward
numerical computation thanks to the method in {\bf Step 2} in \S \ref{CSyn-sect}.
Let $X_a$ be the fiber at $\lambda=a$.
Let $\sigma_\cano(\lambda)=\cano^{1-p}\lambda^p$, and
\[
\reg_\syn(\xi|_{\lambda=\cano})=\ve^{(\sigma_\cano)}_1(\cano)\frac{dx}{y}
+\ve^{(\sigma_\cano)}_2(\cano)\frac{xdx}{y}.
\]

We now discuss a particular case $\cano=4$.
In this case, $X_4$ is isogenous to $X_0(24)$ and
the N\'eron differential $\omega_{X_\cano}$ is $dx/y$.
\begin{thm}\label{X4-thm}
Let $p$ be a prime at which $X_4$ has a good reduction
(including a supersingular reduction).
Then
Conjecture \ref{PREC} for $X_4$ is true in case $n=0$ and $\ord_p(\gamma)<1$.
\end{thm}
\begin{pf}
Noticing
\[
\lim_{\lambda\to+0}\lambda\frac{d}{d\lambda}\langle\reg_\R(\xi),u^-\rangle=2
\lim_{\lambda\to+0}\int_1^{(1-\lambda)^{-1}}\frac{dx_0}{y_0}=2\pi i 
\]
one can show that
the Beilinson regulator
is
\[
R_\infty(X_\cano,\xi)=\frac{1}{2\pi i}
\langle\reg_\R(\xi),u^-\rangle=
\mathrm{Re}\left[-\log 16+\log \cano+\frac{\cano}{4}\,
{}_4F_3\left({\frac{3}{2},\frac{3}{2},1,1\atop 2,2,2};\cano\right)\right]
\]
in the same way as in \S \ref{exp2-sect} (the orientation of $u^-$ is chosen suitably).
By \cite[Lemma 3.5]{A} (or \cite[Example 5.3]{As-Ross1}), this is replaced with
\[
-2\mathrm{Re}\left[\cano^{-\frac{1}{2}}{}_3F_2\left({\frac{1}{2},\frac{1}{2},\frac{1}{2}\atop
1,\frac{3}{2}};\frac{1}{\cano}\right)
\right]
\]
when $\cano\in\R_{>0}$.
A formula of Rogers and Zudilin (\cite[Theorem 2, p.399 and (6), p.386]{RZ})
yields
\[
\frac{\pi^2}{12}\,{}_3F_2\left({\frac{1}{2},\frac{1}{2},\frac{1}{2}\atop \frac{3}{2},1};
\frac{1}{4}\right)=L(E_{24},2)=\frac{\pi^2}{6}L'(E_{24},0).
\]
We thus have
\begin{equation}\label{exp3-eq2}
R_\infty(X_4,\xi)/L'(X_4,0)=-2.
\end{equation}
On the other hand, Theorem \ref{compare-thm} in Appendix B yields
that $\xi$ agrees with the Beilinson-Kato element $4z_{\mathrm{BK},X_4}$.
Since $X_4$ does not have CM, one can apply Brunault's formula
\eqref{remark-BK-eq1} (\cite[Th\'eor\`eme 38 (115)]{Br}).
We thus have
\begin{equation}\label{exp3-eq2-p}
R_{p,\gamma}(X_4,\xi)/L_{p,\gamma}(X_4,\omega^{-1},0)=-2.
\end{equation}
Hence Conjecture \ref{PREC} for the  case $n=0$ and $\operatorname{ord}_p(\gamma )<1$ follows.
\end{pf}
\begin{rem}
\eqref{exp3-eq2-p} also gives an example of \cite[Conjecture 4.30]{New},
\[
(1-p\alpha^{-1})\cF_{\frac12,\frac12}^{(\sigma_4)}(t)|_{t=4}=-L_{p,\alpha}(X_4,\omega^{-1},0)
\]
for $p$ such that $X_4$ has a good ordinary reduction.
\end{rem}

We also discuss the critical slope case (i.e. $\gamma = \beta$, where $\beta$ is a non-unit root with $\operatorname{ord}_p (\beta)=1$) for the elliptic curve $X_4$.
In the same way as in \S \ref{exp2-sect}, one has the table:
\begin{center}
 \begin{tabular}{c|c|c}
 $p$&$\ve^{(\sigma_\cano)}_1(\cano)|_{a=4}$&$\ve^{(\sigma_\cano)}_2(\cano)|_{a=4}$\\
 \hline
$5$&$  14203707162(p^{15})$&$  2689502325(p^{15})$\\
$7$&$3834481697564(p^{15})$&$1564968189963(p^{15})$\\
$11$&$19820018884(p^{10})$&
$ 8395655104 (p^{10})$\\
$13$&$ 114689392802(p^{10})$&$ 107876860898(p^{10})$\\
$17$&
$ 159074202(p^8)$&$ 5132326047(p^8)$
\end{tabular} 
\end{center}
The $p$-adic regulators
\[
R_{p,\beta}(X_4,\xi)
=(1-p\beta^{-1})\frac{\Tr(\reg_\syn(\xi|_{\lambda=4})\cup v_{\beta})}
{\Tr(\omega_{X_4}\cup v_{\beta})}
\]
are as follows:
\begin{center}
 \begin{tabular}{c|c}
 $p$&$R_{p,\beta}(X_4,\xi)$\\
 \hline
$5$&$ 276342201 (p^{13})$\\
$7$&
(supersingular)\\
$11$&
$11^{-1}\cdot 114484571 (p^{8})$\\
$13$&$13^{-1}\cdot 12923433 (p^{8})$\\
$17$&$17^{-1}\cdot   23698690(p^6)$
\end{tabular} 
\end{center}
On the other hand,
$L_{5,\beta}(X_\cano,\omega^{-1},0)=305462$ mod $5^8$ and 
$L_{11,\beta}(E,\omega,0)=1453013467/11$ mod $11^8$ (computed by Robert Pollack).
This indicates
\begin{equation}\label{exp3-eq1}
R_{p,\beta}(X_4,\xi)/L_{p,\beta}(X_4,\omega^{-1},0)=-2
\end{equation}
as desired.


\section{Appendix A: Gauss-Manin connection for elliptic fibrations}\label{Appendix-sect}
For a smooth scheme $Y$ over $T$, we denote by $\Omega^q_{Y/T}=\os{q}{\wedge}_{\O_Y}
\Omega^1_{Y/T}$ the sheaf of relative differential $q$-forms on $Y$ over $T$. 

\medskip

Let $K$ be a field of characteristic zero.
Let $U/S$ be a smooth projective family of curves over the base field $K$.
Let
\[
\nabla:R^1f_*\Omega^\bullet_{U/S}\lra \Omega^1_S\ot
R^1f_*\Omega^\bullet_{U/S}
\]
be the Gauss-Manin connection, where we
simply write $\Omega^q_{S}=\Omega^q_{S/K}$.
This is defined to be the connecting homomorphism
\begin{equation}\label{gm00}
R^1f_*\Omega^\bullet_{U/S}\to R^2f_*(f^*\Omega^1_S\ot
\Omega^{\bullet-1}_{U/S})\cong \Omega^1_S\ot R^2f_*(
\Omega^{\bullet-1}_{U/S})\cong \Omega^1_S\ot R^1f_*\Omega^\bullet_{U/S}
\end{equation}
which arises from an exact sequence
\[
0\lra
f^*\Omega^1_S\ot\Omega^{\bullet-1}_{U/S}
\lra
\bar{\Omega}_U^\bullet
\lra \Omega^\bullet_{U/S}\lra 0,\quad \textup{where} \quad
\bar{\Omega}^\bullet_U:=\Omega^\bullet_U/\operatorname{Im}(f^*\Omega^2_S\ot\Omega^{\bullet-2}_U)
\]
(cf. \cite{h} Ch.III, \S 4).
Here the first isomorphism in \eqref{gm00} is the projection formula,
and the second one is due to the identification 
$R^qf_*\Omega^{\bullet-1}_{U/S}\cong R^{q-1}f_*\Omega^{\bullet}_{U/S}$ with
which we should be careful about ``sign''.
Indeed
the differential of the complex $\Omega^{\bullet-1}_{U/S}$ is ``$-d$''
\[
\Omega^{\bullet-1}_{U/S}:\O_U\os{-d}{\lra} \Omega^1_{U/S}
\quad\mbox{(the first term is placed in degree 1)}
\]
so that we need to arrange
 the sign to make an isomorphism between $R^qf_*\Omega^{\bullet}_{U/S}$
and $R^{q+1}f_*\Omega^{\bullet-1}_{U/S}$.
We make it by a commutative diagram 
\begin{equation}\label{gm-diag}
\xymatrix{
\O_U\ar[r]^d\ar[d]_{-\id}&\Omega^1_{U/S}\ar[d]^{\id}\\
\O_U\ar[r]^{-d}&\Omega^1_{U/S}.
}
\end{equation}
Then $\nabla$ satisfies the usual Leibniz rule
\[
\nabla(gx)=dg\ot x+g\nabla(x),\quad x\in \vg(S,R^1f_*\Omega^\bullet_{U/S}),~g\in\O(S).
\]

\subsection{Family of hyperelliptic curves}\label{expPF-sect}
Let $S$ be an irreducible affine smooth variety over $K$.
Let $f(x)\in \O_S(S)[x]$ be a polynomial of degree $2g+1$ or $2g+2$
which has no multiple roots over any geometric points
$\bar{x}\in S$.
Then it defines a smooth family of hyperelliptic curves $f:U\to S$ 
defined by the Weierstrass equation $y^2=f(x)$. 
To be more precise, let $z=1/x,~u=y/x^{g+1}$ and put $g(z)=z^{2g+2}f(1/z)$.
Let
\[
U_0=\operatorname{Spec} \O_S(S)[x,y]/(y^2-f(x)),\quad
U_\infty=\operatorname{Spec} \O_S(S)[z,u]/(u^2-g(z)).
\]
Then $U$ is obtained by gluing
$U_0$ and $U_\infty$ via identification 
$z=1/x,~u=y/x^{g+1}$.
We assume that there is a section $e:S\to U$.
Then
\begin{equation}\label{he0}
x^{i-1}\frac{dx}{y}=-z^{g-i}\frac{dz}{u},\quad 
\frac{y}{x^i}=\frac{u}{z^{g+1-i}},\quad 1\leq i \leq g.
\end{equation}
We shall compute the Gauss-Manin connection
\begin{equation}\label{gm4}
\nabla:H^1_\dR(U/S)\lra \Omega^1_S\otimes H^1_\dR(U/S),
\quad H^q_\dR(U/S):=H^q_\zar(U,\Omega^\bullet_{U/S})
\end{equation}
(we use the same symbol ``$\Omega^1_S$'' for $\vg(S,\Omega^1_S)$,
since it will be clear from the context which is meant).
To do this, we describe the de Rham cohomology
in terms of the Cech complex.
Write
\[
\check{C}^0({\mathscr F}):=\vg(U_0,{\mathscr F})\op\vg(U_\infty,{\mathscr F}),\quad
\check{C}^1({\mathscr F}):=\vg(U_0\cap U_\infty,{\mathscr F})
\]
for a (Zariski) sheaf $\mathscr F$.
Then the double complex
\[
\xymatrix{
\check{C}^0(\O_U)\ar[r]^d\ar[d]_\delta&
\check{C}^0(\Omega^1_{U/S})\ar[d]^\delta&(x_0,x_\infty)\ar[d]^\delta
\ar[r]^d&(dx_0,dx_\infty)\\
\check{C}^1(\O_U)\ar[r]^d&
\check{C}^1(\Omega^1_{U/S})&x_0-x_\infty
}
\]
gives rise to the total complex
\[
\check{C}^\bullet(U/S):\check{C}^0(\O_U)\os{\delta\times d}{\lra} 
\check{C}^1(\O_U)\times \check{C}^0(\Omega^1_{U/S})
\os{(-d)\times\delta}{\lra} \check{C}^1(\Omega^1_{U/S})
\]
of $\O(S)$-modules starting from degree 0, and the cohomology of it is the de Rham cohomology
$H^\bullet_\dR(U/S)$:
\[
H^q_\dR(U/S)=H^q(\check{C}^\bullet(U/S)),\quad q\geq 0.
\]
Elements of $H^1_\dR(U/S)$ are represented by cocycles
\[
(f)\times (x_0,x_\infty)\quad \mbox{with }df=x_0-x_\infty.
\]
\begin{lem}\label{he-basis}
Suppose
\[
f(x)=a_0+a_1x+\cdots+a_nx^n,\quad a_i\in \O(S)
\]
with $n=2g+1$ or $2g+2$. Put
\begin{equation}\label{he1}
\omega_i:=(0)\times\left(\frac{x^{i-1}dx}{y},-\frac{z^{g-i}dz}{u}\right),
\end{equation}
\begin{equation}\label{he2}
\omega_i^*:=\left(\frac{y}{x^i}\right)\times
\left(\left(\sum_{m>i}(m/2-i)a_mx^{m-i-1}\right)\frac{dx}{y},
\left(\sum_{m\leq i}(m/2-i)a_mz^{g-m+i}\right)\frac{dz}{u}\right)
\end{equation}
for $1\leq i \leq g$.
Then they give a basis of $H^1_\dR(U/S)$.
Moreover \eqref{he1} span the image of $\vg(U,\Omega^1_{U/S})\hra H^1_\dR(U/S)$.
\end{lem}
\begin{pf}
Exercise.
\end{pf}
\begin{lem}\label{cech-equiv}
There are the following equivalence relations:
\[
(x^iy^j)\times(0,0)\equiv(0)\times(-d(x^iy^j),0)\mod \operatorname{Im}\check{C}^0(\O_U),
\]
\[
(z^iu^j)\times(0,0)\equiv(0)\times(0,d(z^iu^j))\mod \operatorname{Im}\check{C}^0(\O_U).
\]
\end{lem}
\begin{pf}
Straightforward from the definition .
\end{pf}

\subsection{Computation of Gauss-Manin connection}
Let us compute $\nabla(\omega_i)$ and $\nabla(\omega_i^*)$.
Recall that there is the exact sequence
\begin{equation}\label{gm0}
0\lra \check{C}^\bullet(f^*\Omega^1_S\ot\Omega^{\bullet-1}_{U/S})\lra
\check{C}^\bullet(\bar{\Omega}_U^\bullet)
\lra \check{C}^\bullet(\Omega^\bullet_{U/S})\lra 0
\end{equation}
and it gives rise to the connecting homomorphism 
\[
\delta:H^1_\dR(U/S)=H^1(\check{C}^\bullet(\Omega^\bullet_{U/S}))\lra H^2(U,
f^*\Omega^1_S\ot\Omega^{\bullet-1}_{U/S})
=H^2(\check{C}^\bullet(f^*\Omega^1_S\ot\Omega^\bullet_{U/S})).
\]
Recall the isomorphism
\[
\check{C}^\bullet(f^*\Omega^1_S\ot\Omega^{\bullet-1}_{U/S})
\os{\cong}{\lra}
\Omega^1_S\ot\check{C}^\bullet(\Omega^{\bullet}_{U/S})
\]
induced from \eqref{gm-diag}. It induces the isomorphism 
\[
\iota:H^2(U,f^*\Omega^1_S\ot\Omega^{\bullet-1}_{U/S})\os{\cong}{\lra}
\Omega^1_S\ot H^1_\dR(U/S).
\]
By definition we have $\nabla=\iota\delta$, where $\nabla$ is
the Gauss-Manin connection \eqref{gm4}.
Let us write down the maps $\delta$ and $\iota$ in terms of Cech cocycles.
The differential operator $\cD$ on the total complex of the middle term of \eqref{gm0} is 
given as follows
\[\cD:
\check{C}^1(\O_U)\times \check{C}^0(\bar{\Omega}^1_{U})
\lra \check{C}^1(\bar{\Omega}^1_{U})\times \check{C}^0(\bar{\Omega}^2_{U}),
\]
\[
(\alpha)\times(\beta_0,\beta_\infty)\longmapsto
(-d\alpha+\beta_0-\beta_\infty)\times(d\beta_0,d\beta_\infty).
\]
We denote a lifting of $(z_0,z_\infty)\in \check{C}^0({\Omega}^1_{U/S})$
by $(\hat{z}_0,\hat{z}_\infty)\in \check{C}^0(\bar{\Omega}^1_U)$.
Then for $(\alpha)\times (z_0,z_\infty)\in H^1_\dR(U/S)$ one has
\begin{align}
(\alpha)\times (z_0,z_\infty)&\os{\delta}{\longmapsto} 
\cD((\alpha)\times (\hat{z}_0,\hat{z}_\infty))\\
&=(-d\alpha+\hat{z}_0-\hat{z}_\infty)\times(d\hat{z}_0,d\hat{z}_\infty)\\
&\in \check{C}^1(f^*\Omega^1_S)\times \check{C}^0(f^*\Omega^1_S\ot\Omega^1_{U/S}).
\end{align}
The isomorphism $\iota$ is given by
\begin{equation}\label{gm3}
(gdt)\times(dt\wedge z_0,dt \wedge z_\infty )
\os{\iota}{\longmapsto}
dt\ot [(-g)\times(z_0,z_\infty)]
\end{equation}
(the ``sign'' appears in the above due to \eqref{gm-diag}).

To compute $\nabla(\omega_i)$ and $\nabla(\omega_i^*)$ for the basis
in Lemma \ref{he-basis}, there remains to compute lifting of $dx/y$ and $dz/u$.

\begin{lem}\label{ABCD}
Let $A,B\in\O(S)[x]$ and $C,D\in \O(S)[z]$ satisfying
\[
Af+B\frac{\partial f}{\partial x}=1,\quad Cg+D\frac{\partial g}{\partial z}=1.
\]
Put differential 1-forms
\[
\widehat{\frac{dx}{y}}:=\frac{Afdx+Bdf}{y}=Aydx+2Bdy
\in\vg(U_0,\Omega^1_U),
\]
\[
\widehat{\frac{dz}{u}}:=\frac{Cgdz+Ddg}{u}=Cudz+2Ddu\in\vg(U_\infty,\Omega^1_U).
\]
Then
\[
x^i\widehat{\frac{dx}{y}}\in\vg(U_0,\Omega^1_U),\quad
z^i\widehat{\frac{dz}{u}}\in\vg(U_\infty,\Omega^1_U)
\]
are liftings of 
$x^idx/y\in\vg(U_0,\Omega^1_{U/S})$ and
$z^idz/u\in\vg(U_\infty,\Omega^1_{U/S})$ respectively.
\end{lem}
\begin{pf}
Straightforward.
\end{pf}

By using the liftings in Lemma \ref{ABCD}, one can compute 
the map $\delta$.
With use of Lemma \ref{cech-equiv}, 
one finally obtains the connection matrix of $\nabla$.

\medskip

Here is an explicit formula in case of elliptic fibration.

\begin{thm}\label{A1}
Let $S$ be a smooth affine curve and $f:U\to S$ a projective smooth family of
elliptic curves whose affine form is given by a Weierstrass equation
$y^2=4x^3-g_2x-g_3$ with $\Delta:=g_2^3-27g_3^2\in \O_S(S)^\times$.
Suppose that $\Omega^1_S$ is a free $\O_S$-module with a base 
$d\lambda\in \vg(S,\Omega^1_{S})$.
For $f\in \O_S(S)$, we define $f'\in\O_S(S)$ by $df=f'd\lambda$.
Let
\begin{equation}\label{A1-1}
\omega:=(0)\times (\frac{dx}{y},-\frac{dz}{u})
\end{equation}
\begin{equation}\label{A1-2}
\eta:=(\frac{y}{2x})\times (\frac{xdx}{y},\frac{(g_2z+2g_3z^2)dz}{4u})
\end{equation}
be elements in $H^1_\dR(U/S)$. 
Then we have
\[
\nabla\omega=-\frac{\Delta^\prime}{12\Delta}d\lambda\ot\omega+\frac{3(\EE)}{2\Delta}
d\lambda\ot\eta ,
\]
\[
\nabla\eta=-\frac{g_2(\EE)}{8\Delta}d\lambda\ot\omega
+\frac{\Delta^\prime}{12\Delta}d\lambda\ot\eta.
\]
\end{thm}

\begin{prop}\label{A-sympl-1}
Let $g_2,g_3\in K[[\lambda]]$ be power series such that
$g_2(0)\ne0$, $g_3(0)\ne0$ and $\Delta=g_2^3-27g_3^2\in \lambda K[[\lambda]]$.
Let $f:U\to S=\operatorname{Spec} K[[\lambda]]$ be the family of elliptic curves given by an affine equation
\[
y^2=4x^3-g_2x-g_3
\]
with semistable reduction at $\lambda=0$.
Let $(12g_2)^{-1/4}=c_0+c_1\lambda+\cdots\in \ol{K}[[\lambda]]$ denote the power series 
such that the initial term satisfies
$c_0^2=-18^{-1}g_2(0)g_3(0)^{-1}$.
Let $F(\lambda)$ be defined by
\[
\rho^*\omega=F(\lambda)\frac{du}{u}.
\]
Then
\begin{equation}\label{S-B-formula}
F(\lambda)=(12g_2)^{-\frac{1}{4}}{}_2F_1\left({\frac{1}{12},\frac{5}{12}\atop 1};J^{-1}\right),
\quad \text{where } J:=g_2^3/(g_2^3-27g_3^2).
\end{equation}
\end{prop}
\begin{pf}
We may assume $K=\C$. Denote by $U_\lambda^{an}=f^{-1}(\lambda)$ the analytic torus for $|\lambda |\ll1$. 
Let $\delta_\lambda\in H_1(U_\lambda^{an},\Z)$ be the vanishing cycle. 
Then
\[
F(\lambda)=\frac{1}{2\pi i}\int_{\delta_\lambda}\omega.
\]
This is the unique solution of the Picard-Fuchs equation \cite[(1.3)]{S-B},
and then \eqref{S-B-formula} is proven in \cite[(1.5)]{S-B}.
\end{pf}
\begin{prop}\label{A-sympl-2}
Let the notation be as in Proposition \ref{A-sympl-1}.
Put 
\[
H(\lambda):=\frac{2\Delta}{3(\EE)}\left(F'(\lambda)+\frac{\Delta'}{12\Delta}F(\lambda)\right)
\]
and 
\begin{equation}\label{S-B-formula-1}
\wh\omega:=\frac{1}{F(\lambda)}\omega,\quad\wh\eta:=-H(\lambda)\omega+F(\lambda)\eta
\end{equation}
\[
\left(\Longleftrightarrow\, \omega=F(\lambda)\wh\omega,\quad \eta=
H(\lambda)\wh\omega+\frac{1}{F(\lambda)}\wh\eta\right) .
\]
Then $\{\wh\omega,\wh\eta\}$ forms a de Rham symplectic basis which satisfies
\[
\nabla(\wh\omega)=\frac{dq}{q}\ot\wh\eta ,
\]
where $q=\kappa\lambda^n+\cdots\in \lambda K[[\lambda]]$ is the multiplicative period.
\end{prop}
\begin{pf}
By \eqref{A1-1} and \eqref{A1-2},
it is straightforward to see that 
\[
\nabla(\wh\eta)=0, \quad
\nabla(\wh\omega)=\frac{3(\EE)}{2\Delta F(\lambda)^2}d\lambda\ot\wh\eta.
\]
Since 
\[
\nabla(\omega\cup\eta)=\nabla(\omega)\cup\eta+\omega\cup\nabla(\eta)=0,
\]
$\wh\omega\cup\wh\eta=C\ne0$ is a nonzero constant.
Hence $\{\wh\omega,C^{-1}\wh\eta\}$ forms a de Rham symplectic basis.
It follows from Proposition \ref{local-GM} that one has
\begin{equation}\label{A-sympl-eq1}
\frac{dq}{q}=C\frac{1}{F(\lambda)^2}\frac{3(\EE)}{2\Delta}d\lambda=
-C\frac{3g_2g_3}{2F(\lambda)^2\Delta}\left(3\frac{g'_2}{g_2}-2\frac{g'_3}{g_3}\right)d\lambda.
\end{equation}
Recall the formula \eqref{S-B-formula}
\[
F(\lambda)=(12g_2)^{-\frac{1}{4}}{}_2F_1\left({\frac{1}{12},\frac{5}{12}\atop 1};J^{-1}\right)
=c_0+c_1\lambda+\cdots,
\quad J:=g_2^3/(g_2^3-27g_3^2)
\]
with $c_0^2=-18^{-1}g_2(0)g_3(0)^{-1}$.
Take the residue of \eqref{A-sympl-eq1} at $\lambda=0$.
Since
\[
\dlog(g_3^{-2}\Delta+27)=\dlog\frac{g_2^3}{g_3^2}=
\left(3\frac{g'_2}{g_2}-2\frac{g'_3}{g_3}\right)d\lambda
=\left(\frac{cng_3(0)^{-2}}{27}\lambda^{n-1}+(\text{higher terms)}\right)d\lambda
\]
where $\Delta=c\lambda^n+\cdots$ with $c\ne0$,
one has
\[
-C\frac{3g_2(0)g_3(0)}{2c_0^2}\frac{ng_3(0)^{-2}}{27}=n
\quad\Longleftrightarrow\quad
C=1.
\]
\end{pf}
In the above proof, we showed the following.
\begin{prop}\label{A-sympl-3}
Let $\tau(\lambda)\in \lambda K[[\lambda]]$ be defined by $q=\kappa \lambda^n\exp(\tau(\lambda))$. Then
\begin{equation}\label{A-sympl-3-eq1}
\frac{dq}{q}=
\frac{3(\EE)}{2\Delta F(\lambda)^2}d\lambda,\quad
\frac{d\tau(\lambda)}{d\lambda}=
\frac{3(\EE)}{2\Delta F(\lambda)^2}-\frac{n}{\lambda}.
\end{equation}
\end{prop}
\subsection{Deligne's canonical extension}\label{de-sect} 
Let $\ol{S}$ be a smooth curve over $\C$, and $j:S:=\ol S-\{P\}\hra \ol S$ an inclusion. 
Let $(\cH,\nabla)$ be a vector bundle with integrable connection
over $S$. 
Then there is unique subbundle $\cH_e\subset j^{an}_*\cH$ satisfying the following conditions
(cf. \cite{zucker} (17)):
\begin{itemize}
\item
The connection extends to a map
$\nabla:\cH_e\to\Omega^1_{\ol S}(\log P)\ot\cH_e$,
\item
each eigenvalue $\alpha$ of $\Res_P(\nabla)$ satisfies $0\leq \mathrm{Re}(\alpha)<1$.
\end{itemize}
The extended bundle $(\cH_e,\nabla)$ is called {\it Deligne's canonical extension}.
The inclusion map 
\[
[\cH_e\os{\nabla}{\to}\Omega^1_{\ol S}(\log P)\ot\cH_e]\lra
[j_*\cH\os{\nabla}{\to}\Omega^1_{S}\ot j_*\cH]
\]
is a quasi-isomorphism of complexes of sheaves. 
Besides $\exp(-2\pi i\Res_P(\nabla))$ coincides with the monodromy
operator on $H_\C=\ker(\nabla^{\mathit an})$ around $P$
(cf. \cite{steenbrink}, (2.21)).

\medskip

Let $f:\ol U\to \ol S$ be a projective flat morphism of nonsingular algebraic varieties such that
$f$ is smooth over $S$ and $D=f^{-1}(P)$ is a NCD in $\ol U$.
Let $U=f^{-1}(S)$.
Let $(\cH,\nabla)$ be the higher direct image $(R^if_*\Omega^\bullet_{U/S},\nabla)$.
Then all eigenvalues of $\Res_P(\nabla)$ are rational numbers.
There is the natural isomorphism
\begin{equation}\label{A2-eq1}
R^if_*\Omega^\bullet_{\ol U/\ol S}(\log D)\cong \cH_e.
\end{equation}
In particular, the left hand side is a locally free $\O_{\ol S}$-module.
If $f$ is a morphism of smooth $K$-schemes, then Deligne's extension 
$(R^if_*\Omega^\bullet_{\ol U/\ol S}(\log D),\nabla)$ can be obviously defined over $K$.

\medskip

We have seen how to compute a connection matrix of the Gauss-Manin connection
for a family of hyperelliptic curves.
In case of an elliptic fibration, they are simply given as follows.
\begin{prop}\label{A2}
Let $f:\cE\to \operatorname{Spec} K[[\lambda]]$ be an elliptic fibration defined by a minimal
Weierstrass equation $y^2=4x^3-g_2x-g_3$ with $g_2,g_3\in K[[\lambda]]$, $\Delta:=g_2^3-27g_3^2
\ne0$.
Then $\cH_e$ has a basis $\{\omega,\eta\}$ (resp. $\{\lambda\omega,\eta\}$)
if $f$ has a semistable or smooth reduction (resp. additive reduction).
\end{prop}
\begin{pf}
Since we have Theorem \ref{A1},
we can show the above by case-by-case analysis based on $(\ord(g_2),\ord(g_3))$.
The detail is left to the reader.
\end{pf}

\subsection{The de Rham cohomology of elliptic curves over a ring}\label{int-DR-sect}
Let $U$ be a smooth scheme over a commutative ring $A$.
For a point $s$ of $\operatorname{Spec} A$, we denote the residue field by $k(s)$, and write $U_s:=
U\times_A k(s)$.
We define the de Rham cohomology $H^\bullet_\dR(U/A)$ in the same way as before:
\[
H^i_\dR(U/A):={\mathbb H}^i_\zar(U,\Omega^\bullet_{U/A}),\quad i\in \Z_{\geq0}.
\] 
\begin{lem}\label{A3-lem}
Let $A$ be a noetherian integral domain.
Let $U\to\operatorname{Spec} A$ be a smooth projective morphism of relative dimension one
with geometrically connected fibers. 
Then, for any $i,j,l\geq 0$,
 $H^i(U,\Omega^j_{U/A})$ and $H^l_\dR(U/A)$ are locally free $A$-modules
and there are the canonical isomorphisms
\begin{equation}\label{A3-lem-eq2}
H^i(U,\Omega^j_{U/A})\ot_A k(s)\os{\cong}{\lra}H^i(U_s,\Omega^j_{U_s/k(s)}),
\end{equation}
\begin{equation}\label{A3-lem-eq3}
H^l_\dR(U/A)\ot_A k(s)\os{\cong}{\lra} H^l_\dR(U_s/k(s))
\end{equation}
for all points $s\in \operatorname{Spec} A$.
Moreover there is a (non-canonical) isomorphism
\begin{equation}\label{A3-lem eq4}
H^1_\dR(U/A)\cong H^0(\Omega^1_{U/A})\op H^1(\O_U).
\end{equation}
\end{lem}
\begin{pf}
Since $\dim_{k(s)}H^i(\Omega^j_{U_s/k(s)})$ is constant with respect to $s$,
one can apply \cite[III,12.9]{Ha}, so that 
$H^i(U,\Omega^j_{U/A})$ is a locally free
$A$-module and the isomorphism \eqref{A3-lem-eq2} follows.
The natural map $H^1(\Omega^1_{U/A})\to H^2_\dR(U/A)$ is surjective as $H^2(\O_U)=0$.
There are the trace maps $\Tr:H^1(\Omega^1_{U_s/k(s)})\to k(s)$ (cf. \cite[VII,4.1]{RD})
and $\Tr:H^2_\dR(U_s/k(s))\to k(s)$ (cf. \cite[Proposition (2.2)]{h}) which are compatible.
This implies that 
the map $H^1(\Omega^1_{U_s/k(s)})\to H^2_\dR(U_s/k(s))$ is a non-zero
surjective map, and hence bijective as $\dim_{k(s)}H^1(\Omega^1_{U_s/k(s)})=1$.
Let
\[
\xymatrix{
H^1(\Omega^1_{U/A})\ot_Ak(s)\ar[r]^{\text{surj.}}\ar[d]_\cong&
H^2_\dR(U/A)\ot_Ak(s)\ar[d]\\
H^1(\Omega^1_{U_s/k(s)})\ar[r]^\cong& H^2_\dR(U_s/k(s))
}
\]
be a commutative diagram.  Hence all arrows are bijective and all terms are one-dimensional for all points $s\in \operatorname{Spec} A$. 
We thus have \eqref{A3-lem-eq3} in case $l=2$ and that
$H^2_\dR(U/A)$ is locally free of rank $1$, since
$A$ is a noetherian integral domain (\cite[II, 8.9]{Ha}).
Moreover we have an isomorphism $H^1(\Omega^1_{U/A})\os{\cong}{\to} H^2_\dR(U/A)$ by Nakayama's lemma, so that we have an exact sequence
\begin{equation}\label{A3-lem-eq5}
\xymatrix{
0\ar[r]&H^0(\Omega^1_{U/A})\ar[r]&H^1_\dR(U/A)\ar[r]
&H^1(\O_U)\ar[r]&0.
}
\end{equation}
Since the right term is a projective $A$-module, the sequence splits.
Hence 
one has a (non-canonical) isomorphism \eqref{A3-lem eq4}, and also
that $H^1_\dR(U/A)$ is a locally free $A$-module.
Finally, the isomorphism \eqref{A3-lem-eq3} in case $l=1$ follows from a commutative
diagram 
\[
\xymatrix{
0\ar[r]&H^0(\Omega^1_{U/A})\ot k(s)\ar[d]_\cong\ar[r]&H^1_\dR(U/A)\ot k(s)\ar[r]\ar[d]
&H^1(\O_U)\ot k(s)\ar[r]\ar[d]^\cong&0\\
0\ar[r]&H^0(\Omega^1_{U_s/k(s)})\ar[r]&H^1_\dR(U_s/k(s))\ar[r]
&H^1(\O_{U_s})\ar[r]&0,
}
\]
where the right and left isomorphisms follow from \eqref{A3-lem-eq2}.
\end{pf}

Let $W$ be a commutative ring of characteristic zero, and let $f:U\to S$ be a projective smooth morphism
of smooth $W$-schemes.
Then one can define the Gauss-Manin connection
\[
\nabla:H^i_\dR(U/S)\lra \Omega^1_S\ot H^i_\dR(U/S)
\]
in the same way as before, where $\Omega^\bullet_S=\Omega^\bullet_{S/W}$.
Moreover if $f$ extends to a projective flat family $\ol f:\ol U\to \ol S=S\cup\{P\}$ such that 
$P$ is a $W$-rational point and
$D=f^{-1}(P)$ is a relative NCD, then one can also define Deligne's canonical extension
\[
(R^if_*\Omega^\bullet_{\ol U/\ol S}(\log D),\nabla).
\]
The sheaf $R^if_*\Omega^\bullet_{\ol U/\ol S}(\log D)$ is a coherent $\O_{\ol S}$-module,
but not necessarily a locally free $\O_{\ol S}$-module.
It seems difficult to ask whether it is locally free or not in a general
situation.
The following gives a partial answer in case of elliptic fibrations.
\begin{prop}\label{A3}
Suppose that $W$ is a noetherian integral domain of characteristic zero
such that $2$ is invertible in $W$. Put $K:=\operatorname{Frac} W$.
Let $f:\cE\to \operatorname{Spec} W[[\lambda]]$ be an elliptic fibration given by a Weierstrass equation
$y^2=4x^3-g_2x-g_3$ with $g_2,g_3\in W[[\lambda]]$, $\Delta:=g_2^3-27g_3^2\in W((\lambda))^\times$,
such that the central fiber $D=f^{-1}(0)$ is a relative NCD.
Suppose that $y^2=4x^3-g_2x-g_3$ is minimal as a Weierstrass equation over $K[[\lambda]]$.
Let $E:=\cE\times_{W[[\lambda]]} W((\lambda))$ and put $\omega:=dx/y$ and $\eta:=xdx/y$.
Then
\begin{equation}\label{A3-eq0}
H^1_\dR(E/W((\lambda)))= W((\lambda))\omega+W((\lambda))\eta
\end{equation}
is a free $W((\lambda))$-module of rank $2$.
Put
\[
H:=\operatorname{Im} [H^1(\cE,\Omega^\bullet_{\cE/W[[\lambda]]}(\log D))\to 
H^1_\dR(E/W((\lambda)))].
\]
Let $\omega_0=\omega=dx/y$ (resp. $\omega_0=\lambda\omega$)
if $f$ has a smooth or multiplicative reduction
(resp. an additive reduction).
Then
\begin{equation}\label{A3-eq1}
H\subset W[[\lambda]]\omega_0+W[[\lambda]]\eta.
\end{equation}
The equality holds if the multiplicity of any component of $D$
is invertible in $W$.
\end{prop}
\begin{pf}
Write $E_s:=E\times_{W((\lambda))}k(s)$ for $s\in \operatorname{Spec} W((\lambda))$.
Then $\{dx/y|_{E_s},xdx/y|_{E_s}\}$ is a $k(s)$-basis of $H^1_\dR(E_s/k(s))$.
Therefore it follows from Lemma \ref{A3-lem}
\eqref{A3-lem-eq3} and Nakayama's lemma that $H^1_\dR(E/W((\lambda)))$ is a 
free $W((\lambda))$-module with basis $\{\omega,\eta\}$.
Next we show \eqref{A3-eq1}.
Write $\cE_K:=\cE\times_{W[[\lambda]]}K[[\lambda]]$, $E_K:=E\times_{W[[\lambda]]}K[[\lambda]]$ etc. Put
\[
H_K:=\Image [H^1(\cE_K,\Omega^\bullet_{\cE_K/K[[\lambda]]}(\log D))\to 
H^1_\dR(E_K/K((\lambda)))].
\]
By Proposition \ref{A2},
$H_K$ is a free $K[[\lambda]]$-module with basis $\{\omega_0,\eta\}$.
Hence
\[
H\subset \bigg(W((\lambda))\omega_0+W((\lambda))\eta\bigg)\cap H_K
=W[[\lambda]]\omega_0+W[[\lambda]]\eta
\]
as required.
The last statement is more delicate.
Recall the descriptions \eqref{A1-1} and \eqref{A1-2} 
\[
\omega_0=(0)\times (\lambda^e\frac{dx}{y},-\lambda^e\frac{dz}{u}),\quad
\eta=(\frac{y}{2x})\times (\frac{xdx}{y},\frac{(g_2z+2g_3z^2)dz}{4u})\in H^1_\dR(\cE/W((\lambda)))
\]
in terms of Cech cocycles,
where we put $e=0$ if $D$ is multiplicative and $e=1$ if additive. 
Let $U^*_0:=\operatorname{Spec} W[[\lambda]][x,y]/(y^2-4x^3+g_2x+g_3)$ and 
$U^*_\infty:=\operatorname{Spec} W[[\lambda]][z,u]/(u^2-4z+g_2z^3+g_3z^4)$, and
$\rho:\cE\to \cE^*:=U^*_0\cup U^*_\infty$ the blow-ups. Let $U_0:=\rho^{-1}(U^*_0)$
and $U_\infty:=\rho^{-1}(U^*_\infty)$.
To prove $\omega_0,\eta\in H$, we show 
that the above cocycles belong to
\[\vg(U_0\cap U_\infty,\O_{\cE})\times
\vg(U_0,\Omega^1_{\cE/W[[\lambda]]}(\log D))\op \vg( U_\infty,\Omega^1_{\cE/W[[\lambda]]}(\log D)).
\]
Since $y/2x\in \vg(U_0\cap U_\infty,\O_\cE)$, it is enough to show that
\begin{equation}\label{A3-d}
\lambda^e\frac{dx}{y}\in \vg(U_0,\Omega^1_{\cE/W[[\lambda]]}(\log D)),\quad
\lambda^e\frac{dz}{u}\in \vg(U_\infty,\Omega^1_{\cE/W[[\lambda]]}(\log D)),
\end{equation}
and
\begin{equation}\label{A3-d-1}
\frac{xdx}{y}\in \vg(U_0,\Omega^1_{\cE/W[[\lambda]]}(\log D)),\quad
\frac{(g_2z+2g_3z^2)dz}{u}\in \vg(U_\infty,\Omega^1_{\cE/W[[\lambda]]}(\log D)).
\end{equation}
We first show \eqref{A3-d}.
We know that
\begin{equation}\label{A3-eq2}
\lambda^e\frac{dx}{y}\in \vg(U_0,\Omega^1_{\cE_K/K[[\lambda]]}(\log D)),\quad
\lambda^e\frac{dz}{u}\in \vg(U_\infty,\Omega^1_{\cE_K/K[[\lambda]]}(\log D))
\end{equation}
by Proposition \ref{A2}.
On the other hand, since $\vg(\cE,\Omega^1_{\cE/W((\lambda))})$
is a free $W((\lambda))$-module with basis $dx/y=-dz/u$, one has
\begin{equation}\label{A3-eq3}
\lambda^{e+m}\frac{dx}{y}\in \vg(U_0,\Omega^1_{\cE/W[[\lambda]]}(\log D)),\quad
\lambda^{e+m}\frac{dz}{u}\in \vg(U_\infty,\Omega^1_{\cE/W[[\lambda]]}(\log D))
\end{equation}
for some large integer $m$.
By the assumption, the multiplicity of any component of $D$
is invertible in $W$. One sees that
$\Omega^1_{\cE/W[[\lambda]]}(\log D)$ is locally free as $f$ is locally given by $\lambda\mapsto x_1^{r_1}x_2^{r_2}$, where $r_i$ denotes the multiplicity.
One can check that the map $a$ in the following diagram is injective:
\[
\xymatrix{
0\ar[r]&\Omega^1_{\cE/W[[\lambda]]}(\log D)\ar[r]^{\lambda^m}\ar[d]
&\Omega^1_{\cE/W[[\lambda]]}(\log D)\ar[r]\ar[d]&
\Omega^1_{\cE/W[[\lambda]]}(\log D)/\lambda^m\ar[r]\ar[d]^a&0\\
0\ar[r]&\Omega^1_{\cE_K/K[[\lambda]]}(\log D)\ar[r]^{\lambda^m}&\Omega^1_{\cE_K/K[[\lambda]]}(\log D)\ar[r]&
\Omega^1_{\cE_K/K[[\lambda]]}(\log D)/\lambda^m\ar[r]&0.
}
\]
Now the diagram chase with use of \eqref{A3-eq2} and \eqref{A3-eq3} implies \eqref{A3-d}.

Next we show \eqref{A3-d-1}.
It follows from \eqref{A3-d} that
\[
\lambda^e\frac{xdx}{y}\in \vg(U_0,\Omega^1_{\cE/W[[\lambda]]}(\log D)),\quad
\lambda^e\frac{(g_2z+2g_3z^2)dz}{u}\in \vg(U_\infty,\Omega^1_{\cE/W[[\lambda]]}(\log D)).
\]
If we show
\[
\frac{xdx}{y}\in \vg(U_{0,K},\Omega^1_{\cE_K/K[[\lambda]]}(\log D_K)),\quad
\frac{(g_2z+2g_3z^2)dz}{u}\in \vg(U_{\infty,K},\Omega^1_{\cE_K/K[[\lambda]]}(\log D_K)),
\]
then, in the same way as above, one can show \eqref{A3-d-1} by the diagram chase, which completes the proof of Proposition \ref{A3}.
Since $U_\infty\to\operatorname{Spec} W[[\lambda]]$ is smooth, 
one has $dz/u\in \vg(U_{\infty},\Omega^1_{\cE/W[[\lambda]]})$, and hence the latter follows.
The rest is to show
\begin{equation}\label{A3-eq3-1}
\frac{xdx}{y}\in \vg(U_{0,K},\Omega^1_{\cE_K/K[[\lambda]]}(\log D_K)).
\end{equation}
If $e=0$, there is nothing to prove. Suppose $e=1$, namely $D_K$ is additive.
Let $T_K$ be the singular points of $D_K$ and put $\cE^\circ_K:=\cE_K\setminus T_K$.
It is enough to show
\begin{equation}\label{A3-eq3-2}
\frac{xdx}{y}\in \vg(U_{0,K}\setminus T_K,\Omega^1_{\cE_K^\circ/K[[\lambda]]}(\log D_K))
\end{equation}
as $\operatorname{codim}(T_K)=2$.
Let $\operatorname{Spec} K[[\lambda_0]]\to \operatorname{Spec} K[[\lambda]]$ be given by $\lambda_0^{n}=\lambda$ and
consider a cartesian diagram
\[
\xymatrix{
\wt\cE_{0,K}^\circ\ar[r]^j
\ar[dr]_{\wt f}&\cE^\circ_{0,K}\ar[r]\ar[d]\ar@{}[rd]|{\square}
&\cE_K^\circ\ar[d]^f\\
&K[[\lambda_0]]\ar[r]&K[[\lambda]]
}\]
with $j$ the normalization.
The morphism
$f$ is locally given by $(u,v)\mapsto \lambda=u^{n_i}$ around a component of the central fiber,
where $n_i$ is the multiplicity.
Take $n=\mathrm{lcm}(n_i)_i$.
Then 
$\cE_{0,K}^\circ$ is locally defined by $u^{n_i}=\lambda^n$ in $(u,v,\lambda)$-space,
and hence $\wt\cE_{0,K}^\circ$ is regular and $\wt f$ is smooth.
Let $\wt U_{0,K}\subset \wt\cE^\circ_{0,K}$ be the inverse image of $U_{0,K}$,
and $\wt D_K$ the central fiber of $\wt f$.
Then, to show \eqref{A3-eq3-2}, it is enough to show
\[
\frac{xdx}{y}\in \vg(\wt U_{0,K},\Omega^1_{\wt\cE^\circ_{0,K}/K[[\lambda_0]]}(\log \wt D_K)).
\]
Let $x_0=x/\lambda_0^m$, $y_0=y/\lambda_0^l$, $g_{2,0}=g_2/\lambda_0^{2m}$ and 
$g_{3,0}=g_3/\lambda_0^{3m}$ with $3m=2l$ such that
$y_0^2=4x_0^3-g_{2,0}x_0-g_{3,0}$ is minimal over $K[[\lambda_0]]$.
Since the Kodaira type of $\wt f$ is I$_s$ with $s\geq0$, 
one has
\[
\frac{xdx}{y}=\lambda_0^{2m-l}\frac{x_0dx_0}{y_0}\in 
\vg(\wt U_{0,K},
\Omega^1_{\wt\cE_{0,K}/K[[\lambda_0]]})\subset
\vg(\wt U_{0,K},
\Omega^1_{\wt\cE_{0,K}/K[[\lambda_0]]}(\log \wt D_K))
\]
as $2m-l>0$. This completes the proof.
\end{pf}


\newcommand{\N}{\mathbf{N}}
\newcommand{\Zhat}{\hat{\Z}}
\newcommand{\Qb}{\overline{\Q}}
\newcommand{\kb}{\overline{k}}
\newcommand{\h}{\mathcal{H}}
\newcommand{\A}{\mathbf{A}}
\newcommand{\p}{\mathbf{P}}

\section{Appendix B : Comparing elements in $K_2$ of elliptic curves}\label{AppendixB-sect}
\begin{center}
{\large by\, Fran\c{c}ois Brunault}
\end{center}

In this appendix, we use results of Goncharov and Levin \cite{goncharov-levin} to compare the element $\xi$ in \eqref{RZ-eq2}
and the Beilinson--Kato element in $K_2$ of the elliptic curve $X_0(24)$.
This comparison is used in the proof of Theorem 5.2.

\subsection{Describing $K_2$ of elliptic curves}

Goncharov and Levin gave an explicit description of $K_2$ of an elliptic curve using the so-called elliptic Bloch group. In this section we recall this construction and state the result needed later to compare the two elements in $K_2$ (see Theorem \ref{GL thm beta 2}).

Let $E$ be an elliptic curve defined over a number field $k$. To describe Quillen's $K$-group $K_2(E)$, our starting point is the localisation map $K_2(E) \to K_2(k(E))$, where $k(E)$ is the function field of $E$. The group $K_2(k(E))$ can be described using Matsumoto's theorem: for any field $F$, we have an isomorphism
\begin{equation*}
K_2(F) \cong \frac{F^\times \otimes_{\Z} F^\times}{\langle x \otimes (1-x) : x \in F \setminus \{0,1\}\rangle}.
\end{equation*}
The class of $x \otimes y$ in $K_2(F)$ is denoted by $\{x,y\}$ and is called a Milnor symbol. The relations $\{x,1-x\}=0$ are called the Steinberg relations.

Let $\Z[E(\kb)]$ be the group algebra of $E(\kb)$. Consider the Bloch map
\begin{align*}
\beta : \kb(E)^\times \times \kb(E)^\times & \longrightarrow \Z[E(\kb)] \\
(f , g) & \longmapsto \sum_{i,j} m_i n_j (p_i - q_j),
\end{align*}
where $\mathrm{div}(f) = \sum_i m_i (p_i)$ and $ \mathrm{div}(g) = \sum_j n_j (q_j)$ are the divisors of $f$ and $g$. The map $\beta$ is bilinear, so it induces a map
\[
\xymatrix{
\kb(E)^\times \otimes_{\Z} \kb(E)^\times \ar[r] & \Z[E(\kb)],
}
\]
which we still denote by $\beta$.

Let $I$ be the augmentation ideal of $\Z[E(\kb)]$. The group $P$ of principal divisors on $E/\kb$ is generated by the divisors
\begin{equation*}
(p+q)-(p)-(q)+(0) = \bigl((p)-(0)\bigr) \cdot \bigl((q)-(0)\bigr) \qquad \textrm{with } p,q \in E(\kb),
\end{equation*}
so we have $P=I^2$ and $I/I^2 \cong E(\kb)$. It follows that $\beta$ takes values in $I^4$, and the image of $\beta$ generates $I^4$. Following Goncharov and Levin, we now define the elliptic Bloch group of $E$.

\begin{defn} Let $R_3^*(E/\kb)$ be the subgroup of $\Z[E(\kb)]$ generated by the divisors $\beta(f, 1-f)$ with $f \in \kb(E)$, $f \neq 0,1$. The elliptic Bloch group of $E/\kb$ is
\begin{equation*}
B_3^*(E/\kb) = \frac{I^4}{R_3^*(E/\kb)}.
\end{equation*}
The elliptic Bloch group of $E/k$ is
\begin{equation*}
B_3^*(E) = B_3^*(E/\kb)^{\Gal(\kb/k)}.
\end{equation*}
\end{defn}

By Matsumoto's theorem and the definition of the elliptic Bloch group, the map $\beta$ gives rise to a commutative diagram
\[
\xymatrix{
K_2(\kb(E)) \ar[r] & B_3^*(E/\kb) \\
K_2(k(E)) \ar[r] \ar[u] & B_3^*(E). \ar[u]
}
\]

Goncharov and Levin \cite{goncharov-levin} proved the following fundamental result.

\begin{thm} \label{GL thm beta}
The composite map
\[
\xymatrix{
K_2(E) \otimes \Q \ar[r] & K_2(k(E)) \otimes \Q \ar[r] & B_3^*(E) \otimes \Q
}
\]
is injective.
\end{thm}

In fact, Goncharov and Levin showed that modulo torsion, $K_2(E)$ is isomorphic to the kernel of an explicit map
\[
\xymatrix{
\delta_3 : B_3^*(E) \ar[r] & \bigl(\kb^\times \otimes E(\kb) \bigr)^{\Gal(\kb/k)}.
}
\]
\par\smallskip\noindent\emph{Proof of Theorem \ref{GL thm beta}.}
By Quillen's localisation theorem, there is a long exact sequence
\[
\xymatrix{
\cdots  \ar[r] & \bigoplus_{p \in E} K_2(k(p))  \ar[r] & K_2(E)  \ar[r] & K_2(k(E))  
\ar[r]^{\partial} & \bigoplus_{p \in E} k(p)^\times  \ar[r] & \cdots,
}\]
where $p$ runs over the closed points of $E$; see the exact sequence in the proof of \cite[Theorem 5.4]{Qui73}, with $i=2$ and $p=0$. Tensoring with $\Q$ and using the fact that $K_2$ of a number field is a torsion group \cite{garland}, we get an isomorphism $K_2(E) \otimes \Q \cong \ker(\partial) \otimes \Q$.

By \cite[Theorem 3.8]{goncharov-levin}, the natural map $K_2(k(E)) \to B_3^*(E)$ induces an isomorphism
\begin{equation*}
\biggl(\frac{H^0(E,\mathcal{K}_2)}{\operatorname{Tor}(k^\times,E(k)) + K_2(k)}\biggr) \otimes \Q \cong \ker(\delta_3) \otimes \Q,
\end{equation*}
where $H^0(E,\mathcal{K}_2)=\ker(\partial)$. Since $\operatorname{Tor}(k^\times,E(k))$ is a torsion group \cite[Proposition 3.1.2(a)]{weibel} and $K_2(k)$ is also torsion, we get the result.
\hfill\qed\par\smallskip

We will need to work with the full group of divisors $\Z[E(\kb)]$, using (a modification of) the group $B_3(E)$ introduced in \cite[Definition 3.1]{goncharov-levin}. The difference is that we use the $m$-distribution relations only for $m=-1$. 

\begin{defn} Let $R_3(E/\kb)$ be the subgroup of $\Z[E(\kb)]$ generated by the divisors $\beta(f, 1-f)$ with $f \in \kb(E)$, $f \neq 0,1$ as well as the divisors $(p)+(-p)$ with $p \in E(\kb)$. We define
\begin{equation*}
B_3(E/\kb) = \frac{\Z[E(\kb)]}{R_3(E/\kb)}, \qquad \textrm{and} \qquad B_3(E) = B_3(E/\kb)^{\Gal(\kb/k)}.
\end{equation*}
\end{defn}

Goncharov and Levin proved the following result (compare \cite[Proposition 3.19(a)]{goncharov-levin}).

\begin{prop}\label{pro B3 B3*}
The canonical map $B_3^*(E) \otimes \Q \to B_3(E) \otimes \Q$ is injective.
\end{prop}

\begin{pf}
It suffices to establish the result for $E/\kb$. Let $D = \sum n_i [p_i] \in I^4$ be a divisor belonging to $R_3(E/\kb)$. Write $D=D' + D''$ with $D' \in R_3^*(E/\kb)$ and $D'' = \sum_j m_j ((q_j)+(-q_j))$. The divisor $D''$ belongs to $I^4$ and is invariant under the elliptic involution $\sigma: p\mapsto -p$ on $E$. Thus we can write
\begin{equation*}
2D'' = D'' + \sigma(D'') = \beta \Bigl( \sum_\ell (f_\ell \otimes g_\ell) + \sigma^* (f_\ell \otimes g_{\ell})\Bigr)
\end{equation*}
for some rational functions $f_\ell$, $g_\ell$. By \cite[Lemma 3.21]{goncharov-levin}, for any rational functions $f,g$ on $E/\kb$, we have $\sigma^* \{f,g\} = - \{f,g\}$ in $K_2(\kb(E))/\{\kb^\times,\kb(E)^\times\}$. It follows that $(f \otimes g) + \sigma^*(f \otimes g)$ is a linear combination of Steinberg relations and of elements $\lambda \otimes h$ with $\lambda \in \kb^\times$ and $h \in \kb(E)^\times$. Applying the map $\beta$ and noting that $\beta(\lambda \otimes h)=0$, we get $2D'' \in R_3^*(E/\kb)$ as desired.
\end{pf}


Putting together Theorem \ref{GL thm beta} and Proposition \ref{pro B3 B3*}, we get the following result.

\begin{thm}\label{GL thm beta 2}
The composite map
\[
\xymatrix{
\overline{\beta} : \quad K_2(E) \otimes \Q \ar[r] 
& K_2(k(E)) \otimes \Q \ar[r] & B_3(E) \otimes \Q,
}\]
sending an element $\sum_i \{f_i, g_i\}$ to the class of the divisor $\sum_i \beta(f_i,g_i)$, is injective.
\end{thm}

Thanks to Theorem \ref{GL thm beta 2}, any equality in $K_2(E) \otimes \Q$ can be proved by applying the map $\overline{\beta}$ and comparing the divisors. Of course, the difficult part is to find the necessary Steinberg relations. In the following sections, we use this strategy to compare an ``elliptic'' and a ``modular'' element in $K_2$ of the elliptic curve $X_0(24)$.

\subsection{Special elements in $K_2$ of $X_0(24)$}

\subsubsection{The minimal model}

The curve $E_4$ in the Legendre family is given by the equation $y^2 = x(1-x)(1+3x)$. A minimal Weierstrass equation is
\begin{equation*}
E : Y^2=X^3-X^2-4X+4,
\end{equation*}
obtained with the change of variables $(X,Y)=(1-3x, -3y)$.
This is the curve $24a1$ in the Cremona database \cite{cremona}. The N\'eron 
differential is (up to sign)
\begin{equation*}
\omega_E = \frac{dX}{2Y} = \frac{dx}{2y}.
\end{equation*}
The torsion subgroup of $E$ is isomorphic to $\Z/4\Z \times \Z/2\Z$, generated by the points $p_1 = (0,2)$ of order 4, and $p_2 = (1,0)$ of order 2.

\subsubsection{The modular parametrisation}

The curve $E$ is in fact isomorphic to the modular curve $X_0(24)$. We denote by
\begin{equation*}
\varphi : X_0(24) \to E
\end{equation*}
the modular parametrisation, normalised so that $\varphi(\infty)=0$ and $\varphi^*(\omega_E) = \omega_f = 2\pi if(z)dz$, where $f$ is the unique newform of weight 2 and level $\Gamma_0(24)$.

The modular curve $X_0(24)$ has 8 cusps: $\infty, 0, \frac12, \frac13, \frac14, \frac16, \frac18, \frac{1}{12}$. They are all rational, and therefore correspond via $\varphi$ to the 8 rationals points on $E$. We now make explicit this bijection. Let $\Lambda \subset \C$ be the lattice of periods of $\omega_E = dX/(2Y)$. We have a canonical isomorphism
\begin{equation*}
\gamma : \C/\Lambda \xrightarrow{\cong} E(\C)
\end{equation*}
such that $\gamma^*(\omega_E) = dz$ (so that $\gamma^{-1}(p) = \int_0^p \omega_E$). The idea is the following: given a point $\tau \in \h \cup \p^1(\Q)$, we have
\begin{equation*}
\varphi(\tau) = \gamma \left(\int_0^{\varphi(\tau)} \omega_E \right) = \gamma \left( \int_\infty^{\tau} \omega_f \right).
\end{equation*}
The last integral, as well as the map $\gamma$, can be computed using \textsc{Pari/GP} \cite{PARI2}. Note that although the computation is numerical, we know that if $\tau$ is a cusp, then $\int_\infty^\tau \omega_f$ belongs to the lattice $\frac14 \Lambda$, hence its value can be ascertained.

We used the following \textsc{Pari/GP} code to compute the images of the cusps of $X_0(24)$ by $\varphi$.

{\footnotesize
\begin{verbatim}
E = ellinit("24a1");
mf = mfinit([24, 2]);
f = mfeigenbasis(mf)[1];
symb = mfsymbol(mf, f);
phiE(c) = ellztopoint(E, polcoef(mfsymboleval(symb, [oo, c])*2*Pi*I, 0));
apply(phiE, [oo, 0, 1/2, 1/3, 1/4, 1/6, 1/8, 1/12])
\end{verbatim}
}
The results are shown in the following table.

\begin{equation*}
\begin{tabular}{c|c|c|c|c|c|c|c|c}
$c$ & $\infty$ & $0$ & $\frac12$ & $\frac13$ & $\frac14$ & $\frac16$ & $\frac18$ & $\frac{1}{12}$ \\[5pt]
\hline
& & & & & & & & \\[-10pt]
$\varphi(c)$ & $0$ & $3p_1+p_2$ & $p_1+p_2$ & $3p_1$ & $2p_1+p_2$ & $p_1$ & $p_2$ & $2p_1$ \\[5pt]
\hline
& & & & & & & & \\[-10pt]
$(X,Y)$ & $(0:1:0)$ & $(4,-6)$ & $(4,6)$ & $(0,-2)$ & $(-2,0)$ & $(0,2)$ & $(1,0)$ & $(2,0)$
\end{tabular}
\end{equation*}

The sign of the functional equation of the $L$-function $L(E,s)$ is $+1$, hence the Atkin-Lehner involution $W_{24} : \tau \mapsto -1/(24\tau)$ satisfies $W_{24}(f) = -f$. This implies that the map $W_{24} : E \to E$ has the form $p \mapsto p_0-p$ for some rational point $p_0$ on $E$, and the table above gives $p_0=3p_1+p_2$.

\subsubsection{The Beilinson--Kato element}

Recall the definition of the Beilinson--Kato element $z_E$ in $K_2(E) \otimes \Q$ (see \cite[D\'efinition 9.3]{Br}):
\begin{equation*}
z_E = \varphi_* \Bigl(\frac12 \{u_{24},W_{24} (u_{24})\}'\Bigr),
\end{equation*}
where $u_N$ is the following product of Siegel units
\begin{equation*}
u_N = \prod_{\substack{a \in (\Z/N\Z)^\times \\ b \in \Z/N\Z}} g_{a,b},
\end{equation*}
and the superscript $'$ means addition of Milnor symbols $\{\lambda,v\}$ with $\lambda \in \Q^\times$ and $v \in \mathcal{O}(Y_0(24))^\times$ in order to obtain an element of $K_2(X_0(24)) \otimes \Q$. Since the symbols $\{\lambda,v\}$ are killed by $\beta$, we can safely ignore them in the computation.

We wish to compute the divisor of the modular unit $u_{24}$. Let us work more generally with $u_N$. From the definition of Siegel units as infinite products, we know that the order of vanishing of $g_{a,b}$ at the cusp $\infty$ of $X(N)$ is equal to $N B_2(\{\frac{a}{N}\})/2$, where $B_2(x)=x^2-x+1/6$ is the second Bernoulli polynomial and $\{t\} = t-\lfloor t \rfloor$ is the fractional part of $t$. Moreover, we have the transformation formula $g_{a,b} \circ \alpha = g_{(a,b)\alpha}$ in $\mathcal{O}(Y(N))^\times \otimes \Q$ for any $\alpha \in \SL_2(\Z)$. Using this, we can compute the order of vanishing of $g_{a,b}$ at any cusp.

Since we are working with $X_0(N)$ instead of $X(N)$, we need to take into account the widths of the cusps of $X_0(N)$.
The width of the cusp $1/d \in X_0(N)$ is
\begin{equation*}
w(1/d) = \frac{N}{d \cdot \gcd(d,N/d)}.
\end{equation*}
Since $1/d = (\begin{smallmatrix} 1 & 0 \\ d & 1 \end{smallmatrix})\infty$, we have
\begin{align*}
\ord_{1/d}(u_N) & = \sum_{\substack{a \in (\Z/N\Z)^\times \\ b \in \Z/N\Z}} \ord_{1/d}(g_{a,b}) \\
& = \sum_{\substack{a \in (\Z/N\Z)^\times \\ b \in \Z/N\Z}} w(1/d) \ord_\infty(g_{a+db,b}) \\
& = \frac{w(1/d)}{2} \sum_{\substack{a \in (\Z/N\Z)^\times \\ b \in \Z/N\Z}} B_2\bigl(\bigl\{\frac{a+db}{N}\bigr\}\bigr) \\
& = \frac{d \varphi(N)}{2 \gcd(d,N/d) \varphi(d)} \sum_{a \in (\Z/d\Z)^\times} B_2\bigl(\bigl\{\frac{a}{d}\bigr\}\bigr).
\end{align*}
Here we used the distribution relation for the periodic Bernoulli polynomials,
\begin{equation*}
B_n(\{mt\}) = m^{n-1} \sum_{k=0}^{m-1} B_n \bigl(\bigl\{t+\frac{k}{m}\bigr\}\bigr) \qquad (m \geq 1).
\end{equation*}
We deduce the order of vanishing of $u_{24} $ at each cusp of $X_0(24)$, and therefore its divisor:
\begin{equation*}
 \mathrm{div}(u_{24}) = \frac16 (\infty) + \frac23 (0) - \frac13 (1/2) - \frac23 (1/3) - \frac16 (1/4) + \frac13 (1/6) - \frac16 (1/8) + \frac16 (1/12).
\end{equation*}
The fractions appearing here mean that $u_{24}$ is only an element of $\mathcal{O}(Y_0(24))^\times \otimes \Q$, in other words some power of $u_{24}$ is a modular unit.

Applying the modular parametrisation $\varphi$, we get
\begin{equation} \label{div u24}
 \mathrm{div}(u_{24}) = \frac16 \bigl( (0) +4 (3p_1+p_2) - 2(p_1+p_2) - 4 (3p_1) + 2(p_1) - (2p_1+p_2) - (p_2) + (2p_1) \bigr).
\end{equation}
Applying the Atkin-Lehner involution, we also have
\begin{equation} \label{div W24 u24}
 \mathrm{div}(W_{24} (u_{24})) = \frac16 \bigl( (3p_1+p_2) +4 (0) - 2(2p_1) - 4 (p_2) + 2(2p_1+p_2) - (p_1) - (3p_1) + (p_1+p_2) \bigr).
\end{equation}

\subsubsection{The symbol $\xi$ in \eqref{RZ-eq2}}

Recall 
\begin{equation}
\xi = \{f,g\} \in K_2(E_4) \otimes \Q, \qquad f = \frac{y-x+1}{y+x-1}, \qquad g = -\frac{(x-1)^2}{4x^2}.
\end{equation}
Using \textsc{Magma},
we can find the divisors of $f$ and $g$. Here is the code:
{\footnotesize
\begin{verbatim}
A2<x,y> := AffinePlane(Rationals());
C := Curve(A2, y^2-x*(1-x)*(1+3*x)); 
Cbar := ProjectiveClosure(C);
E, phi := EllipticCurve(Cbar);
Emin, psi := MinimalModel(E);
F<x,y> := FunctionField(Cbar);
f := (y-x+1)/(y+x-1);
g := -(x-1)^2/(4*x^2);
div_f := Decomposition(Divisor(f));
div_g := Decomposition(Divisor(g));
print "div(f) =", [<p[2], psi(phi(RepresentativePoint(p[1])))> : p in div_f];
print "div(g) =", [<p[2], psi(phi(RepresentativePoint(p[1])))> : p in div_g];
\end{verbatim}
}
We obtain
\begin{align} \label{div fg}
 \mathrm{div}(f) & = -(3p_1) + (p_1) - (3p_1+p_2) + (p_1+p_2) \\
\nonumber  \mathrm{div}(g) & = 4(2p_1+p_2) - 4(p_2).
\end{align}

\subsection{Comparing the divisors} \label{apply beta}

We are now going to apply $\overline{\beta}$ to the two elements in $K_2(E) \otimes \Q$, and compare the results.

For the Beilinson--Kato element, we find using \eqref{div u24} and \eqref{div W24 u24} that
\begin{align*}
\beta(u_{24}, W_{24}(u_{24})) & = \frac{1}{36} \bigl( 8(0) - 8(p_2) + 28 (p_1) - 28 (p_1+p_2) \\
 & \qquad \quad + 8 (2p_1) - 8 (2p_1+p_2) - 44(3p_1) + 44(3p_1+p_2)\bigr).
\end{align*}

In the group $B_3(E) \otimes \Q$, we have the relation $(p)+(-p)=0$ for any point $p$, hence $(p)=0$ if $p$ is $2$-torsion. So we can remove the 2-torsion points from the above divisor. In fact, we can express everything in terms of $p_1$ and $p_1+p_2$ alone. We find
\begin{equation*}
\overline{\beta}\bigl(\{u_{24} , W_{24}(u_{24})\}\bigr) = 2 (p_1) - 2 (p_1+p_2),
\end{equation*}
and thus
\begin{equation}\label{beta zE}
\overline{\beta}(z_E) = (p_1) - (p_1+p_2).
\end{equation}

We proceed similarly for the element $\xi$ in \eqref{RZ-eq2}. Using \eqref{div fg}, we compute
\begin{equation*}
\beta(f, g) = -8(p_1) -8(p_1+p_2) + 8 (3p_1) + 8(3p_1+p_2),
\end{equation*}
which gives
\begin{equation}\label{beta xi}
\overline{\beta}(\xi) = - 16(p_1) -16(p_1+p_2).
\end{equation}

The divisors $\overline{\beta}(z_E)$ and $\overline{\beta}(\xi)$ are not proportional, which suggests that there should be a non-trivial Steinberg relation involving $p_1$ and $p_1+p_2$. We can determine it experimentally by computing the elliptic dilogarithm of these points. Let us denote by $D_E : E(\C) \to \R$ the Bloch elliptic dilogarithm. Using \textsc{Pari/GP}, we find numerically
\begin{equation}
5 D_E(p_1) + 3 D_E(p_1+p_2) \approx 0.
\end{equation}
This means that we should have $5(p_1)+ 3(p_1+p_2) = 0$ in the group $B_3(E) \otimes \Q$. We will prove that this is indeed the case, by exhibiting a Steinberg relation.

We search for a rational function $h$ on $E$ such that the zeros and poles of both $h$ and $1-h$ are among the 8 torsion points of $E$. To do this, we use Mellit's technique of \emph{incident lines} \cite{mellit}; see also \cite[Proof of Lemma 3.29]{goncharov-levin}.

We view $E$ as a non-singular plane cubic. We generate all the lines passing only through the 8 torsion points of $E$. Say we have found three such lines $\ell_1, \ell_2, \ell_3$ which, moreover, meet at a point $p_0$ of $\mathbb{P}^2$. If the lines are pairwise distinct, we may choose equations $f_1,f_2,f_3$ for these lines satisfying $f_1+f_2=f_3$. Then $h=f_1/f_3$ has the property that the divisors of $h$ and $1-h$ are supported at the torsion points. In particular $\beta(h, 1-h)$ is also supported at the torsion points, which gives a relation in $B_3(E) \otimes \Q$.

If the intersection point $p_0$ lies on the curve, then the above relation is trivial: it is contained in the subgroup generated by the divisors $(p)+(-p)$. If, however, $p_0$ does not lie on the curve, then we usually get something interesting. It turns out that this method of incident lines works remarkably well in practice.

Using a computer, it is possible to search for all incident lines, and determine the associated Steinberg relations. In the present situation, we find the lines $\ell_1$, $\ell_2$ defined by the equations
\begin{equation*}
f_1 = -\frac14 (X+Y-2) \qquad f_2 = \frac14 (X+Y+2).
\end{equation*}
We have $f_1+f_2=1$, so that the lines are parallel (taking $\ell_3$ to be the line at infinity, the lines $\ell_1,\ell_2,\ell_3$ are incident, so this is a particular case of the situation above). The divisors of these functions are given by
\begin{align*}
 \mathrm{div}(f_1) & = 2(p_1) + (2p_1) - 3(0) \\
 \mathrm{div}(f_2) & = (3p_1) + (2p_1+p_2) + (3p_1+p_2) - 3(0)
\end{align*}
and the associated Steinberg relation is
\begin{align*}
\beta(f_1, f_2) & = 9(0) + (p_2) - 9(p_1) - 3(p_1+p_2) - (2p_1) - (2p_1+p_2) + (3p_1) + 3(3p_1+p_2)\\
& \equiv -10 (p_1) - 6(p_1+p_2) \qquad \textrm{in } B_3(E) \otimes \Q.
\end{align*}
This shows that indeed $5(p_1)+3(p_1+p_2)=0$ in $B_3(E) \otimes \Q$. Thus \eqref{beta zE} and \eqref{beta xi} simplify:
\begin{equation*}
\overline{\beta}(z_E) = (p_1) +\frac53 (p_1) = \frac83 (p_1)
\end{equation*}
and
\begin{equation*}
\overline{\beta}(\xi) = -16 (p_1) - 16 \times -\frac53 (p_1) = \frac{32}{3} (p_1).
\end{equation*}
Using Theorem \ref{GL thm beta 2}, we deduce 
\begin{thm}\label{compare-thm}
$\xi = 4z_E$ in $K_2(E) \otimes \Q$.
\end{thm}

\noindent
Hokkaido University, Department of Mathematics, 
\par\noindent
Sapporo 060-0810, Japan.
\smallskip

\noindent
e-mail: \texttt{asakura@math.sci.hokudai.ac.jp}

\bigskip

\noindent
Tokyo Denki University, School of Science and Technology for Future Life,
\par\noindent
Tokyo 120-8551, Japan.

\smallskip

\noindent
e-mail: \texttt{chida@mail.dendai.ac.jp}

\bigskip

\noindent
\'ENS Lyon, Unit\'e de math\'ematiques pures et appliqu\'ees, 
\par\noindent
46 all\'ee d'Italie, 69007 Lyon, France.

\smallskip

\noindent
e-mail: \texttt{francois.brunault@ens-lyon.fr}

\end{document}